\documentclass[a4paper,reqno,10pt]{amsart}
\usepackage[utf8]{inputenc}
\usepackage{amssymb,amsmath,dsfont,tikz-cd}
\usepackage{mathtools}
\usepackage{hyperref}
\usepackage[all]{xy}
\usepackage{upgreek}
\usepackage{colonequals}
\usepackage{todonotes}
\usepackage{amsthm}
\usepackage{mathrsfs}
\usepackage{stmaryrd}
\usepackage{bm}

\newtheorem{innercustomthm}{Conjecture}
\newenvironment{customconj}[1]
  {\renewcommand\theinnercustomthm{#1}\innercustomthm}
  {\endinnercustomthm}

\usepackage{etoolbox}

\usepackage{relsize}

\renewcommand{\theta}{\uptheta}
\renewcommand{\alpha}{\upalpha}
\renewcommand{\beta}{\upbeta}
\renewcommand{\gamma}{\upgamma}
\renewcommand{\delta}{\updelta}
\renewcommand{\zeta}{\upzeta}
\renewcommand{\pi}{\uppi\hspace{0.05em}}
\renewcommand{\rho}{\uprho}
\renewcommand{\xi}{\upxi}
\renewcommand{\chi}{\upchi}
\renewcommand{\sigma}{\upsigma}
\renewcommand{\Lambda}{\Uplambda}
\renewcommand{\Gamma}{\Upgamma}
\renewcommand{\phi}{\upphi}
\renewcommand{\psi}{\uppsi}
\renewcommand{\nu}{\upnu}
\renewcommand{\tau}{\uptau}
\renewcommand{\mu}{\upmu}
\renewcommand{\eta}{\upeta}

\newtheorem{theorem}{Theorem}[section]
\newtheorem{thmx}{Theorem}

\newtheorem{proposition}[theorem]{Proposition}

\newtheorem{lemma}[theorem]{Lemma}
\newtheorem{conjecture}[theorem]{Conjecture}
\newtheorem{pconjx}{``Conjecture''}

\newtheorem{conjx}{Conjecture}

\newtheorem{corollary}[theorem]{Corollary}

\theoremstyle{remark}
\newtheorem{remark}[theorem]{Remark}

\theoremstyle{definition}
\newtheorem{definition}[theorem]{Definition}

\newcommand{\BA}{{\mathbb{A}}}

\newcommand{\BC}{{\mathbb{C}}}

\newcommand{\BD}{{\mathbb{D}}}

\DeclareMathOperator{\HO}{\mathbf{H}}

\DeclareMathOperator{\BM}{BM}
\DeclareMathOperator{\Hit}{\mathtt{h}}

\newcommand{\kac}{\mathtt{a}}

\newcommand{\BN}{{\mathbb{N}}}

\newcommand{\BP}{{\mathbb{P}}}
\newcommand{\BQ}{{\mathbb{Q}}}

\newcommand{\BZ}{{\mathbb{Z}}}

\newcommand{\CF}{{\mathcal F}}
\newcommand{\CG}{{\mathcal G}}
\newcommand{\CH}{{\mathcal H}}

\newcommand{\CL}{{\mathcal L}}
\newcommand{\CM}{{\mathcal M}}
\newcommand{\CN}{{\mathcal N}}

\newcommand{\CP}{{\mathcal P}}

\newcommand{\CU}{{\mathcal U}}

\DeclareMathOperator{\vdim}{vdim}

\newcommand{\Fg}{{\mathfrak{g}}}
\newcommand{\Fh}{{\mathfrak{h}}}

\newcommand{\Fn}{{\mathfrak{n}}}

\newcommand{\Fp}{{\mathfrak{p}}}

\newcommand{\FF}{{\mathfrak{F}}}

\newcommand{\FH}{{\mathfrak{H}}}

\newcommand{\FJ}{{\mathfrak{J}}}

\newcommand{\FL}{{\mathfrak{L}}}
\newcommand{\FM}{{\mathfrak{M}}}
\newcommand{\FN}{{\mathfrak{N}}}

\let\ol=\overline
\let\ul=\underline

\DeclareMathOperator{\EXP}{EXP}

\DeclareMathOperator{\JH}{\mathtt{JH}}
\newcommand{\Up}{{}^{\mathfrak{p}\!}}

\newcommand{\ICS}{\mathcal{IC}}

\newcommand{\DTS}{\mathcal{BPS}}

\newcommand{\phim}[1]{\phi^{\mon}_{#1}}
\newcommand{\phin}[1]{\Up\phi_{#1}}

\newcommand{\sstab}{\textrm{sst}}

\DeclareMathOperator{\crit}{crit}

\DeclareMathOperator{\Free}{Free}
\DeclareMathOperator{\Bor}{\mathbf{Bor}}
\DeclareMathOperator{\Lie}{Lie}

\DeclareMathOperator{\Hom}{Hom}

\DeclareMathOperator{\mon}{mon}

\DeclareMathOperator{\Ext}{Ext}

\DeclareMathOperator{\MMHM}{MMHM}

\DeclareMathOperator{\Ker}{Ker}

\DeclareMathOperator{\Gr}{Gr}
\DeclareMathOperator{\FoM}{FM}
\DeclareMathOperator{\FFM}{\mathfrak{F}M}

\DeclareMathOperator{\Lvert}{\mathlarger{\mathlarger{\mathlarger{\lvert}}}}

\DeclareMathOperator{\Id}{Id}

\DeclareMathOperator{\dd}{\mathbf{d}}
\DeclareMathOperator{\ee}{\mathbf{e}}

\DeclareMathOperator{\EP}{E}

\DeclareMathOperator{\Perv}{Perv}

\DeclareMathOperator{\MHM}{MHM}

\DeclareMathOperator{\ICA}{IC}
\DeclareMathOperator{\DT}{BPS}
\DeclareMathOperator{\Sym}{Sym}

\DeclareMathOperator{\Spec}{Spec}
\DeclareMathOperator{\Gl}{GL}
\DeclareMathOperator{\id}{id}

\DeclareMathOperator{\Jac}{Jac}

\DeclareMathOperator{\Tr}{Tr}

\DeclareMathOperator{\Pic}{Pic}
\DeclareMathOperator{\pt}{pt}

\DeclareMathOperator{\vir}{vir}
\DeclareMathOperator{\Ob}{Ob}

\DeclareMathOperator{\LLL}{\mathscr{T}}
\DeclareMathOperator{\Ho}{\mathcal{H}}

\DeclareMathOperator{\B}{B\!}

\DeclareMathOperator{\Dol}{\scriptscriptstyle{Dol}}

\DeclareMathOperator{\Coha}{\mathcal{A}}
\DeclareMathOperator{\RCoha}{\widetilde{\mathcal{A}}}

\newcommand{\Db}{\mathcal{D}^{b}}
\newcommand{\Dub}{\mathcal{D}}

\DeclareMathOperator{\Coh}{\mathbf{Coh}}

\newcommand{\csdot}{\!\cdot\!}

\newcommand{\UP}{{}^{\Fp}\!}

\DeclareMathOperator{\Betti}{\scriptscriptstyle{Betti}}
\let \ol=\overline
\let \ul=\underline

\DeclareMathOperator{\Mat}{Mat}

\DeclareMathOperator{\simp}{simp}

\DeclareMathOperator{\Sta}{st}

\DeclareMathOperator{\ana}{an}
\newcommand{\Msp}{\mathcal{M}}
\newcommand{\Mst}{\mathfrak{M}}

\DeclareMathOperator{\supp}{supp}

\DeclareMathOperator{\UEA}{\mathbf{U}}

\title{Nonabelian Hodge theory for stacks and a stacky P=W conjecture}
\author{Ben Davison}

\begin{document}

\setcounter{tocdepth}{2}

\let\oldtocsection=\tocsection

\let\oldtocsubsection=\tocsubsection

\let\oldtocsubsubsection=\tocsubsubsection

\renewcommand{\tocsection}[2]{\hspace{0em}\oldtocsection{#1}{#2}}
\renewcommand{\tocsubsection}[2]{\hspace{3em}\oldtocsubsection{#1}{#2}}
\renewcommand{\tocsubsubsection}[2]{\hspace{2em}\oldtocsubsubsection{#1}{#2}}

\makeatletter
\patchcmd{\@tocline}
  {\hfil}
  {\leaders\hbox{\,.\,}\hfil}
  {}{}
\makeatother

\maketitle

\begin{abstract}
We introduce a version of the P=W conjecture relating the Borel--Moore homology of the stack of representations of the fundamental group of a genus $g$ Riemann surface with the Borel--Moore homology of the stack of degree zero semistable Higgs bundles on a smooth projective complex curve of genus $g$.  In order to state the conjecture we propose a construction of a canonical isomorphism between these Borel--Moore homology groups.  We relate the stacky P=W conjecture with the original P=W conjecture concerning the cohomology of smooth moduli spaces, and the PI=WI conjecture concerning the intersection cohomology groups of singular moduli spaces.  In genus zero and one, we prove the conjectures that we introduce in this paper.

\end{abstract}

\tableofcontents

\section{Introduction}
\label{sec:introduction}
In this paper we explain a framework for extending nonabelian Hodge theory and the P=W conjecture to the Borel--Moore homology\footnote{Equivalently, dual compactly supported cohomology.} of the singular Betti and Dolbeault stacks associated to a smooth projective complex curve.  We use the decomposition theorem for stacks of objects in 2-Calabi--Yau (2CY) categories, along with Hall algebra constructions, to begin the extension of the nonabelian Hodge correspondence to the Borel--Moore homology of moduli stacks, and formulate (and prove for genus less than two) a stacky P=W conjecture (the ``PS=WS conjecture'').
\smallbreak
Furthermore, we explain how (conjecturally) the famous P=W conjecture\footnote{Now a theorem, proved independently in \cite{PW1,PW2}} due to de Cataldo, Hausel and Migliorini (\cite{deCat12}) comparing the weight filtration on the cohomology of the twisted character variety $\Msp^{\Betti}_{g,r,d}$ (with $(r,d)=1$) with the perverse filtration on the cohomology of $\Msp^{\Dol}_{r,d}(C)$ for $C$ a genus $g$ smooth projective curve, is equivalent to two other P=W conjectures.  One of these (the PI=WI conjecture due to de Cataldo and Maulik \cite{dCM20}) regards intersection cohomology.  We refer to \cite{FM20,Ma21} for further discussion and results on the PI=WI conjecture.  The other -- the PS=WS conjecture introduced here -- compares the weight filtration on the Borel--Moore homology of $\Mst^{\Betti}_{g,r,0}$ (the stack of $r$-dimensional representations of $\pi_1(\Sigma_g)$) with a perverse filtration on the Borel--Moore homology of $\Mst^{\Dol}_{r,0}(C)$ (the stack of rank $r$ degree zero semistable Higgs bundles on $C$).  We remark at the outset that stating a reasonable PS=WS conjecture requires some care; indeed, we will see that the most naive formulation is false.  The key issue is that the perverse filtration on the Dolbeault side has to be taken to be a mixture of two different perverse filtrations, both of which are defined with the help of the decomposition theorem for moduli stacks of objects in 2CY categories \cite{Da21a}; see \S \ref{MixPerv} for a more specific formulation.

\subsection{Nonabelian Hodge theory for moduli schemes}
\label{NAHTschemes}
Before coming to stacks, we briefly recall the more familiar version of the theory, involving smooth moduli schemes.  Let $C$ be a smooth projective complex curve of genus $g$.  A \textit{Higgs bundle} on $C$ is a pair
\[
(\CF,\eta\colon \CF\rightarrow \CF\otimes\omega_C)
\]
where $\CF$ is a locally free coherent sheaf on $C$.  We define the rank and degree of $(\CF,\eta)$ to be the rank and degree of the underlying locally free sheaf $\CF$.  A Higgs bundle is called \textit{semistable} if for all coherent sheaves $\CG\subset \CF$ satisfying $\eta(\CG)\subset \CG\otimes\omega_C$ the slope of $\CG$ is no greater than that of $\CF$.  We denote by $\Msp^{\Dol}_{r,d}(C)$ the moduli space of rank $r$, degree $d$ semistable Higgs bundles.  This moduli space is a smooth variety if $(r,d)=1$, and is a singular variety if the genus of $C$ is nonzero and $(r,d)\neq 1$.
\smallbreak
Let $C'$ be obtained by removing a point $x$ from $C_{\ana}$, the underlying analytic space of $C$.  We denote by $\Msp^{\Betti}_{g,r,d}$ the smooth moduli space of those $r$-dimensional representations of the fundamental group of $C'$ for which the monodromy around the puncture is given by $\exp(2d\pi \sqrt{-1}/r)\Id_{r\times r}$.  By the nonabelian Hodge correspondence \cite{Hi87,Co88,Do87,Si94} there is a homeomorphism
\begin{equation}
\label{ClNAHT}
\Gamma_{r,d}\colon \Msp^{\Betti}_{g,r,d}\rightarrow \Msp^{\Dol}_{r,d}(C).
\end{equation}
For now we consider $d=1$ so that $(r,d)=1$.  By Deligne's construction \cite{Deligne71, Deligne74}, the cohomology of the schemes $\Msp^{\Betti}_{g,r,1}$ and $\Msp^{\Dol}_{r,1}(C)$ carry mixed Hodge structures.  In particular, each of the cohomology groups
\[
\HO^n(\Msp^{\Betti}_{g,r,1},\BQ),\quad \quad \HO^n(\Msp^{\Dol}_{r,1}(C),\BQ)
\]
carry ascending \textit{weight filtrations}.  The morphism 
\begin{equation}
\label{UpsDef1}
\Upsilon\coloneqq \HO(\Gamma_{r,d})\colon \HO(\Msp^{\Betti}_{g,r,d},\BQ)\xrightarrow{\cong}\HO(\Msp^{\Dol}_{r,d}(C),\BQ)
\end{equation}
does not respect these mixed Hodge structures; the homeomorphism $\Gamma_{r,d}$ is not induced by a morphism in the category of algebraic varieties.  To describe the filtration on $\HO^n(\Msp^{\Dol}_{r,1}(C),\BQ)$ that corresponds to the weight filtration on $\HO^n(\Msp^{\Betti}_{g,r,1},\BQ)$ under $\Upsilon$ we next recall the construction of \textit{perverse filtrations}.
\smallbreak
The moduli space $\Msp^{\Dol}_{r,1}(C)$ comes equipped with the \textit{Hitchin map}
\begin{align}
\label{hmap}
\Hit_{r,1}\colon &\Msp^{\Dol}_{r,1}(C)\rightarrow \Lambda_r\coloneqq \prod_{i=1}^r\HO^0(C,\omega_C^{\otimes i})\\ \nonumber
&(\CF,\eta)\mapsto \mathrm{char}(\eta)= (\Tr(\eta),\Tr(\wedge^2\eta)),\ldots)
\end{align}
sending a Higgs bundle $(\CF,\eta)$ to the generalised eigenvalues of $\eta$.  We will drop the subscripts from $\Hit$ where there is no danger of confusion.  The Hitchin map is projective, and so since $\Msp^{\Dol}_{r,1}(C)$ is smooth, the Beilinson--Bernstein--Deligne--Gabber decomposition theorem \cite{BBD} implies that there is a direct sum decomposition into shifted perverse cohomology sheaves\footnote{Here, and throughout the paper, for $g$ a morphism of stacks or varieties, we denote by $g_*$ the \textit{derived} direct image along $g$.}
\begin{equation}
\label{ClaDec}
\Hit_{*}\!\BQ_{\Msp^{\Dol}_{r,1}(C)}\cong \bigoplus_{i\in\BZ}\Up\!\Ho^i\!\left(\Hit_{*}\!\BQ_{\Msp^{\Dol}_{r,1}(C)}\right)[-i].
\end{equation}
As a consequence of this decomposition, the natural morphism
\[
{}^{\Fp}\!\tau^{\leq i}\Hit_{*}\!\BQ_{\Msp^{\Dol}_{r,1}(C)}\rightarrow \Hit_{*}\!\BQ_{\Msp^{\Dol}_{r,1}(C)}
\]
has a left inverse.  Therefore, taking the derived direct image to a point, the morphism
\begin{equation}
\label{PIN}
\HO\!\left(\Lambda_r,{}^{\Fp}\!\tau^{\leq i}\Hit_{r,1,*}\!\BQ_{\Msp^{\Dol}_{r,1}(C)} \right)\rightarrow \HO\!\left(\Lambda_r,\Hit_{r,1,*}\!\BQ_{\Msp^{\Dol}_{r,1}(C)} \right)=\HO(\Msp^{\Dol}_{r,1}(C),\BQ)
\end{equation}
is an injection of cohomologically graded vector spaces.  The images of the inclusions \eqref{PIN} provide an ascending \textit{perverse filtration} defined by setting
\begin{equation}
\label{PFdef}
\FH^i\HO(\Msp^{\Dol}_{r,1}(C),\BQ)=\HO\!\left(\Lambda_r,{}^{\Fp}\!\tau^{\leq i+\dim(\Lambda_r)}\Hit_*\!\BQ_{\Msp^{\Dol}_{r,1}(C)}\right).
\end{equation}
The following theorem, relating the perverse and the weight filtrations, has recently been proved independently by Maulik and Shen, and by Hausel, Mellit, Minets and Schiffmann:
\begin{theorem}[\cite{PW1,PW2, MSY23b} (P=W conjecture \cite{deCat12})]
\label{PW}
With $\Upsilon$ defined as in \eqref{UpsDef1},
\[
\Upsilon\!\left(W_{2i}\HO^{n}(\Msp^{\Betti}_{g,r,1},\BQ)\right)=\FH^i\HO^n(\Msp^{\Dol}_{r,1}(C),\BQ),
\]
where $W_{2i}\HO^{n}(\Msp^{\Betti}_{g,r,1},\BQ)$ is the $2i$th piece of the weight filtration of the natural mixed Hodge structure on $\HO^{n}(\Msp^{\Betti}_{g,r,1},\BQ)$.  
\end{theorem}
Stated roughly (i.e. in a way that blurs out the overall shift in perverse degrees introduced in \eqref{PFdef}) the theorem claims that the passage from the weight filtration on the Betti side to the perverse filtration on the Dolbeault side \textit{halves} degrees.
\smallbreak
For the rest of \S \ref{NAHTschemes} we assume $g\geq 2$.  If now we set $d=0$, then the moduli spaces $\Msp^{\Betti}_{g,r,d}$ and $\Msp^{\Dol}_{r,d}(C)$ are singular, but their underlying topological spaces (for the analytic topology) are still homeomorphic via $\Gamma_{r,0}$.  By \cite[Thm.11.1]{Si94} these varieties are irreducible.  For an irreducible variety $X$ we denote by $\ICS_X$ the intermediate extension of the constant perverse sheaf $\BQ_{X^{\mathrm{sm}}[\dim(X)]}$ from the smooth locus, and define the normalised intersection cohomology via
\[
\ICA(X)\coloneqq \HO(X,\ICS_X).
\]
If $X$ is an even-dimensional variety, there is a natural lift of $\ICS_X$ to a pure weight zero mixed Hodge module, on $X$ and we endow $\ICA(X)$ with the resulting natural mixed Hodge structure.  Mimicking the definition in the smooth situation, we set
\[
\mathfrak{h}^i\ICA(\Msp^{\Dol}_{r,0}(C))\coloneqq\HO(\Lambda_r,\Up\tau^{\leq i}\Hit_*\ICS_{\Msp^{\Dol}_{r,0}(C)}).
\]
Then $\Gamma_{r,0}$ induces an isomorphism 
\[
\Upsilon=\ICA(\Gamma_{r,0})\colon \ICA(\Msp^{\Betti}_{g,r,0})\xrightarrow{\cong}\ICA(\Msp^{\Dol}_{r,0}(C))
\]
in intersection cohomology.
\begin{conjecture}[PI=WI conjecture \cite{dCM20}]
Assume that $g\geq 2$.  Then
\label{IPIW}
\[
\Upsilon(W_{2i}\ICA(\Msp_{g,r,0}^{\Betti}))=\mathfrak{h}^i\ICA(\Msp^{\Dol}_{r,0}(C)).
\]
\end{conjecture}
A consequence of the conjectures in this paper will be that the claim that the PI=WI conjecture given here is equivalent to the original P=W conjecture (Theorem \ref{PW}).
\subsection{Nonabelian Hodge theory for moduli stacks}
\label{cnaht}
Now we turn to moduli stacks, the main focus of this paper, and the setting in which we prove new results.  
\smallbreak
Let $\Mst^{\Betti}_{g,r,0}$ denote the stack of $r$-dimensional representations of $\pi_1(C_{\ana})$, and let $\Mst^{\Dol}_{r,0}(C)$ denote the moduli stack of rank $r$ degree zero semistable Higgs bundles on $C$.  These stacks are singular if $r\geq 2$ and $g\geq 1$.  Moreover, the stabiliser groups of certain objects (e.g. direct sums of the trivial $\pi_1(C_{\ana})$-representation) are large, in contrast to the stacky versions of $\Msp^{\Betti}_{g,r,1}$ and $\Msp^{\Dol}_{r,1}(C)$, which are $\B \BC^*$-gerbes over their respective coarse moduli schemes.  As such, it is an open question whether there exists an isomorphism
\begin{equation}
\label{NAHTS}
\Mst^{\Betti}_{g,r,0}\overset{?}{\cong}\Mst^{\Dol}_{r,0}(C)
\end{equation}
in the category of (topological) stacks.  Rather than address this question directly, in this paper we attempt to construct a natural isomorphism
\[
\Upsilon\colon \HO^{\BM}(\Mst^{\Betti}_{g,r,0},\BQ_{\vir})\cong \HO^{\BM}(\Mst^{\Dol}_{r,0}(C),\BQ_{\vir})
\]
out of the ``classical'' homeomorphism \eqref{ClNAHT}, and ideas arising from the study of cohomological Donaldson--Thomas theory.
\smallbreak
Some words on the notation are due here.  Firstly, we define the \textit{virtual dimension} of the stacks $\Mst^{\Betti}_{g,r,d}$ and $\Mst^{\Dol}_{r,d}(C)$ to be $2r^2(g-1)$.  Given a stack $\FN$ for which we have defined a locally constant virtual dimension $\vdim$, we define $\BQ_{\FN,\vir}\lvert_{\FN'}=\BQ_{\FN'}[\vdim(\FN')]$ on each connected component $\FN'$.  Then
\[
\HO^{\BM}(\FN,\BQ_{\vir})\coloneqq t_*\BD\BQ_{\FN,\vir}
\]
where $t\colon \FN\rightarrow \pt$ is the morphism to a point, $t_*$ is the derived direct image functor, and $\BD$ is the Verdier duality functor.  If $\FN$ is smooth and $d$-dimensional, then 
\[
\BD\BQ_{\FN,\vir}\cong \BQ_{\FN}[2d-\vdim]
\]
and 
\begin{align*}
\HO^{\BM}_{-\bullet}(\FN,\BQ_{\vir})\cong& \HO^{\bullet+ 2d-\vdim}(\FN,\BQ).
\end{align*}
So in the original setting of nonabelian Hodge theory, the Borel--Moore homology recovers the usual singular cohomology (with a shift).  Once we turn our attention to singular stacks, the Borel--Moore homology can be entirely different from the singular cohomology.  
\smallbreak
We define the morphism
\[
\Hit^{\Sta}_{r,0}\colon \Mst^{\Dol}_{r,0}(C)\rightarrow \Lambda_r
\]
in the same way as the original Hitchin map \eqref{hmap}, as the morphism taking a Higgs bundle to the generalised eigenvalues of the Higgs field.  Again, we will omit the subscripts from $\Hit^{\Sta}$ if there is no risk of ambiguity.  Then we define\footnote{As compared to \eqref{PFdef}, the overall shift in the definition of the perverse filtration has disappeared.  The fact that the comparison between the perverse and the weight filtration no longer has an overall shift in this paper is a consequence of systematically twisting/shifting by virtual dimensions.}
\begin{equation}
\label{PSH}
\FH^i\HO^{\BM}(\Mst^{\Dol}_{r,0}(C),\BQ_{\vir})\coloneqq\HO\left(\Lambda_r,{}^{\Fp}\!\tau^{\leq i}\:\Hit^{\Sta}_*\BD\BQ_{\Mst^{\Dol}_{r,0}(C),\vir}\right).
\end{equation}
Again, taking the derived global sections functor applied to the natural transformation
\[
\Up\tau^{\leq i}\rightarrow \id
\]
yields the morphism
\begin{equation}
\label{SPF}
\FH^i\!\HO^{\BM}(\Mst^{\Dol}_{r,0}(C),\BQ_{\vir})\rightarrow \HO^{\BM}(\Mst^{\Dol}_{r,0}(C),\BQ_{\vir}).
\end{equation}
\smallbreak
Since $\Mst^{\Betti}_{g,r,0}$ is a global quotient stack, one may define a mixed Hodge structure on its Borel--Moore homology (following \cite{Deligne74}); we denote this mixed Hodge structure $\HO^{\BM}(\Mst^{\Betti}_{g,r,0},\ul{\BQ})$.  In particular, the Borel--Moore homology of the stacks that we consider carry natural weight filtrations.  For a stack $\Mst$ with virtual dimension $\vdim$, we define the mixed Hodge structure\footnote{For now the reader may assume that the virtual dimension is always even.}
\begin{equation}
\label{virHS}
\HO^{\BM}\!\left(\Mst,\underline{\BQ}_{\vir}\right)\coloneqq \bigoplus_{\Mst'\in\pi_0(\Mst)}\HO^{\BM}\!\left(\Mst',\underline{\BQ}\right)\otimes\LLL^{-\vdim(\Mst')/2}\\
\end{equation}
where $\LLL$ is the pure cohomologically graded mixed Hodge structure
\[
\LLL\coloneqq \HO_c(\BA^1,\underline{\BQ}).
\]
Note that the underlying cohomologically graded vector space of \eqref{virHS} is identified with $\HO^{\BM}(\Mst,\BQ_{\vir})$.  When we talk of e.g. the weight filtration on $\HO^{\BM}(\Mst,\BQ_{\vir})$ we mean the weight filtration coming from the mixed Hodge structure \eqref{virHS}.  
Alternatively, in the language of mixed Hodge modules, we can write
\[
\HO^{\BM}\!\left(\Mst,\underline{\BQ}_{\vir}\right)=\HO\!\left(\Mst,\BD\underline{\BQ}_{\Mst,\vir}\right)
\]
where $\underline{\BQ}_{\Mst,\vir}\lvert_{\Mst'}=\underline{\BQ}_{\Mst'}\otimes \LLL^{-\vdim(\Mst')/2}$, and $\underline{\BQ}_{\Mst'}$ is the constant mixed Hodge module complex, for which the underlying complex of perverse sheaves is the constant complex $\BQ_{\Mst'}$.  
\smallbreak
After these preliminaries, we may formulate the following naive attempt at a ``stacky'' P=W conjecture:
\begin{pconjx}
\label{fcoj}
There is a natural isomorphism
\begin{equation}
\label{UpsDef}
\Upsilon\colon \HO^{\BM}(\Mst^{\Betti}_{g,r,0},\BQ_{\vir})\rightarrow \HO^{\BM}(\Mst^{\Dol}_{r,0}(C),\BQ_{\vir})
\end{equation}
such that $\Upsilon(W_{2i}(\HO^{\BM}(\Mst^{\Betti}_{g,r,0},\BQ_{\vir})))=\FH^i\!\HO^{\BM}(\Mst^{\Dol}_{r,0}(C),\BQ_{\vir})$.
\end{pconjx}

There are numerous problems with proving or even making sense of the above conjecture:
\begin{enumerate}
\item
Since the morphism $\Hit^{\Sta}\colon \Mst^{\Dol}_{r,0}(C)\rightarrow \Lambda_r$ is no longer a projective morphism of varieties, and its domain is not smooth, it is not clear that the perverse filtration defined with respect to it is well-behaved.  I.e. we would like to show that the morphisms \eqref{SPF} are injective.
\item
We do not have an isomorphism \eqref{NAHTS} which we can apply $\HO^{\BM}(-)$ to in order to define the isomorphism \eqref{UpsDef}.
\end{enumerate}
Assuming that these problems can be resolved, a more disconcerting one arises:
\begin{enumerate}
\setcounter{enumi}{2}
\item
``Conjecture'' \ref{fcoj} is false!
\end{enumerate}
To illustrate the third problem, we set $g=0$, i.e. $C=\BP^1$.  Then there are isomorphisms of stacks
\[
\Mst^{\Dol}_{r,0}(C)\cong \pt/\Gl_r\cong \Mst^{\Betti}_{0,r,0}.
\]
So in this case there is no problem defining $\Upsilon$: we simply apply $\HO^{\BM}(-,\BQ)$ to the isomorphism of algebraic stacks $\Mst^{\Dol}_{r,0}(C)\cong\Mst^{\Betti}_{0,r}$.  Since $\vdim$ in this example is twice the actual dimension of the stacks, which are smooth, there are isomorphisms
\begin{equation}
\label{BMforH}
\BD\BQ_{\FM^{\Dol}_{r,0}(\BP^1),\vir}\cong \BQ_{\FM^{\Dol}_{r,0}(\BP^1)},\quad\quad \BD\underline{\BQ}_{\FM^{\Betti}_{0,r,0},\vir}\cong \underline{\BQ}_{\FM^{\Betti}_{0,r,0}}
\end{equation}
and thus there are isomorphisms
\[
\HO^{\BM}(\FM^{\Dol}_{r,0}(\BP^1),\BQ_{\vir})\cong\HO(\FM^{\Dol}_{r,0}(\BP^1),\BQ),\quad\quad \HO^{\BM}(\FM^{\Betti}_{0,r,0},\underline{\BQ}_{\vir})\cong\HO(\FM^{\Betti}_{0,r,0},\underline{\BQ}).
\]
Moreover the Hitchin base $\Lambda_r$ is a point, so that the perverse filtration defined with respect to the morphism $\Mst^{\Dol}_{r,0}(C)\rightarrow \Lambda_r$ makes sense, and is in fact just the filtration by cohomological degree.  On the other hand, the natural mixed Hodge structure on $\HO(\pt/\Gl_r,\BQ)$ is pure by \cite{Deligne74}, i.e. the $i$th cohomology is pure of weight $i$.  So we find that the weight filtration of $\HO(\Mst^{\Betti}_{0,r},\BQ)$ also agrees with the filtration by cohomological degree.  In conclusion, we find that
\[
\Upsilon(W_{i}(\HO^{\BM}(\Mst^{\Betti}_{g,r,0},\BQ_{\vir})))=\FH^i\!\HO^{\BM}(\Mst^{\Dol}_{r,0}(C),\BQ_{\vir})
\]
contradicting Conjecture \ref{fcoj}, since there is no halving of degrees when passing from the weight to the perverse filtration.

\subsection{The PS=WS conjectures}
\label{MixPerv}
Happily the resolution of problem (3) follows naturally from a careful consideration of problems (1) and (2), and this is what we carry out in this paper.  
\smallbreak
Firstly, problem (1) is resolved via the decomposition theorem for good moduli stacks of objects in 2CY categories, established in \cite{Da21a} (recalled as Theorem \ref{p2CY}).  As a corollary of this decomposition theorem, along with Saito's version of the decomposition theorem for direct images of pure complexes of mixed Hodge modules, there is a decomposition 
\[
\Hit^{\Sta}_*\BD\BQ_{\Mst^{\Dol}_{r,0}(C),\vir}\cong \bigoplus_{i\in\BZ}\Up\!\Ho^i\!\left(\Hit^{\Sta}_*\BD\BQ_{\Mst^{\Dol}_{r,0}(C),\vir}\right)[-i]
\]
analogous to \eqref{ClaDec}, for the direct image of the complex $\BD\BQ_{\Mst^{\Dol}_{r,0}(C),\vir}$ to the Hitchin base $\Lambda_r$.  It follows that the morphisms \eqref{SPF} are indeed injective.  
\smallbreak
We denote by
\[
\JH\colon \Mst^{\Dol}_{r,0}(C)\rightarrow \Msp^{\Dol}_{r,0}(C)
\]
the morphism to the coarse moduli space.  A second application of the decomposition theorem from \cite{Da21a} presents a way to resolve problem (3).  The theorem states that there is a decomposition into shifts of semisimple perverse sheaves
\begin{equation}
\label{DolDC}
\JH_*\BD\BQ_{\Mst^{\Dol}_{r,0}(C),\vir}\cong \bigoplus_{i\in 2\cdot\mathbb{Z}_{\geq 0}} \!\!\!\Up\!\Ho^i\!\left( \JH_*\BD\BQ_{\Mst^{\Dol}_{r,0}(C),\vir}\right)[-i]
\end{equation}
meaning that if we write
\[
\FL^i\!\HO(\Mst^{\Dol}_{r,0}(C),\BQ_{\vir})\coloneqq \HO\left(\Msp^{\Dol}_{r,0}(C),\Up\tau^{\leq i}\JH_*\BD\BQ_{\Mst^{\Dol}_{r,0}(C),\vir}\right)
\]
the natural morphism $\FL^i\!\HO(\Mst^{\Dol}_{r,0}(C),\BQ_{\vir})\rightarrow \HO(\Mst^{\Dol}_{r,0}(C),\BQ_{\vir})$ is an injection.  The resulting filtration is called the \textit{less perverse filtration}, and respects mixed Hodge structures.  We obtain in the same way the less perverse filtration on $\HO^{\BM}(\Mst^{\Betti}_{g,r,0},\underline{\BQ}_{\vir})$ by mixed Hodge structures.  This different perverse filtration provides the correction required to state a reasonable version of P=W for stacks:
\setcounter{conjx}{1}
\begin{customconj}{$\textrm{B}^0$}[$\FL^0$(stacky nonabelian Hodge isomorphism \& PS=WS)]
\label{conj_b}
There is a natural isomorphism
\begin{equation}
\FL^0\Upsilon\colon \FL^0\!\HO^{\BM}\!\left( \Mst^{\Betti}_{g,r,0},\BQ_{\vir}\right) \xrightarrow{\cong} \FL^0\!\HO^{\BM}\!\left(\Mst^{\Dol}_{r,0}(C),\BQ_{\vir}\right)
\end{equation}
such that
\begin{equation}
\label{Lfilt}
\FL^0\Upsilon\left(W_{2i} \FL^0\!\HO^{\BM}(\Mst^{\Betti}_{g,r,0},\BQ_{\vir})\right)=\FH^i\FL^0\!\HO^{\BM}\left(\Mst^{\Dol}_{r,0}(C),\BQ_{\vir}\right).
\end{equation}
\end{customconj}
Conjecture \ref{conj_b} only concerns a piece of the Borel--Moore homology of the stacks $\Mst^{\Dol}_{r,0}(C)$ and $\Mst^{\Betti}_{g,r,0}$: the zeroth piece of the less perverse filtrations.  We give a closely related formulation, concerning the associated graded objects $\Gr^{\FL}_{\bullet}\HO^{\BM}\!\left( \Mst^{\Betti}_{g,r,0},\BQ_{\vir}\right)$ and $\Gr^{\FL}_{\bullet}\HO^{\BM}\!\left(\Mst^{\Dol}_{r,0}(C),\BQ_{\vir}\right)$ of the full Borel--Moore homology with respect to the less perverse filtration, for which we consider the filtrations induced by the weight filtration and the perverse filtration with respect to the morphism $\Hit^{\Sta}$, respectively.

\begin{conjx}[stacky nonabelian Hodge isomorphism \& PS=WS]
\label{conj_c}
There is a natural isomorphism
\begin{equation}
\label{SNAHT}
\Upsilon\colon \HO^{\BM}\!\left( \Mst^{\Betti}_{g,r,0},\BQ_{\vir}\right) \xrightarrow{\cong}\HO^{\BM}\!\left(\Mst^{\Dol}_{r,0}(C),\BQ_{\vir}\right)
\end{equation}
respecting the less perverse filtrations, such that
\[
\Upsilon(W_{2i-j}\!\Gr^{\FL}_{j}\HO^{\BM}(\Mst^{\Betti}_{g,r,0},\BQ_{\vir}))= \FH^{i-j}\Gr^{\FL}_{j}\HO^{\BM}\!\left(\Mst^{\Dol}_{r,0}(C),\BQ_{\vir}\right).
\]
\end{conjx}
The above conjecture can be informally stated as saying that the weight filtration on $\HO^{\BM}(\Mst^{\Betti}_{g,r,0},\BQ_{\vir})$ records the perverse cohomological degrees of the direct image $\JH_*\BD\BQ_{\Mst^{\Dol}_{r,0}(C),\vir}$, along with twice the extra shift in perverse degrees induced by taking the direct image of this complex along the (non-stacky) Hitchin map $\Hit_{r,0}\colon\Msp^{\Dol}_{r,0}(C)\rightarrow \Lambda_r$.
\smallbreak
Finally, we say something about problem (2), by considering the summands appearing in the decomposition \eqref{DolDC}.  For now, we assume that $g\geq 2$.  By \cite[Thm.6.6]{Da21a}, for every $i\in 2\csdot\BZ_{\geq 0}$ there are inclusions
\begin{align*}
\ICS_{\Msp^{\Dol}_{r,0}(C)}\subset \Up\!\Ho^i\!\left( \JH_*\BD\BQ_{\Mst^{\Dol}_{r,0}(C),\vir}\right)\\
\ICS_{\Msp^{\Betti}_{g,r,0}}\subset \Up\!\Ho^i\!\left( \JH_*\BD\BQ_{\Mst^{\Betti}_{g,r,0},\vir}\right).
\end{align*}
Next, by the classical nonabelian Hodge correspondence, there is a homeomorphism of singular spaces
\[
\Gamma\colon \Msp^{\Betti}_{g,r,0}\rightarrow \Msp^{\Dol}_{r,0}(C)
\]
and since intersection cohomology is a topological invariant, there are isomorphisms between the intersection cohomology of the domain and the target of $\Gamma$.  These facts together at least provide an isomorphism between certain canonically defined summands of $\HO^{\BM}\!\left(\Mst^{\Dol}_{r,0}(C),\BQ_{\vir}\right)$ and $\HO^{\BM}\!\left( \Mst^{\Betti}_{g,r,0},\BQ_{\vir}\right)$:
\[
\xymatrix{
\HO^{\BM}\!\left(\Mst^{\Dol}_{r,0}(C),\BQ_{\vir}\right)\ar@{.>}[rr]_{\Upsilon}^{\exists ?}&&\HO^{\BM}\!\left( \Mst^{\Betti}_{g,r,0},\BQ_{\vir}\right)
\\
\ar@{^{(}->}[u] \ICA(\Msp^{\Dol}_{r,0}(C))\ar[rr]^{\cong}_{\Upsilon'=\Upsilon\lvert_{\ICA(\Msp^{\Dol}_{r,0}(C))}}&&\ar@{^{(}->}[u]\ICA(\Msp^{\Betti}_{g,r,0}).
}
\]

  However, there are other summands corresponding to simple perverse sheaves that do not have full support on the coarse moduli spaces, as well as nonzero non-vanoshoing perverse cohomological degrees, so that we have not yet defined the isomorphism $\Upsilon$ demanded by \eqref{SNAHT}.  In order to boot-strap a canonical definition of $\Upsilon$ out of the isomorphism $\Upsilon'$, we will utilise the theory of cohomological Hall algebras (CoHAs).  In brief, we propose to define $\Upsilon$ as the morphism of algebras induced by $\Upsilon'$.  
\smallbreak
The claim that this definition of $\Upsilon$ makes sense, and indeed provides an isomorphism, rests on various results regarding these CoHAs.  We prove analogues of these conjectures adapted to the cases $g\leq 1$.  For genus $g\leq 1$ we then go on to prove the PS=WS conjectures stated above, and this is our main theorem
\begin{thmx}
\label{ThmA}
Assume that the genus of the smooth projective curve $C$ is less than or equal to one.  Then Conjectures \ref{conj_b} and \ref{conj_c} are true
\end{thmx}

Finally, using our proposed definition of $\Upsilon$, we show how the PS=WS conjecture becomes equivalent to the original P=W conjecture and the PI=WI conjectures:

\begin{thmx}
\label{thm_b}
Let $g(C)\geq 2$.  Under the assumptions of Theorem \ref{fsta}, and assuming Conjecture \ref{BCI} (``Betti $\chi$-independence''), Theorem \ref{PW} (P=W) implies Conjectures \ref{IPIW} (PI=WI), \ref{conj_b} ($\FL^0$(PS=WS)), and \ref{conj_c} (PS=WS).
\end{thmx}
\subsection{An update}
Since this paper first appeared, proofs of Conjectures \ref{Freeg} and \ref{HGC}, along with the construction of the isomorphism  \eqref{SNAHT} required for Conjectures \ref{conj_b} and \ref{conj_c}, have all been provided in \cite{DHSM22}.  In addition, the isomorphism \eqref{BCICI} is provided in \cite[Thm.14.10]{DHSM22}.  At this stage, all that remains, amongst the assumptions of Theorem \ref{thm_b} is to prove that this ``Betti $\chi$-independence''  isomorphism preserves mixed Hodge structures.  Moreover, the original P=W conjecture has been proven \cite{PW1,PW2, MSY23b}, so that this Betti $\chi$-independence conjecture is now the only missing component to finish proving all of the conjectures mentioned in this paper.
\subsection{Acknowledgements}
This research was supported by “Categorified Donaldson–Thomas theory” No. 759967 of the European Research Council, and by a Royal Society university research fellowship.  I would like to thank Tasuki Kinjo, Andrei Okounkov, Olivier Schiffmann and Sebastian Schlegel Mejia for numerous helpful conversations relating to the paper, and the anonymous referees who made a large number of helpful suggestions.  Finally, special thanks are due to Woonam Lim who spotted a mistake in the published version of Conjecture \ref{conj_c}, which is corrected in this version.

\section{Background from cohomological DT theory}
\label{DT_sec}
\subsection{Quivers and potentials}
Fix a quiver $Q$, with finite vertex set $Q_0$ and finite arrow set $Q_1$, and two morphisms $s,t\colon Q_1\rightarrow Q_0$ taking an arrow to its source and target respectively.  We will assume for now that $Q$ is symmetric, meaning that for any pair $i,j\in Q_0$ there are as many arrows from $i$ to $j$ as there are from $j$ to $i$.  
\smallbreak
We denote by $\BC Q$ the free path algebra, having as basis over $\BC$ the set of paths in $Q$.  We include paths $e_i$ of length zero, for each $i\in Q_0$, beginning and ending at $i$.  A $\BC Q$-module $\CF$ is defined by a set of complex vector spaces $e_i\csdot \CF$ for each $i\in Q_0$ along with morphisms $e_{s(a)}\csdot \CF\rightarrow e_{t(a)}\csdot \CF$ for each $a\in Q_1$.  The \textit{dimension vector} of a module $\CF$ is the $Q_0$-tuple $(\dim(e_i\csdot \CF))_{i\in Q_0}$.
\smallbreak
We assume that we are given a subset of arrows $E\subset Q$ which we formally invert in $\BC Q$ to form the partially localised path algebra $B$.  We consider the category of $B$-modules as a Serre subcategory of the category of $\BC Q$-modules, i.e. the subcategory containing those modules $\CF$ for which the action of each $a\in E$ is an invertible morphism from $e_{s(a)}\csdot\CF$ to $e_{t(a)}\csdot \CF$.  For $\dd\in\mathbb{N}^{Q_0}$ a dimension vector, we denote by $\Mst_{\dd}(B)$ the stack of $\dd$-dimensional $B$-modules.  Note that if $E$ contains all of the edges of a tree $T$, and $T$ contains all of the vertices $Q_0$, then $\Mst_{\dd}(B)=\emptyset$ unless $\dd$ is constant, i.e. $\dd=(d,\ldots,d)$ for some $d\in\mathbb{Z}_{\geq 0}$.
\smallbreak
Throughout the paper we use the notational convention that, where an object $\mathcal{F}_{\alpha}$ is defined for each $\alpha$ in some monoid $N$, we denote by $\mathcal{F}_{\bullet}$ the coproduct over all values of $\alpha$.  So
\[
\Mst_{\bullet}(B)\coloneqq \coprod_{\dd\in\BN^{Q_0}}\Mst_{\dd}(B)
\]
is the stack of all finite-dimensional $B$-modules.  
\smallbreak
There is an open inclusion of stacks $\Mst_{\dd}(B)\subset \Mst_{\dd}(\BC Q)$, and so smoothness of $\Mst_{\dd}(B)$ follows from the smoothness of $\Mst_{\dd}(\BC Q)$.
\smallbreak
We choose a \textit{potential} $W\in \BC Q/[\BC Q,\BC Q]$.  Given the data $(B,W)$, we define the noncommutative derivatives of $W$ with respect to $a\in Q_1$ as follows.  If $W=a_1\ldots a_l$ is a single cyclic word, then
\[
\frac{\partial W}{\partial a}=\sum_{i\lvert a_i=a}a_{i+1}\ldots a_la_1\ldots a_{i-1}.
\]
We extend to arbitrary $W$ by linearity.  Then we define the Jacobi algebra
\[
\Jac(B,W)\coloneqq B/\langle \frac{\partial W}{\partial a}\;\lvert\; a\in Q_1\rangle.
\]
The element $W$ defines a function $\Tr(W)$ on $\Mst_{\bullet}(B)$: its value at a point $x\in \Mst_{\bullet}(B)$ is given by taking the trace of the endomorphism $W\cdot $ on the underlying vector space of the module that $x$ represents.  
\smallbreak
Given a pair of dimension vectors $\dd,\ee\in\BN^{Q_0}$ we define the \textit{Euler form}
\[
\chi_Q(\dd,\ee)=\sum_{i\in Q_0} \dd_i\ee_i-\sum_{a\in Q_1}\dd_{s(a)}\ee_{t(a)}.
\]
If $\Mst_{\dd}(B)$ is non-empty, its dimension is given by $-\chi_Q(\dd,\dd)$.
\subsection{Critical CoHAs}
We define $\BQ_{\Mst_{\bullet}(B),\vir}\in\Ob(\Perv(\Mst_{\bullet}(B)))$ to be the constant perverse sheaf, by setting $\vdim=\dim$, i.e. 
\begin{align*}
\BQ_{\Mst_{\bullet}(B),\vir}\lvert_{\Mst_{\dd}(B)}=&\BQ_{\Mst_{\dd}(B)}[\dim(\Mst_{\dd}(B))]\\
=&\BQ_{\Mst_{\dd}(B)}[-\chi_Q(\dd,\dd)].
\end{align*}

\smallbreak
There is an equality of susbstacks of $\Mst(B)$
\[
\crit(\Tr(W))=\Mst(\Jac(B,W)).
\]
Given a function $f$ on a complex stack $X$ we denote by $\phin{f}=\phi_f[-1]$ the shifted vanishing cycles functor.  This shifted vanishing cycles functor is exact with respect to the perverse t-structure (see \cite[Cor.10.3.13]{KSsheaves}).  Note that, since $\Mst_{\bullet}(B)$ is smooth, the perverse sheaf of vanishing cycles $\phin{\Tr(W)}\BQ_{\Mst_{\bullet}(Q),\vir}$ is supported on the critical locus of $\Tr(W)$, which we identify with the stack of $\Jac(B,W)$-modules.
\smallbreak
We consider the $\BN^{Q_0}$-graded vector space
\begin{align*}
\Coha_{B,W,\bullet}\coloneqq &\HO(\Mst_{\bullet}(B),\phin{\Tr(W)}\BQ_{\Mst_{\bullet}(B),\vir})\\
\cong&\HO(\Mst_{\bullet}(\Jac(B,W)),\phin{\Tr(W)}\BQ_{\Mst_{\bullet}(B),\vir})\\
\cong&\bigoplus_{\dd\in\BN^{Q_0}}\HO(\Mst_{\dd}(B),\phin{\Tr(W)}\BQ_{\Mst_{\dd}(B),\vir}).
\end{align*}
There is a morphism of stacks $\tau\colon \Mst_{\dd}(B)\rightarrow \pt/\BC^*$ obtained by sending a $B$-module to the determinant of its underlying vector space, which induces an action 
\begin{equation}
\label{hpta}
\HO(\pt/\BC^*,\BQ)\otimes \Coha_{B,W,\dd}\rightarrow \Coha_{B,W,\dd}
\end{equation}
for each $\dd\in\BN^{Q_0}$, since the function $\Tr(W)$ factors through the morphism 
\begin{equation}
\label{LUlt}
\Mst_{\dd}(B)\xrightarrow{\tau\times\id}\pt/\mathbb{C}^*\times\Mst_{\dd}(B).
\end{equation}
Explicitly, we obtain \eqref{hpta} by taking derived global sections of the morphism 
\begin{equation}
\phin{0\boxplus\Tr(W)}\left( \BQ_{\pt/\mathbb{C}^*\times \Mst_{\dd}(B)}\rightarrow (\tau\times\id)_*\BQ_{\Mst_{\dd}(B)}\right)[-\chi_Q(\dd,\dd)].
\end{equation}
\smallbreak
The multiplication operation for $\Coha_{B,W,\bullet}$, introduced by Kontsevich and Soibelman \cite{KS2} is defined by taking push-forward and pull-back of vanishing cycle cohomology in the usual correspondence diagram
\begin{equation}
\xymatrix{
\mathfrak{M}_{\bullet}(B)\times\mathfrak{M}_{\bullet}(B)&&\ar[ll]_-{\pi_1\times\pi_3}\mathfrak{E}\mathrm{xact}(B)\ar[rr]^-{\pi_2}&&\mathfrak{M}_{\bullet}(B).
}
\end{equation}
We have denoted by $\mathfrak{E}\mathrm{xact}(B)$ the stack of short exact sequences of finite-dimensional $B$-modules
\[
0\rightarrow \CF'\rightarrow \CF\rightarrow\CF''\rightarrow 0
\] 
and by $\pi_1,\pi_2,\pi_3$ the morphisms sending such a short exact sequence to $\CF',\CF,\CF''$ respectively.  Properness of $\pi_2$ enables us to define a pushforward map in vanishing cycle cohomology: see \cite[Sec.7]{KS2} for details of the definition of the cohomological Hall algebra defined this way.  
\smallbreak
We denote by 
\[
\JH_{\dd}\colon \Mst_{\dd}(B)\rightarrow \Msp_{\dd}(B)
\]
the morphism to the coarse moduli space, which is just the affinization of the stack $\Mst_{\dd}(B)$.  Closed points of $\Msp_{\dd}(B)$ are in natural bijection with isomorphism classes of $\dd$-dimensional semisimple $B$-modules.  The scheme $\Msp_{\bullet}(B)$ is a monoid in the category of schemes, i.e. there is a morphism
\[
\oplus\colon \Msp_{\bullet}(B)\times \Msp_{\bullet}(B)\rightarrow \Msp_{\bullet}(B)
\]
which at the level of points sends a pair of semisimple $B$-modules to their direct sum.  The morphism is finite by \cite[Lem.2.1]{Meinhardt14}, and thus the category $\Perv(\Msp_{\bullet}(B))$ acquires a \textit{convolution} tensor structure defined by
\[
\CF\boxtimes_{\oplus}\CG\coloneqq \oplus_*\!\left(\CF\boxtimes\CG\right).
\]

An \textit{algebra object} in $\Perv(\Msp_{\bullet}(B))$ is a triple $(\mathcal{F},m,i)$, where $\mathcal{F}\in\Ob(\Perv(\Msp_{\bullet}(B)))$, $m$ is a morphism $\CF\boxtimes_{\oplus}\CF\rightarrow \CF$, and $i\colon \mathbb{Q}_{\Msp_0(B)}\rightarrow \CF$ is a morphism from the monoidal unit of the tensor category $\Perv(\Msp_{\bullet}(B))$.  The morphisms $m,i$ are required to satisfy the usual conditions of an associative algebra, i.e. $m\circ (m\boxtimes_{\oplus} \Id_{\CF})=m\circ (\Id_{\CF}\boxtimes_{\oplus}m)\circ \mathtt{a}$, where $\mathtt{a}$ is the associator natural isomorphism for $\boxtimes_{\oplus}$, and $\Id_{\CF}=m\circ (\Id_{\CF}\boxtimes_{\oplus} i)=m\circ (i\boxtimes_{\oplus}\Id_{\CF})$.  We define the category of Lie algebra objects in $\Perv(\Msp_{\bullet}(B))$ similarly.
\smallbreak
The convolution tensor product defines a monoidal structure on $\Dub^+(\Perv(\Msp_{\bullet}(B)))$ and we define algebra objects in this category the same way.  A morphism of algebra objects is a morphism of underlying objects respecting the algebra structures.
\smallbreak
Set
\[
\RCoha_{B,W,\dd}\coloneqq \JH_{\dd,*}\!\phin{\Tr(W)}\BQ_{\Mst_{\dd}(B),\vir}.
\]
Then $\RCoha_{B,W,\bullet}$ carries a Hall algebra structure
\begin{equation}
\label{RCM}
\RCoha_{B,W,\bullet}\boxtimes_{\oplus}\RCoha_{B,W,\bullet}\rightarrow \RCoha_{B,W,\bullet}
\end{equation}
lifting the algebra structure on $\Coha_{B,W,\bullet}$ in the following sense: the natural $\BN^{Q_0}$-graded isomorphism of cohomologically graded vector spaces
\begin{equation}
\label{AlgIn}
\HO\!\left(\Msp_{\bullet}(B),\RCoha_{B,W,\bullet}\right)\cong \Coha_{B,W,\bullet}
\end{equation}
respects the multiplication structures, where the multiplication operation on the LHS of \eqref{AlgIn} is obtained by applying the derived global sections functor to the multiplication morphism \eqref{RCM} for $\RCoha_{B,W,\bullet}$ (see \cite{Da13,QEAs} for details).
\smallbreak
As explained in \cite[Sec.4.1]{QEAs}, it follows from the BBDG decomposition theorem and the fact that $\JH_{\dd}$ is \textit{approximated by projective maps} that there are isomorphisms
\begin{align}
\label{VCDT}
\RCoha_{B,W,\dd}\cong \bigoplus_{i\in \BZ}\UP\!\Ho^i\!\left(\RCoha_{B,W,\dd}\right)[-i].
\end{align}
We define 
\begin{equation}
\label{JHPF}
\FJ^i\!\Coha_{B,W,\dd}\coloneqq \HO(\Msp_{\dd}(B),\UP\tau^{\leq i}\RCoha_{B,W,\dd}).
\end{equation}
Because of the existence of the decomposition \eqref{VCDT} the morphism $\FJ^i\!\Coha_{B,W,\dd}\rightarrow \Coha_{B,W,\dd}$ induced by the natural transformation $\Up\tau^{\leq i}\rightarrow \id$ is an injection, and we define the \textit{perverse filtration} with respect to $\JH$, by setting the $i$th piece to be \eqref{JHPF}.
\smallbreak
Since the CoHA multiplication on $\Coha_{B,W,\bullet}$ lifts to the multiplication \eqref{RCM} for the complex $\RCoha_{B,W,\bullet}$, and the convolution tensor structure on $\Perv(\Msp_{\bullet}(B))$ is bi-exact for the perverse t-structure, the perverse filtration is respected by the multiplication on $\Coha_{B,W,\bullet}$, i.e. it restricts to operations
\[
\FJ^{i'}\!\!\Coha_{B,W,\dd'}\boxtimes_{\oplus} \:\FJ^{i''}\!\!\Coha_{B,W,\dd''}\rightarrow \FJ^{i'+i''}\!\!\Coha_{B,W,\dd'+\dd''}.
\]
In particular, we may consider the associated graded algebra with $i$-th graded piece
\[
\Gr^{\FJ}_{i}\!\Coha_{B,W,\bullet}\coloneqq \left(\FJ^{i}\!\Coha_{B,W,\bullet}\right)/\left(\FJ^{i+1}\!\Coha_{B,W,\bullet}\right).
\]
Using the (localised) coproduct on $\Coha_{B,W,\bullet}$ constructed in \cite{Da13}, it is shown in \cite{QEAs} that the algebra $\Gr^{\FJ}_{\bullet}\!\Coha_{B,W,\bullet}$ is a free supercommutative algebra\footnote{Strictly speaking, this is only true after twisting the product, restricted to $\Coha_{B,W,\dd}\otimes\Coha_{B,W,\ee}$, by $(-1)^{\psi(\dd,\ee)}$, where $\psi$ is a bilinear form satisfying a certain relation mod 2.  In this paper we will always consider $B$ such that we may take $\psi=0$ and so we ignore this complication; see \cite[Sec.1.6]{QEAs} for details.}.  Furthermore, there are equalities $\Gr^{\FJ}_{i}\!\Coha_{B,W,\bullet}=0$ for $i\leq 0$.  Thus the $\BN^{Q_0}$-graded subspace
\begin{align}
\label{bpsladef}
\Fg_{B,W,\bullet}\coloneqq \FJ^1\!\Coha_{B,W,\bullet}\subset \Coha_{B,W,\bullet}
\end{align}
is closed under the commutator Lie bracket, and is by definition the \textit{BPS Lie algebra}.  We define the \textit{BPS sheaves} to be the perverse sheaves
\begin{equation}
\label{BPSS}
\DTS_{B,W,\dd}\coloneqq (\Up\tau^{\leq 1}\!\RCoha_{B,W,\dd})[1].
\end{equation}
Then by construction, there are isomorphisms $\Fg_{B,W,\dd}\cong \HO(\Msp_{\dd},\DTS_{B,W,\dd})[-1]$.  Likewise, we define the monodromic mixed Hodge modules
\begin{equation}
\label{BPSMHM}
\underline{\DTS}_{B,W,\dd}\coloneqq \Ho^1(\JH_{\dd,*}\phim{\Tr(W)}\underline{\BQ}_{\Mst_{\dd}(B),\vir})\in\Ob(\MMHM(\Msp_{\dd}(B))).
\end{equation}
A monodromic mixed Hodge module (MMHM) can be thought of roughly as a usual mixed Hodge module, with a monodromy operator.  Within applications to nonabelian Hodge theory this monodromy will generally be trivial, up to a half Tate twist, and we refrain from giving a full account of these objects, instead referring the interested reader to \cite[Sec.2.1]{QEAs} and  \cite[Sec.2.9]{BBBD}.  The MMHMs \eqref{BPSMHM} are lifts of the BPS sheaves, in the sense that the underlying perverse sheaf of \eqref{BPSMHM} is naturally identified with \eqref{BPSS}.  
\begin{remark}
Throughout the paper, wherever $\mathcal{F}$ is an object in the derived category of constructible sheaves, or a cohomologically graded vector space which has a natural lift to a complex of MMHMs, or a cohomologically graded monodromic mixed Hodge structure respectively, we denote this lift by underlining: $\underline{\mathcal{F}}$.  
\end{remark}

We define
\[
\underline{\Fg}_{B,W,\bullet}\coloneqq \HO(\Msp_{\dd}(B),\underline{\DTS}_{B,W,\dd})[-1],
\]
providing a lift of the Lie algebra $\mathfrak{g}_{B,W,\bullet}$ to a Lie algebra object in the category of monodromic mixed Hodge structures.
\smallbreak
Given a cohomologically graded $\BQ$-vector space $V$, we denote by 
\[
\Sym^n(V)\subset \underbrace{V\otimes\cdots\otimes V}_{n\textrm{ times}}
\]
the invariant summand under the action of the symmetric group on $n$ elements, generated by the transpositions $\sigma_i$.  Note that because $V$ carries a cohomological grading, the Koszul sign rule will make an appearance here, i.e. if $a_1,\ldots,a_n$ are homogeneous elements of cohomological degrees $\lvert a_1\lvert,\ldots,\lvert a_r\lvert$, then
\[
\sigma_i(a_1\otimes\cdots \otimes a_n)=(-1)^{\lvert a_i\lvert \lvert a_{i+1}\lvert}(a_1\otimes\cdots\otimes a_{i-1}\otimes a_{i+1}\otimes a_i\otimes a_{i+2}\otimes\cdots\otimes a_n).
\]
We set
\begin{equation}
\label{SymDef}
\Sym(V)\coloneqq \bigoplus_{n\geq 0}\Sym^n(V).
\end{equation}
If $V$ is $\mathbb{N}^{Q_0}$-graded, then $\Sym(V)$ carries a natural $\mathbb{N}^{Q_0}$-grading.  We will not\footnote{This is why we have not denoted the LHS of \eqref{SymDef} $\Sym^{\bullet}(V)$, in mild violation of our stated conventions regarding the use of the symbol $\bullet$.} consider the extra grading on $\Sym(V)$ given by the direct sum decomposition \eqref{SymDef}.  For $\CF=\bigoplus_{i\in\mathbb{Z}}\CF^i[-i]$ a cohomologically perverse sheaf or MMHM\footnote{I.e. we assume that each $\CF^i$ is a perverse sheaf.  We allow for infinitely many $i$ for which $\CF^i$ is nonzero in the definition.}  on $\Msp_{\bullet}(B)$, we define
\[
\Sym_{\oplus}^n(\CF)\subset \underbrace{\CF\boxtimes_{\oplus}\cdots\boxtimes_{\oplus} \CF}_{n\textrm{ times}}
\]
and $\Sym_{\oplus}(\CF)$ as in \eqref{SymDef}.  
\smallbreak
With these preliminaries out of the way, we can now state the PBW theorem for cohomological Hall algebras of Jacobi algebras:
\begin{theorem}\cite[Thm.C]{QEAs}
\label{CIT}
(Relative integrality/PBW theorem): The morphism of complexes
\[
\Sym_{\oplus}\!\left(\bigoplus_{0\neq \dd\in\BN^{Q_0}}\underline{\DTS}_{B,W,\dd}\otimes \HO(\pt/\BC^*,\BQ)[-1]\right)\rightarrow \underline{\RCoha}_{B,W,\bullet}
\]
induced by the natural morphism $\underline{\DTS}_{B,W,\bullet}[-1]\hookrightarrow \underline{\RCoha}_{B,Q,\bullet}$, the lift of the $\HO(\pt/\BC^*,\BQ)$-action \eqref{hpta}, and the CoHA multiplication, is an isomorphism\footnote{Note: this is typically \textit{not} a morphism of algebra objects.}.
\smallbreak
(Absolute integrality/PBW theorem): Taking derived global sections, the morphism
\[
\Sym\left(\bigoplus_{\dd\in\BN^{Q_0}\setminus \{0\}}\underline{\Fg}_{B,W,\dd}\otimes \HO(\pt/\BC^*,\BQ)\right)\rightarrow \underline{\Coha}_{B,W,\bullet}
\]
is an isomorphism.
\end{theorem}

\section{Representations of preprojective algebras}
\label{propsec}
\subsection{Preprojective CoHAs}
The Donaldson--Thomas theory of preprojective algebras informs much of what follows in this paper, in two different ways.  Firstly, it is a warmup: a more linear, and tractable 2CY category in which we can prove results, the analogues of which are harder for the categories we encounter later in nonabelian Hodge theory.  Secondly, the theory for the 2CY categories appearing in nonabelian Hodge theory is locally modelled by the theory of representations of preprojective algebras, in a sense made precise by Theorem \ref{nbhd2CY}.
\smallbreak
We start, as in \S \ref{DT_sec}, with a finite quiver $Q$.  We define the doubled quiver $\overline{Q}$ by adding an arrow $a^*$ to $Q_1$ for each of the original arrows $a\in Q_1$, and setting $s(a^*)=t(a)$ and $t(a^*)=s(a)$.  Then we define the \textit{preprojective algebra}
\[
\Pi_Q\coloneqq \BC \overline{Q}/\langle \sum_{a\in Q_1} [a,a^*]\rangle.
\]
For $\dd\in \BN^{Q_0}$ we define $\Mst_{\dd}(\Pi_Q)$ to be the stack of $\dd$-dimensional $\Pi_Q$-modules.  
\smallbreak
On $\Mst_{\bullet}(\Pi_Q)$ we consider the complex $\BD\BQ_{\Mst_{\bullet}(\Pi_Q),\vir}\in\Dub(\Perv(\Mst_{\bullet}(\Pi_Q)))$.  This complex is defined as in the introduction by setting $\vdim=-2\chi_Q(\dd,\dd)$, i.e. we define
\[
\BQ_{\Mst_{\bullet}(\Pi_Q),\vir}\lvert_{\Mst_{\dd}(\Pi_Q)}=\BQ_{\Mst_{\dd}(\Pi_Q)}[-2\chi_Q(\dd,\dd)].
\]
This choice of $\vdim$ is justified by the observation that $\Mst_{\dd}(\Pi_Q)$ is the classical truncation of the derived stack of $\dd$-dimensional representations of the derived preprojective algebra, and this derived stack has virtual dimension $-2\chi_Q(\dd,\dd)$.  
\smallbreak

We define
\begin{align*}
\Coha_{\Pi_Q,\bullet}\coloneqq & \HO^{\BM}(\Mst_{\bullet}(\Pi_Q),\BQ_{\vir})\\
\coloneqq &\HO(\Mst_{\bullet}(\Pi_Q),\BD\BQ_{\Mst_{\bullet}(\Pi_Q),\vir})\\
=&\bigoplus_{\dd\in\BN^{Q_0}}\HO^{\BM}(\Mst_{\dd}(\Pi_Q),\BQ_{\vir}).
\end{align*}
This cohomologically graded $\BN^{Q_0}$-graded vector space carries a Hall algebra multiplication defined by Schiffmann and Vasserot \cite[Sec.4]{ScVa13} (see also \cite{YZ18}).  
\smallbreak
We next recall how to deduce many of the fundamental properties of the algebra $\Coha_{\Pi_Q,\bullet}$ using the cohomological DT theory of Jacobi algebras from \S \ref{DT_sec}.  Let $\tilde{Q}$ be the \textit{tripled quiver}, obtained from the doubled quiver $\overline{Q}$ by adding a loop $\omega_i$ at each vertex $i\in Q_0$.  The canonical cubic potential for $\BC \tilde{Q}$ is given by
\[
\tilde{W}=\left(\sum_{i\in Q_0}\omega_i\right)\left(\sum_{a\in Q_1} [a,a^*]\right).
\]
There is a natural isomorphism
\[
\alpha\colon\Jac(\BC \tilde{Q},\tilde{W})\cong \Pi_Q[x]
\]
with $\alpha^{-1}(x)=\sum_{i\in Q_0} \omega_i$.  We consider the commutative diagram
\begin{equation}
\label{ppcd}
\xymatrix{
\Mst_{\bullet}(\Pi_Q[x])\ar[r]^-{\varpi}\ar[d]^{\tilde{\JH}}&\Mst_{\bullet}(\Pi_Q)\ar[d]^{\JH}\\
\Msp_{\bullet}(\Pi_Q[x])\ar[r]^-{\varpi'}&\Msp_{\bullet}(\Pi_Q).
}
\end{equation}
in which the horizontal maps are the natural forgetful morphisms.  
\smallbreak
Recall that the perverse sheaf $\phin{\Tr(\tilde{W})}\BQ_{\Mst_{\bullet}(\BC \tilde{Q}),\vir}\in\Perv(\Mst_{\bullet}(\BC \tilde{Q}))$ is supported on the critical locus of $\Tr(\tilde{W})$, i.e. $\Mst_{\bullet}(\Pi[x])\subset \Mst_{\bullet}(\BC \tilde{Q})$.  Via the dimensional reduction isomorphism \cite[App.A]{Da13} there is a natural isomorphism
\begin{equation}
\beta\colon\BD\BQ_{\Mst_{\bullet}(\Pi_Q),\vir}\cong\varpi_*\phin{\Tr(\tilde{W})}\BQ_{\Mst_{\bullet}(\BC \tilde{Q}),\vir}.
\end{equation}
Taking derived global sections, there is an isomorphism
\[
\HO^{\BM}(\Mst_{\bullet}(\Pi_Q),\BQ_{\vir})\cong \HO(\Mst_{\dd}(\BC \tilde{Q}),\phin{\Tr(\tilde{W})}\BQ_{\Mst_{\bullet}(\BC \tilde{Q}),\vir}),
\]
i.e. an isomorphism of cohomologically graded, $\BN^{Q_0}$-graded vector spaces
\begin{equation}
\label{FID}
F\colon \Coha_{\Pi_Q,\bullet}\cong \Coha_{\BC \tilde{Q},\tilde{W},\bullet}.
\end{equation}
After twisting the definition of $F$ by a sign depending on $\dd\in\BN^{Q_0}$, \eqref{FID} is also an isomorphism of algebras \cite[Appendix]{RS17}, \cite{YZ16}, where the domain is given the Schiffmann--Vasserot product, and the target is given the Kontsevich--Soibelman product.
\smallbreak
The morphism $\varpi'$ is a morphism of monoids, where both schemes are given the monoid structure given by taking direct sums of semisimple modules.  As such, the complex
$\varpi'_*\RCoha_{\BC \tilde{Q},\tilde{W},\bullet}$ inherits the structure of an associative algebra in the category $\Dub^+(\Perv(\Msp(\Pi_Q)))$, i.e. we consider the unique $m'$ making the following diagram commute
\[
\xymatrix{
\varpi'_*\RCoha_{\BC \tilde{Q},\tilde{W},\bullet}\ar[drr]^-{m'}\ar[rr]^-{\varpi'_*m}&&\varpi'_*\left(\RCoha_{\BC \tilde{Q},\tilde{W},\bullet}\boxtimes_{\oplus}\RCoha_{\BC \tilde{Q},\tilde{W},\bullet}\right)\ar[d]^{\cong}\\
&&\varpi'_*\RCoha_{\BC \tilde{Q},\tilde{W},\bullet}\boxtimes_{\oplus}\varpi'_*\RCoha_{\BC \tilde{Q},\tilde{W},\bullet}
}
\]
with $m$ defined as in \eqref{RCM}.  We define 
\[
\RCoha_{\Pi_Q,\bullet}\coloneqq \JH_*\BD\BQ_{\Mst_{\bullet}(\Pi_Q),\vir}\in\Ob(\Dub^+(\Perv(\Msp_{\bullet}(\Pi_Q)))).
\]
Then via the isomorphism 
\[
\JH_*\beta\colon \RCoha_{\Pi_Q,\bullet}\xrightarrow{\cong}\varpi'_*\RCoha_{\BC \tilde{Q},\tilde{W},\bullet}
\]
the complex $\RCoha_{\Pi_Q,\bullet}$ acquires an algebra structure in the category $\Dub^+(\Perv(\Msp_{\bullet}(\Pi_Q)))$ (equipped with the convolution tensor product), from the algebra structure on $\varpi'_*\RCoha_{\BC \tilde{Q},\tilde{W},\bullet}$.  Taking derived global sections, we recover the associative algebra structure\footnote{As mentioned above, the algebra structure we define this way differs from the Schiffmann-Vasserot product by a system of signs that will not be important to us here, see the Appendix to \cite{RS17} for the geometric origin of these signs.}. on
\begin{align*}
\Coha_{\Pi_Q,\bullet}\coloneqq &\HO^{\BM}(\Mst_{\bullet}(\Pi_Q),\BQ_{\Mst_{\bullet}(\Pi_Q),\vir})\\
\cong &\HO(\Msp_{\bullet}(\Pi_Q),\RCoha_{\Pi_Q,\bullet}).
\end{align*}
\subsection{Properties of the BPS Lie algebra}
\label{BPS_pro}
As explained in \S \ref{DT_sec}, the algebra $\RCoha_{\BC\tilde{Q},\tilde{W},\bullet}$ carries a perverse filtration defined with respect to the morphism $\tilde{\JH}$, and we define the BPS Lie algebra
\begin{align*}
\Fg_{\Pi_Q,\bullet}\coloneqq &\Fg_{\BC\tilde{Q},\tilde{W},\bullet}\\
= &\FJ^1\!\Coha_{\BC\tilde{Q},\tilde{W},\bullet}\\
=&\HO\left(\Msp_{\bullet}(\Pi_Q[x]),\Up\tau^{\leq 1}\RCoha_{\BC \tilde{Q},\tilde{W},\bullet}\right)
\end{align*}
as in \eqref{bpsladef}.
\begin{proposition}\cite{preproj3}
\label{PiBPSp}
The complex 
\[
\DTS_{\Pi_Q,\bullet}\coloneqq \varpi'_*\DTS_{\BC\tilde{Q},\tilde{W},\bullet}[-1]
\]
is a semisimple perverse sheaf.  Moreover, the Lie algebra structure on
\begin{align*}
\Fg_{\Pi_Q,\bullet}=&\HO(\Msp_{\bullet}(\Pi_Q[x]),\DTS_{\BC\tilde{Q},\tilde{W},\bullet}[-1])\\
\cong&\HO(\Msp_{\bullet}(\Pi_Q),\DTS_{\Pi_Q,\bullet})
\end{align*}
lifts to a Lie algebra structure on the object $\DTS_{\Pi_Q,\bullet}$, i.e. there is an anticommutative morphism in $\Perv(\Msp_{\bullet}(\Pi_Q))$
\[
\DTS_{\Pi_Q,\bullet}\boxtimes_{\oplus}\DTS_{\Pi_Q,\bullet}\rightarrow \DTS_{\Pi_Q,\bullet}
\]
satisfying the Jacobi identity, which induces the Lie bracket on $\Fg_{\Pi_Q,\bullet}$ under taking derived global sections.  
\smallbreak
All of these statements lift to statements regarding mixed Hodge modules: the complex $\underline{\DTS}_{\Pi_Q,\bullet}\coloneqq \varpi'_*\underline{\DTS}_{\BC\tilde{Q},\tilde{W},\bullet}[-1]$ is a semisimple mixed Hodge module, pure of weight zero, and the Lie algebra structure on $\underline{\Fg}_{\Pi_Q,\bullet}\coloneqq \HO(\Msp_{\bullet}(\Pi_Q),\underline{\DTS}_{\Pi_Q,\bullet})$ lifts to a Lie algebra structure on $\underline{\DTS}_{\Pi_Q,\bullet}$.
\end{proposition}
\begin{remark}
In contrast with Proposition \ref{PiBPSp}, when we lift the BPS Lie algebra $\Fg_{\Pi_Q,\bullet}$ to a complex of sheaves on $\Msp_{\bullet}(\Pi_Q[x])$, i.e. by considering $\Up\tau^{\leq 1}\RCoha_{\BC \tilde{Q},\tilde{W},\bullet}$, it is situated in perverse degree one, and the Lie algebra structure is determined by an \textit{extension class} between perverse sheaves, not a morphism between them.  The fact that the perverse degree drops when we take the direct image along the morphism $\varpi'\colon \Msp(\Pi_Q[x])\rightarrow \Msp_{\bullet}(\Pi_Q)$ is one reason for calling the filtration of $\RCoha_{\BC \tilde{Q},\tilde{W},\bullet}$ with respect to the base $\Msp_{\bullet}(\Pi_Q)$ the ``less'' perverse filtration.
\end{remark}
The BPS Lie algebra $\Fg_{\Pi_Q,\bullet}$ has a cohomological grading, which is respected by the Lie algebra structure, since the product structure on $\Coha_{\Pi_Q,\bullet}$ respects cohomological degrees.  The next proposition describes the size of the graded pieces.  Before stating it, we recall the definition of the \textit{Kac polynomials} associated to $Q$.  In \cite{Kac83} Kac proved that for any quiver $Q$ and dimension vector $\dd\in\mathbb{N}^{Q_0}$ there is a polynomial $\kac_{Q,\dd}(q)$ satisfying the property that for prime powers $q=p^n$, the number of isomorphism classes of absolutely indecomposable $\dd$-dimensional $\mathbb{F}_q Q$-modules is equal to $\kac_{Q,\dd}(q)$.  
\begin{proposition}\cite{preproj}
\label{Kac1}
There is an equality of polynomials
\begin{align*}
\chi_{q^{1/2}}(\mathfrak{g}_{\Pi_Q,\dd})\coloneqq& \sum_{i\in\mathbb{Z}}(-1)^i\dim(\Fg^i_{\Pi_Q,\dd}) q^{i/2}\\
=&\kac_{Q,\dd}(q^{-1}).
\end{align*}
In particular $\Fg_{\Pi_Q,\dd}$ is situated in even, nonpositive cohomological degrees.
\end{proposition}

The Lie algebra structure on $\Fg_{\Pi_Q,\bullet}$ remains partly mysterious, although the zeroth cohomologically graded piece turns out to be isomorphic to a familiar and fundamental mathematical object:

\begin{proposition}\cite[Thm.6.6]{preproj3}
\label{Kac2}
Let $Q'$ denote the full real subquiver of $Q$, i.e. the quiver obtained by removing any vertices from $Q$ that support 1-cycles, along with any arrows beginning or ending at such vertices.  Then there is an isomorphism of $\BN^{Q_0}$-graded Lie algebras
\[
\bigoplus_{\dd\in\BN^{Q_0}}\Fg^0_{\Pi_Q,\dd}\cong \mathfrak{n}^-_{Q'}
\]
where the Lie algebra on the RHS is the Kac--Moody Lie algebra associated to $Q'$, and we extend the natural $\BN^{Q'_0}$-grading on the RHS to a $\BN^{Q_0}$-grading via the extension-by-zero map $\BN^{Q'_0}\hookrightarrow \BN^{Q_0}$.
\end{proposition}

\begin{remark}
Proposition \ref{Kac1} in particular implies the positivity of the coefficients of $\kac_{Q,\dd}(q)$.  This positivity was originally conjectured by Kac, and proved by Hausel, Letellier and Rodriguez-Villegas \cite{HLRV13}, see also \cite{Da13b,DGT16} for subsequent proofs, both closer to the techniques in the current paper.  The combination of Propositions \ref{Kac1} and \ref{Kac2} implies the Kac constant term conjecture: there is an equality $\dim(\Fg_{Q',\dd})=\kac_{Q,\dd}(0)$.  This conjecture was proved by Hausel \cite{Hau10} using Nakajima quiver varieties, and the proof of this result is one of the main ingredients of the proof of Proposition \ref{Kac2}.
\end{remark}
Since $\Fn^-_{Q'}$ is only one half of the Kac--Moody algebra associated to $Q'$, its presentation as a Lie algebra is especially simple.  Namely, it is the free Lie algebra with generators $f_i$ for $i\in Q'_0$, satisfying the \textit{Serre relations}.  We next recall a generalised form of these.  
\smallbreak 
We define the symmetrized Euler product on dimension vectors
\begin{align*}
(\dd',\dd'')_Q\coloneqq \chi_Q(\dd',\dd'')+\chi_Q(\dd'',\dd').
\end{align*}
\begin{definition}
Let $V_{\bullet}$ be an object in a tensor category of $\BN^{Q_0}$-graded objects; in examples $V_{\bullet}$ will either be a $\BN^{Q_0}$-graded vector space, or a perverse sheaf or mixed Hodge module on $\Msp_{\bullet}(\Pi_Q)$.  In the final two cases we give $V$ the obvious $\mathbb{N}^{Q_0}$-grading induced by the decomposition $V_{\dd}\cong \bigoplus_{\dd\in\mathbb{N}^{Q_0}}V\lvert_{\Msp_{\dd}(\Pi_Q)}$.  We assume that $V_{\dd}=0$ if there is no simple $\dd'$-dimensional $\Pi_Q$-module for any $\dd'\in\BQ\csdot\dd$.  By \cite{CB01} this is a purely combinatorial condition on $\dd$ and $Q$, see \cite[Thm.1.2]{CB01} for an exact statement.  Then we define $\Bor(V_{\bullet})$ to be the free Lie algebra object generated by $V_{\bullet}$, subject to the generalised Serre relations, which state that the morphisms
\begin{align}
\label{SerreRel}
V_{\dd'}^{\otimes(1-(\dd',\dd'')_Q)} \otimes V_{\dd''}&\rightarrow V_{l}\\\nonumber
(v_1,\ldots,v_{1-(\dd',\dd'')_Q},w)&\mapsto [v_1,[v_2,\ldots [v_{1-(\dd',\dd'')_Q},v_w]\ldots]
\end{align}
are the zero morphism for all pairs $\dd',\dd''\in\mathbb{N}^{Q_0}$ such that either 
\begin{itemize}
\item
$\dd'=1_i$ for $i\in Q'_0$ and $\dd'\neq \dd''$, or 
\item
$(\dd',\dd'')=0$.  
\end{itemize}
\end{definition}
For instance $\Fn^-_{Q'}=\Bor(S)$, where $S$ is the $\BN^{Q_0}$-graded vector space with $S_{\dd}=\BQ$ if $\dd=1_i$ with $i\in Q'_0$, and zero otherwise.  For the definition to make sense, we need the following easy result:
\begin{proposition}
Assume that there exist distinct simple $\Pi_Q$-modules $\CF',\CF''$ of dimension vectors $\dd',\dd''$ respectively.  Then $(\dd',\dd'')_Q\leq 0$.
\end{proposition}
\begin{proof}
We will assume that $Q$ is not an ADE quiver, for otherwise the result is trivial, as the only simple $\Pi_Q$-modules are 1-dimensional, and they all have different dimension vectors.  We may then identify $(-,-)_Q$ with the Euler form $\chi_{\mathscr{C}}(-,-)$ on the category of representations of the preprojective algebra $\Pi_Q$ (see \cite{Keller} and references therein), i.e.
\[
(\dim(\CF'),\dim(\CF''))_Q=\chi_{\mathscr{C}}(\CF',\CF'').
\]
By the 2-Calabi-Yau property of $\mathscr{C}$, and since $\CF',\CF''$ are distinct simple $\Pi_Q$-modules, there are equalities
\[
\dim(\Ext^0(\CF',\CF''))=\dim(\Ext^2(\CF',\CF''))=0
\]
and the result follows.
\end{proof}
The proposition ensures that the exponent $(1-(\dd',\dd'')_Q)$ appearing in \eqref{SerreRel} is strictly positive.
\begin{conjecture}
\label{BorC1}
The Lie algebra $\Fg_{\Pi_Q,\bullet}$ is one half of a $\BN^{Q_0}$-graded Borcherds algebra, i.e. there is a $\BN^{Q_0}$-graded vector space $P_{\bullet}$ and an isomorphism of Lie algebras
\[
\Fg_{\Pi_Q,\bullet}\cong\Bor(P_{\bullet}).
\]
Moreover, the Lie algebra object $\DTS_{\Pi_Q,\bullet}\in\Perv(\Msp_{\bullet}(\Pi_Q))$ is a Borcherds algebra object, i.e. there is a perverse sheaf $\CP_{\bullet}\in\Perv(\Msp_{\bullet}(\Pi_Q))$ such that 
\[
\DTS_{\Pi_Q,\bullet}\cong\Bor(\CP_{\bullet}),
\]
and similarly there is a mixed Hodge module $\underline{\CP}_{\bullet}$ such that 
\[
\underline{\DTS}_{\Pi_Q,\bullet}\cong\Bor(\underline{\CP}_{\bullet}).
\]
\end{conjecture}
The conjecture is a strengthening of the main conjecture of \cite{BoSch19}.  Assuming the conjecture, the identities of graded pieces
\begin{align*}
P_{\dd}=&\BQ\\
\CP_{\dd}=&\ICS_{\Msp_{\dd}(\Pi_Q)}
\end{align*}
for $\dd=1_i$ with $i\in Q'_0$ are forced by the definitions.  One of the main results of \cite{preproj3} is that the second equality becomes at least an inclusion for general $\dd$ such that there exists a simple $\dd$-dimensional $\Pi_Q$-module, i.e. if $\underline{\mathcal{P}}'_{\bullet}$ is \textit{any} generating object for the Lie algebra object $\underline{\DTS}_{\Pi_Q}$, there is a canonical (up to scalar multiplication) inclusion
\[
\underline{\ICS}_{\Msp_{\dd}(\Pi_Q)}\subset \underline{\CP}'_{\dd}
\]
for each $\dd$ such that there exists a simple $\dd$-dimensional $\Pi_Q$-module.  This motivates a strengthening of Conjecture \ref{BorC1}, which we come to after introducing the morphism
\begin{align*}
\Delta_{\dd,n}\colon&\Msp_{\dd}(\Pi_Q)\rightarrow \Msp_{n\!\dd}(\Pi_Q)\\
&\CF\mapsto \underbrace{\CF\oplus\ldots\oplus\CF}_{n\textrm{ times}}.
\end{align*}
\begin{conjecture}
\label{BorC2}
In Conjecture \ref{BorC1}, there are isomorphisms
\begin{enumerate}
\item
\[
\underline{\CP}_{1_i}\cong \underline{\BQ}_{\Msp_{1_i}(\Pi_Q)}
\]
for $i\in Q'_0$,
\item
\[
\underline{\CP}_{\dd}\cong \Delta_{\dd',n,*}\underline{\ICS}_{\Msp_{\dd'}(\Pi_Q)}
\]
for $\dd\in\mathbb{N}^{Q_0}$ satisfying $(\dd,\dd)_Q=0$ and $\dd=n\dd'$, with $n=\gcd(\dd_1,\ldots,\dd_{\lvert Q_0\lvert})$,
\item
\[
\underline{\CP}_{\dd}\cong \underline{\ICS}_{\Msp_{\dd}(\Pi_Q)}
\]
for $\dd\in\mathbb{N}^{Q_0}$ satisfying $(\dd,\dd)_Q<0$, and such that there exists\footnote{In fact this implies that $(\dd,\dd)_Q<0$.} a simple $\dd$-dimensional $\Pi_Q$-module,
\item
\[
\underline{\CP}_{\dd}=0
\]
otherwise.
\end{enumerate}
\end{conjecture}
\begin{remark}
\label{trich_rema}
Cases (1)-(3) match the trichotomy $(\dd,\dd)_Q>0$, $(\dd,\dd)_Q=0$, $(\dd,\dd)_Q<0$, i.e. the possible signs of the Euler form on the category of $\Pi_Q$-modules.  For a $\pi_1(\Sigma_g)$-representation $\CF$, the sign of the Euler form $\chi(\CF,\CF)$ is entirely determined by $g$, with the trichotomy matching the cases $g=0$, $g=1$ and $g>1$.  The same comment remains true for Higgs bundles on a complex curve.  As such, in genus zero, all simple roots of BPS Lie algebras will be real (type (1)), in genus one they will all be imaginary and isotropic (type (2)), and in higher genus we expect them to all be imaginary and hyperbolic (type (3)).  In words, if $g\geq 2$ we expect associated BPS Lie algebras on the Betti and Dolbeault sides to be freely generated by intersection cohomology of the coarse moduli spaces.  While this paper was being revised, this statement was proven in \cite{DHSM12}.
\end{remark}
\subsection{Decomposition theorem for moduli of objects in 2CY categories}
The cohomological DT theory of preprojective algebras governs the DT theory of more general 2CY categories via the following two theorems from \cite{Da21a}.  The first is a consequence of formality in 2CY categories (see \cite[Sec.4]{Da21a} for a general statement), and theorems on the \'etale local structure of good moduli spaces due to Alper, Hall and Rydh \cite{AHR20}.  The second theorem is a consequence of the first, and results on the cohomological DT theory of preprojective algebras in \cite{preproj3}.  See \cite{Da21a} for details and references.
\begin{theorem}\cite[Thm.5.11]{Da21a}
\label{nbhd2CY}
Let $\Mst$ be one of the stacks of objects in a 2CY category considered in this paper (i.e. representations of preprojective algebras, semistable Higgs bundles, or representations of $\BC[\pi_1(\Sigma_g)]$).  Slightly abusing notation, we denote by $\JH\colon \Mst\rightarrow \Msp$ be the morphism to the coarse moduli space in any of these examples.  Let a point $x\in\Mst$ represent the semisimple/polystable object
\[
\mathcal{F}=\bigoplus_{i\leq n}\mathcal{F}_i^{\oplus \dd_i}
\]
with $\mathcal{F}_1,\ldots,\CF_n$ distinct simple/stable objects.  Let $Q$ be a quiver with $n$ vertices indexing the objects $\mathcal{F}_i$, such that $\overline{Q}$ is isomorphic to the Ext quiver of $\{\CF_i\}$.  Then there is an affine $\Gl_{\dd}$-equivariant variety $H$ with a $\Gl_{\dd}$-fixed point $y$, and a commutative diagram of pointed stacks
\begin{equation}\label{eq:neighbourhood_diagram}
\xymatrix{
(\FM_{\mathbf{d}}(\Pi_Q),0_{\mathbf{d}})\ar[d]^{\JH}&(H/\Gl_{\mathbf{d}},y)\ar[d]^q\ar[r]^-{\ol{l}}\ar[l]_-{\ol{j}}&(\FM,x)\ar[d]^{\JH}\\
(\CM_{\mathbf{d}}(\Pi_Q),0_{\mathbf{d}})& (H/\!\!/\Gl_{\mathbf{d}},y)\ar[l]_-j\ar[r]^-l&(\CM,p(x))
}
\end{equation}
in which the horizontal arrows are \'etale and the squares are Cartesian.
\end{theorem}
As in the rest of the paper, by the quotient $H/\Gl_{\mathbf{d}}$ we mean the stack-theoretic quotient.

\begin{theorem}\cite[Thm.B]{Da21a}
\label{p2CY}
Let $\JH\colon \Mst\rightarrow \Msp$ be as in Theorem \ref{nbhd2CY}.  Then the mixed Hodge module complex $\JH_*\BD\underline{\BQ}_{\Mst,\vir}$ is pure.  As a consequence, there is an isomorphism in $\Dub^+(\Perv(\Msp))$
\[
\JH_*\BD\BQ_{\Mst,\vir}\cong \bigoplus_{i\in\BZ}{}^{\Fp}\!\Ho^i\!\left(\JH_*\BD\BQ_{\Mst,\vir}\right)[-i]
\]
and isomorphisms in $\Perv(\Msp)$
\[
{}^{\Fp}\!\Ho^i\!\left(\JH_*\BD\BQ_{\Mst,\vir}\right)\cong \bigoplus_{j\in J_i}\ICS_{Z_j}(\mathcal{L}_j)
\]
where $Z_j\subset \Msp$ are locally closed subvarieties indexed by a finite set $J_i$, $\CL_j$ is a local system on $Z_j$, and $\ICS_{Z_j}(\mathcal{L}_j)$ denotes the intermediate extension of the perverse sheaf $\mathcal{L}_j[\dim(Z_j)]$.  These decompositions lift to the derived category of mixed Hodge modules, with $\mathcal{L}_j$ lifting to simple variations of pure Hodge structure. Furthermore,
\begin{align}
\label{LPD}
{}^{\Fp}\!\Ho^i\!\left(\JH_*\BD\BQ_{\Mst,\vir}\right)=0
\end{align}
for $i\notin 2\csdot \BZ_{\geq 0}$.
\end{theorem}

\begin{definition}
Given $\JH\colon \Mst\rightarrow \Msp$ as in Theorem \ref{p2CY}, i.e. the morphism from the stack of semistable objects in a 2CY category to the good moduli space, we define\footnote{We use the symbol $\mathfrak{L}$ since in the case of representations of the preprojective algebra, this filtration is the less perverse filtration introduced in \S \ref{BPS_pro}.} the \textit{perverse filtration} (with respect to $\JH$) by setting
\[
\FL^{n}\!\HO^{\BM}(\Mst,\BQ_{\vir})=\HO(\Msp,{}^{\Fp}\tau^{\leq n}\JH_*\BD\BQ_{\Mst,\vir}).
\]
\end{definition}
As a consequence of Theorem \ref{p2CY} the natural morphism
\[
\FL^{n}\!\HO^{\BM}(\Mst,\BQ_{\vir})\rightarrow \HO(\Msp,\JH_*\BD\BQ_{\Mst,\vir})\cong \HO^{\BM}(\Mst,\BQ_{\vir})
\]
is an inclusion of cohomologically graded vector spaces, and the objects $\FL^{n}\!\HO^{\BM}(\Mst,\BQ_{\vir})$ provide an increasing filtration of $\HO^{\BM}(\Mst,\BQ_{\vir})$.  Furthermore, these morphisms lift to the category of mixed Hodge structures, i.e. there is a filtration of $\HO^{\BM}(\Mst,\underline{\BQ}_{\vir})$ by mixed Hodge structures 
\[
\FL^{n}\!\HO^{\BM}(\Mst,\underline{\BQ}_{\vir})\coloneqq \HO(\Msp,\tau^{\leq n}\JH_*\BD\underline{\BQ}_{\Mst,\vir}).
\]
\smallbreak
By \eqref{LPD}, this perverse filtration begins in degree zero, and the associated graded object with respect to this filtration is concentrated entirely in even non-negative degrees.
\section{The Betti side}
\label{Betti_sec}
Let $\Sigma_g$ be a genus $g$ Riemann surface, and let $\Sigma'_g$ denote $\Sigma_g$ with a single puncture.  The standard presentations of the fundamental groups of $\Sigma_g$ and $\Sigma'_g$ are
\begin{align*}
\pi_1(\Sigma_g)=&\langle a_1,\ldots,a_g,b_1,\ldots,b_g\lvert \prod_{i=1}^g[a_i,b_i]\rangle\\
\pi_1(\Sigma'_g)=&\langle a_1,\ldots,a_g,b_1,\ldots,b_g\rangle.
\end{align*}
We denote by $\Mst^{\Betti}_{g,r,d}$ the moduli stack of $r$-dimensional representations of $\pi_1(\Sigma_g)$ for which the monodromy around the puncture is given by multiplication by $\exp(2d\pi \sqrt{-1}/r)$.  We identify $\Mst^{\Betti}_{g,r,0}$ with the moduli stack of $r$-dimensional representations of $\pi_1(\Sigma_g)$, i.e. 
\[
\Mst^{\Betti}_{g,r,0}=\Mst_r(\BC[\pi_1(\Sigma_g)]).
\]
There is an isomorphism
\[
\Mst^{\Betti}_{g,r,d}\cong R_{g,r,d}/\Gl_r(\BC)
\]
where $R_{g,r,d}\subset \Gl_r(\BC)^{\times 2g}$ is the variety defined by the matrix-valued relation
\[
\prod_{i=1}^g[A_i,B_i]=\exp(2d\pi  \sqrt{-1}/r)\Id_{r\times r},
\]
where the commutator on the left hand side is again the group theoretic one, and $\Gl_r(\BC)$ acts via simultaneous conjugation.  We denote by
\[
\Msp^{\Betti}_{g,r,d}=\Spec(\Gamma(R_{g,r,d})^{\Gl_r(\mathbb{C})})
\]
the coarse moduli space, and by
\[
\JH\colon \Mst^{\Betti}_{g,r,d}\rightarrow \Msp^{\Betti}_{g,r,d}
\]
the affinization map.  If $r$ and $d$ are coprime, $\Msp^{\Betti}_{g,r,d}$ is a smooth variety, $\JH$ is a $\B \mathbb{C}^*$-gerbe, and there is an isomorphism in cohomology
\[
\HO(\Mst^{\Betti}_{g,r,d},\mathbb{Q})\cong \HO(\Msp^{\Betti}_{g,r,d},\mathbb{Q})\otimes\HO(\pt/\mathbb{C}^*,\mathbb{Q}).
\]

\subsection{Betti BPS sheaves}
Set 
\[
\Coha^{\Betti}_{g,r}\coloneqq \HO^{\BM}(\Mst^{\Betti}_{g,r,0},\mathbb{Q}_{\vir}).
\]
In \cite{Da16} it is explained how to construct a partially localised\footnote{I.e. we formally invert some of the arrows in a usual free path algebra $\BC Q_{\Delta}$ for some quiver $Q_{\Delta}$.} quiver path algebra $B_g=\BC\tilde{Q}_{\Delta}$, along with a potential $W\in \BC Q_{\Delta}/[\BC Q_{\Delta},\BC Q_{\Delta}]$, such that there is a natural isomorphism
\[
\Jac(B_g,W)\cong \Mat_{s\times s}(\BC[\pi_1(\Sigma_g)][x])
\]
for $s$ the number of vertices of $Q_{\Delta}$.  The set of edges that we localise forms a maximal tree in $Q_{\Delta}$, so that the only dimension vectors $\dd\in\mathbb{N}^{Q_{\Delta,0}}$ for which $\Mst_{\dd}(B_g)$ is nonempty are constant.  There is a natural equivalence of stacks
\[
\Mst_{(r,\ldots,r)}(\Mat_{s\times s}(\BC[\pi_1(\Sigma_g)][x]))\cong \Mst_r(\BC[\pi_1(\Sigma_g)][x])
\]
and thus we can identify $\Mst_r(\BC[\pi_1(\Sigma_g)][x])$ with the support of the perverse sheaf $\phin{\Tr(W)}\BQ_{\Mst_{(r,\ldots,r)}(B_g),\vir}$.  We abbreviate
\[
\overline{r}=(r,\ldots,r).
\]
The quiver $Q_{\Delta}$ is not symmetric, but since we only need to consider constant dimension vectors $\dd,\dd'$, which obviously satisfy $\chi_{Q_{\Delta}}(\dd,\dd')=\chi_{Q_{\Delta}}(\dd',\dd)$, the results of \S \ref{DT_sec} (in particular Theorem \ref{CIT}) apply without modification.

We denote by 
\begin{align*}
\varpi\colon &\Mst_r(\BC[\pi_1(\Sigma_g)][x])\rightarrow \Mst^{\Betti}_{g,r,0}\\
\varpi'\colon &\Msp_r(\BC[\pi_1(\Sigma_g)][x])\rightarrow \Msp^{\Betti}_{g,r,0}
\end{align*}
the forgetful morphisms, taking $\BC[\pi_1(\Sigma_g)][x]$-modules to their underlying $\BC[\pi_1(\Sigma_g)]$-modules (i.e. forgetting the action of $x$).  Then there is a commutative diagram (analogous to \eqref{ppcd})
\begin{equation}
\xymatrix{
\Mst_r(\BC[\pi_1(\Sigma_g)][x])\ar[d]^{\tilde{\JH}}\ar[r]^-{\varpi}&\Mst^{\Betti}_{g,r,0}\ar[d]^{\JH}\\
\Msp_r(\BC[\pi_1(\Sigma_g)][x])\ar[r]^-{\varpi'}&\Msp^{\Betti}_{g,r,0}
}
\end{equation}
where we denote by
\[
\tilde{\JH}\colon \Mst_r\left(\BC[\pi_1(\Sigma_g)][x]\right)\rightarrow \Msp_r\left(\BC[\pi_1(\Sigma_g)][x]\right)
\]
the canonical morphism from the stack of $\BC[\pi_1(\Sigma_g)][x]$-modules to the coarse moduli space.  Then using the dimensional reduction isomorphism (see \cite[Prop.5.3]{Da16}) there is a natural isomorphism
\begin{equation}
\label{drB}
\varpi_*\phin{\Tr(W)}\BQ_{\Mst_{\overline{r}}(B_g),\vir}\cong \BD\BQ_{\Mst^{\Betti}_{g,r,0},\vir}.
\end{equation}

The complex $\RCoha_{B_g,W,\bullet}=\tilde{\JH}_*\phin{\Tr(W)}\BQ_{\Mst_{\bullet}(B_g),\vir}$ carries the cohomological Hall algebra structure 
\[
m\colon\RCoha_{B_g,W,\bullet}\boxtimes_{\oplus}\RCoha_{B_g,W,\bullet}\rightarrow \RCoha_{B_g,W,\bullet} 
\]
of \S \ref{DT_sec} for a localised path algebra with potential.  Since $\varpi'$ is a morphism of monoids in the category of schemes, it follows that the derived direct image
\begin{align*}
\RCoha^{\Betti}_{g,\bullet}\coloneqq &\bigoplus_{r\geq 0}\JH_*\BD\BQ_{\Mst^{\Betti}_{g,\bullet,0},\vir}\\
\cong&\varpi'_*\RCoha_{B_g,W,\bullet}
\end{align*}
also carries a cohomological Hall algebra structure, i.e. the unique $m'$ making the diagram
\[
\xymatrix{
\ar[drr]^{m'}\varpi'_*\RCoha_{B_g,W,\bullet}\ar[rr]^-{\varpi'_*m}&& \varpi'_*\left(\RCoha_{B_g,W,\bullet}\boxtimes_{\oplus}\RCoha_{B_g,W,\bullet}\right)\ar[d]^{\cong}\\
&&\varpi'_*\RCoha_{B_g,W,\bullet}\boxtimes_{\oplus}\varpi'_*\RCoha_{B_g,W,\bullet}
}
\]
commute.  Similarly, the hypercohomology
\begin{align*}
\Coha^{\Betti}_{g,\bullet}\coloneqq &\bigoplus_{r\geq 0}\HO^{\BM}(\Mst^{\Betti}_{g,r,0},\BQ_{\vir})\\
\cong &\bigoplus_{r\geq 0}\HO(\Msp_{r}(\BC[\pi_1(\Sigma_g)]),\varpi'_*\tilde{\JH}_*\phin{\Tr(W)}\BQ_{\Mst_{\overline{r}}(B_g),\vir})\\\cong&\Coha_{B_g,W,\bullet}
\end{align*}
carries a Hall algebra structure, which moreover lifts to an algebra object in the derived category of mixed Hodge structures.
\smallbreak
By Theorem \ref{p2CY} the right hand side of \eqref{drB} satisfies the BBDG decomposition theorem for the morphism $\JH$, and there is an isomorphism
\[
\JH_*\BD\BQ_{\Mst^{\Betti}_{g,r,0},\vir}\cong \bigoplus_{i\in 2\cdot \BZ_{\geq 0}}\!\!\UP\Ho^i\!\left(\JH_*\BD\BQ_{\Mst^{\Betti}_{g,r,0},\vir}\right)[-i]
\]
with each of the $\UP\Ho^i\!\left(\JH_*\BD\BQ_{\Mst^{\Betti}_{g,r,0},\vir}\right)$ semisimple perverse sheaves.

\smallbreak
Via the (relative) cohomological integrality theorem (Theorem \ref{CIT}), we have the isomorphism
\begin{equation}
\label{ICB}
\tilde{\JH}_*\phin{\Tr(W)}\BQ_{\Mst_{\bullet}(B_g),\vir}\cong \Sym_{\oplus}\left (\bigoplus_{r\geq 1} \DTS_{B_g,W,\overline{r}}\otimes \HO(\pt/\BC^*,\BQ)[-1]\right)
\end{equation}
for certain perverse sheaves $\DTS_{B_g,W,\overline{r}}\in \Perv\left( \Msp_r(\BC[\pi_1(\Sigma_g)][x])\right)$.  If we assume $g\geq 1$, then explicitly:
\[
\DTS_{B_g,W,\overline{r}}\coloneqq \phin{\Tr(W)}\ICS_{\Msp_{\overline{r}}(B_g)}.
\]
Alternatively, noting that $\HO(\pt/\BC^*,\BQ)\cong \bigoplus_{i\geq 0}\BQ[-2i]$ we deduce from \eqref{ICB} that
\[
\DTS_{B_g,W,\overline{r}}\cong {}^{\mathfrak{p}}\!\Ho^1\!\left(\tilde{\JH}_*\phi_{\Tr(W)}\BQ_{\Mst_{\overline{r}}(B_g),\vir}\right).
\]
Applying $\varpi'_*$ to \ref{ICB}, from \eqref{drB} we deduce the cohomological integrality theorem for the Betti moduli stack:
\begin{align}
\label{JAR}
\RCoha^{\Betti}_{g,\bullet}\cong &\Sym_{\oplus}\left (\bigoplus_{r\geq 1} \varpi'_*\DTS_{B_g,W,\overline{r}}\otimes \HO(\pt/\BC^*,\BQ)[-1]\right)
\end{align}
as well as the analogous result at the level of mixed Hodge modules for the lift $\underline{\RCoha}^{\Betti}_{g,\bullet}$ of $\RCoha^{\Betti}_{g,\bullet}$ to the category of complexes of mixed Hodge modules.
\begin{proposition}
\label{2BP}
The complex $\varpi'_*\DTS_{B_g,W,\overline{r}}[-1]$ is a perverse sheaf, for which the lift to a mixed Hodge module is semisimple and pure of weight zero.
\end{proposition}
\begin{proof}
Since the morphism $\varpi'$ is affine, the functor $\varpi'_*$ is right t-exact with respect to the perverse t-structure, so $\varpi'_*\DTS_{B_g,W,\overline{r}}[-1]\in {}^{\mathfrak{p}}\!\mathcal{D}^{\leq 1}(\Perv(\Msp^{\Betti}_{g,r,0}))$.  By \eqref{JAR} there is a morphism of complexes
\begin{equation}
\label{DTBin}
\varpi'_*\DTS_{B_g,W,\overline{r}}[-1]\rightarrow \RCoha^{\Betti}_{g,r}
\end{equation}
which moreover has a left inverse, since $\BQ[-1]$ is a direct summand of 
\[
\HO(\pt/\BC^*,\BQ)[-1]\cong \bigoplus_{i\in\BZ_{\geq 0}}\BQ[-2i-1].
\]
On the other hand, by Theorem \ref{p2CY} the perverse cohomology sheaves of $\RCoha^{\Betti}_{g,r}$ vanish in odd degrees, and strictly negative degrees, and so the morphism 
\[
{}^{\mathfrak{p}}\!\Ho^0\!\left(\RCoha^{\Betti}_{g,r}\right)=\UP\tau^{\leq 0}\!\RCoha^{\Betti}_{g,r}\rightarrow \UP\tau^{\leq 1}\!\RCoha^{\Betti}_{g,r}
\]
is an isomorphism, and ${}^{\mathfrak{p}}\tau^{\leq 1}\!\RCoha^{\Betti}_{g,r}$ is perverse.  It follows that $\varpi'_*\DTS_{B_g,W,\overline{r}}[-1]$ is a perverse subsheaf of the perverse sheaf ${}^{\mathfrak{p}}\!\Ho^0\!\left(\RCoha^{\Betti}_{g,r}\right)$, for which the lift to a mixed Hodge module is semisimple and pure of weight zero by Theorem \ref{p2CY}.
\end{proof}
We have shown that the complex of mixed Hodge modules $\varpi'_*\underline{\DTS}_{B_g,W,\overline{r}}[-1]$ is a semisimple mixed Hodge module.  We define the perverse sheaf (respectively, mixed Hodge module)
\begin{align*}
\DTS^{\Betti}_{g,r}\coloneqq&\varpi'_*\DTS_{B_g,W,\overline{r}}[-1]\\
\underline{\DTS}^{\Betti}_{g,r}\coloneqq&\varpi'_*\underline{\DTS}_{B_g,W,\overline{r}}\otimes\LLL^{1/2}.
\end{align*}
\smallbreak
Combining the MHM version of \eqref{JAR} with Proposition \ref{2BP}, the MMHM complex
\begin{align}
\label{JAR2}
\underline{\RCoha}^{\Betti}_{g,\bullet}\cong &\Sym_{\oplus}\left (\bigoplus_{r\geq 1} \ul{\DTS}^{\Betti}_{g,r}\otimes \HO(\pt/\BC^*,\ul{\BQ})\right)
\end{align}
is a split complex of pure mixed Hodge modules (i.e. it is isomorphic to its total cohomology).
\begin{remark}
A longer, but probably more illuminating way to prove Proposition \ref{2BP} would be by proving an analogue of the support lemma regarding BPS sheaves from \cite{preproj}.  I.e. adapting the proof from \cite{preproj}, or arguing as in the appendix of \cite{KiKo21}, it should be possible to to show that the support of $\DTS_{B_g,W,\overline{r}}$ is contained in the subvariety of $\Msp_r(\BC[\pi_1(\Sigma_g)][x])$ for which the operator $x$ has a single eigenvalue, and is moreover $\mathbb{A}^1$-equivariant for the $\mathbb{A}^1$-action that acts on $\BC[\pi_1(\Sigma_g)][x]$ via the substitution $x\mapsto x+t$, with $t\in\BA^1$.  Proposition \ref{2BP} follows easily from these statements.
\end{remark}
We define $\Fg^{\Betti}_{g,\bullet}=\bigoplus_{r\geq 1} \Fg^{\Betti}_{g,r}$ where
\[
\Fg^{\Betti}_{g,r}\coloneqq \HO\left(\Msp^{\Betti}_{g,r,0},\DTS^{\Betti}_{g,r}\right).
\]
Applying the functor $\HO(\Msp^{\Betti}_{g,r,0},-)$ to the mixed Hodge module version of \eqref{JAR} yields the isomorphism of mixed Hodge modules (not preserving algebra structures)
\begin{align}
\label{JAA}
\underline{\Coha}^{\Betti}_{g,\bullet}\cong &\Sym\!\left(\bigoplus_{r\geq 1} \underline{\Fg}^{\Betti}_{g,r}\otimes \HO(\pt/\BC^*,\ul{\BQ})\right).
\end{align}
We recall that the mixed Hodge structure on $\HO(\pt/\BC^*,\ul{\BQ})$ is pure, i.e.
\[
\HO(\pt/\BC^*,\ul{\BQ})\cong\bigoplus_{i\in\BZ_{\geq 0}}\LLL^i.
\]
The isomorphism \eqref{JAA} tells us that, although the Borel--Moore homology of the stacks $\Mst_{g,r,0}^{\Betti}$ is very large, and impure, it is generated in an elementary way by something much more manageable (the BPS Lie algebra $\underline{\Fg}^{\Betti}_{g,\bullet}$), and its impurity is also controlled by this Lie algebra.
\smallbreak
Since \eqref{DTBin} has a right inverse, the morphism 
\begin{equation}
\label{fgBin}
j\colon \Fg^{\Betti}_{g,\bullet}\rightarrow \Coha^{\Betti}_{g,\bullet}
\end{equation}
obtained by taking hypercohomology is an inclusion, and is by definition the inclusion of the BPS Lie algebra, and factors through an inclusion
\[
j'\colon \Fg^{\Betti}_{g,\bullet}\hookrightarrow \FL^{ 0}\!\Coha^{\Betti}_{g,\bullet}.
\]
The following is proved exactly the same way as \cite[Thm.7.35]{preproj3}.
\begin{proposition}
\label{UEAB}
There is an isomorphism
\[
\FL^{0}\!\Coha^{\Betti}_{g,\bullet}\cong \UEA\!\left(\Fg^{\Betti}_{g,\bullet}\right)
\]
identifying $j'$ with the natural inclusion of a Lie algebra inside its universal enveloping algebra.
\end{proposition}
By \cite[Cor.6.10]{QEAs} the commutator bracket 
\[
[-,-]\colon \Up\tau^{\leq 1}\RCoha_{B_g,W,\bullet}\boxtimes_{\oplus}\Up\tau^{\leq 1}\RCoha_{B_g,W,\bullet}\rightarrow \Up\tau^{\leq 2}\RCoha_{B_g,W,\bullet}
\]
factors through the morphism $\Up\tau^{\leq 1}\RCoha_{B_g,W,\bullet}\rightarrow \Up\tau^{\leq 2}\RCoha_{B_g,W,\bullet}$.  Applying $\varpi'_*$, we deduce that the commutator bracket
\[
\DTS_{g,\bullet}^{\Betti}\boxtimes_{\oplus}\DTS_{g,\bullet}^{\Betti}\rightarrow \RCoha^{\Betti}_{g,\bullet}
\]
factors through the split morphism $\DTS_{g,\bullet}^{\Betti}\rightarrow \RCoha^{\Betti}_{g,\bullet}$.  In other words:
\begin{proposition}
The commutator Lie bracket of the algebra structure on $\RCoha^{\Betti}_{g,\bullet}$ induces the structure of a Lie algebra object on the perverse sheaf $\DTS_{g,\bullet}^{\Betti}$, from which we recover the Lie algebra $\Fg^{\Betti}_{g,\bullet}$ by taking derived global sections.  The same result holds at the level of mixed Hodge modules.
\end{proposition}
\subsection{Generators}
By Theorem \ref{p2CY} the complex 
\[
\mathcal{U}_{\bullet}=\Up\tau^{\leq 0}\RCoha^{\Betti}_{g,\bullet}
\]
is a perverse sheaf, and moreover by applying $\Up\tau^{\leq 0}$ to the CoHA multiplication morphism for $\RCoha^{\Betti}_{g,\bullet}$ we obtain the structure of an algebra object on $\CU_{\bullet}$.  The algebra object $\FL^{0}\!\Coha^{\Betti}_{g,\bullet}$ is obtained by  applying the derived global sections functor to the algebra object $\CU_{\bullet}$.  By Theorem \ref{p2CY} again, the lift of $\CU_{\bullet}$ to the category of mixed Hodge modules is semisimple.  
\smallbreak
Now assume that $g\geq 2$.  By the proof of \cite[Thm.6.6]{Da21a}, for all $r\geq 1$ there is an inclusion
\[
\iota\colon \ICS_{\Msp^{\Betti}_{g,r,0}}\hookrightarrow \CU_{r}
\]
as a \textit{primitive} subspace, meaning that there is a (unique) decomposition of perverse sheaves
\[
\CU_r\cong \ICS_{\Msp^{\Betti}_{g,r,0}}\oplus \mathcal{C}_r
\]
such that for all $r',r''\in\BZ_{>0}$ with $r'+r''=r$ the multiplication morphism
\[
\CU_{r'}\boxtimes_{\oplus} \CU_{r''}\rightarrow \CU_r
\]
factors through $\mathcal{C}_r$.  
\begin{proposition}
Let $g\geq 2$.  The morphism $\iota$ factors through the inclusion of perverse sheaves
\[
i\colon \DTS^{\Betti}_{g,r}\hookrightarrow \CU_r.
\]
Taking hypercohomology, the natural inclusions
\begin{equation}
\label{BCusp}
\HO(\iota)\colon\ICA(\Msp^{\Betti}_{g,r,0})\hookrightarrow \FL^{0}\!\Coha^{\Betti}_{g,r}
\end{equation}
embedding intersection cohomology of the coarse moduli spaces $\Msp^{\Betti}_{g,r}$ as primitive generators of $\Coha^{\Betti}_{g,r,0}$ factor through the inclusion of the BPS Lie algebra $\Fg_{g,\bullet}^{\Betti}\hookrightarrow \FL^{0}\!\Coha^{\Betti}_{g,\bullet}$.  The same result holds at the level of mixed Hodge structures/MHMs.
\end{proposition}
\begin{proof}
The statement for perverse sheaves follows from the statement for mixed Hodge modules, which we prove.  
\smallbreak
Since the morphism $\iota$ is an inclusion of a simple summand inside a semisimple mixed Hodge module, and $i$ is an inclusion of a summand inside a semisimple mixed Hodge module, the proposition follows from two claims; firstly, the claim that $\CU_r$ has a \textit{unique} summand with full support, and secondly that $\underline{\DTS}^{\Betti}_{g,r}$ also contains a summand with full support.  
\smallbreak
The first claim follows from the proof of \cite[Thm.6.6]{Da21a}, so we just need to show that $\underline{\DTS}^{\Betti}_{g,r}$ contains a summand with full support.  We use the characterisation
\[
\underline{\DTS}^{\Betti}_{g,r}=\varpi'_*\Up\tau^{\leq 1}\!\tilde{\JH}_*\phin{\Tr(W)}\underline{\BQ}_{\Mst_{\overline{r}}(B_g),\vir}.
\]
\smallbreak
Let 
\begin{align*}
\Mst^{\simp}_{\overline{r}}(B_g)\subset \Mst_{\overline{r}}(B_g)\\
\Msp^{\simp}_{\overline{r}}(B_g)\subset \Msp_{\overline{r}}(B_g)
\end{align*}
denote the substack (respectively, subvariety) of simple $\overline{r}$-dimensional $B_g$-modules.  The stack $\Mst^{\simp}_{\overline{r}}(B_g)$ is smooth, as it is an open substack of the stack $\Mst_{\overline{r}}(B_g)$, which is itself an open substack of the (smooth) stack of $\mathbb{C}Q_{\Delta}$-modules.  The critical locus of the function $\Tr(W)$ on $\Mst^{\simp}_{\overline{r}}(B_g)$ is isomorphic to the stack $\Mst^{\simp}_r(\BC[\pi_1(\Sigma_g)][x])$ of $r$-dimensional simple $\BC[\pi_1(\Sigma_g)][x]$-modules.  Let $\CF$ be such a module.  Since $x$ is central in $\pi_1(\Sigma_g)[x]$, $\CF$ admits a decomposition
\[
\CF\cong \bigoplus_{\lambda\in\BC}\CF_{\lambda}
\]
as a $\pi_1(\Sigma_g)[x]$-module, where $\CF_{\lambda}$ is the generalised eigenspace for the operator $x\cdot$, with eigenvalue $\lambda$.  By simplicity, there is only one generalised eigenvalue.  Since $\Ker((x\cdot) -\lambda \Id_{\CF})$ is then a submodule of $\CF$, we find that (again by simplicity) $x \cdot$ acts on $\CF$ via multiplication by $\lambda$.  Define the morphism of stacks 
\[
\beta\colon \Mst^{\simp}_r(\BC[\pi_1(\Sigma_g)])\times\BA^1\rightarrow \Mst^{\simp}_r(\BC[\pi_1(\Sigma_g)][x])
\]
taking a pair $(\CF,t)$ to the $\BC[\pi_1(\Sigma_g)][x]$-module given by extending the given $\BC[\pi_1(\Sigma_g)]$-action on $\CF$ by setting $x\cdot= t\Id_{\CF}$.  The domain of $\beta$ is smooth: the dimension of the Zariski tangent space of any point of $\Mst^{\simp}_r(\BC[\pi_1(\Sigma_g)])$ representing a module $\rho$ is given by $2r^2(g-1)-\dim(\Hom(\rho,\rho))-\dim(\Ext^2(\rho,\rho))$, which is constant, since $\Ext^2(\rho,\rho)^{\vee}\cong \Hom(\rho,\rho)\cong \mathbb{C}$.  Then we have shown that the morphism $\beta$ is an isomorphism of stacks, and in particular, $\Mst^{\simp}_r(\BC[\pi_1(\Sigma_g)][x])$ is smooth.
\smallbreak
Since $\crit(\Tr(W))\cap \Mst^{\mathrm{simp}}_{\overline{r}}(B_g)$ is a smooth stack, it follows from the holomorphic Bott-Morse lemma that $\phin{\Tr(W)}\BQ_{\Mst^{\simp}_{\overline{r}}(B_g),\vir}$ is analytically locally isomorphic to the constant perverse sheaf $\BQ_{\Mst^{\simp}_r(\BC[\pi_1(\Sigma_g)][x])}[d]$, where 
\begin{align}
d=&\dim\left(\Mst^{\simp}_r(\BC[\pi_1(\Sigma_g)][x])\right)\\
=&\dim\left(\Mst^{\simp}_r(\BC[\pi_1(\Sigma_g)])\right)+1.\label{dcit}
\end{align}
The projection morphism 
\[
\tilde{p}\colon \Mst^{\simp}_r(\BC[\pi_1(\Sigma_g)][x])\rightarrow \Msp^{\simp}(\BC[\pi_1(\Sigma_g)][x])
\]
has fibres isomorphic to $\pt/\BC^*$, and in particular, the fibres have trivial fundamental groups.  It follows that the monodromy of the perverse sheaf $\phin{\Tr(W)}\BQ_{\Mst^{\simp}_{\overline{r}}(B_g),\vir}$ is trivial after restricting to the fibres of $\tilde{p}$, and so, possibly after pulling back along an analytic cover $l\colon U\rightarrow \Msp^{\simp}(\BC[\pi_1(\Sigma_g)][x])$, we have 
\begin{align*}
l^*\tilde{p}_*\phin{\Tr(W)}\BQ_{\Mst^{\simp}_{\overline{r}}(B_g)}\cong &l^*\BQ_{\Msp_r^{\simp}(\BC[\pi_1(\Sigma_g)][x])}[d]\otimes\HO(\pt/\BC^*)\\
\cong&l^*\bigoplus_{i\geq 0}\BQ_{\Msp_r^{\simp}(\BC[\pi_1(\Sigma_g)][x])}[d-2i].
\end{align*}
Since 
\[
d=\dim(\Msp^{\simp}(\BC[\pi_1(\Sigma_g)][x]))-1
\]
it follows that 
\[
l^*({}^{\Fp}\tau^{\leq 1}\tilde{p}_*\phin{\Tr(W)}\BQ_{\Mst^{\simp}_{\overline{r}}(B_g)})\cong l^*\BQ_{\Msp^{\simp}(\BC[\pi_1(\Sigma_g)][x])}[\dim(\Msp^{\simp}(\BC[\pi_1(\Sigma_g)]))]
\]
in other words, $\DTS_{B_g,\overline{r}}$ is locally isomorphic to the constant perverse sheaf after restriction to $\Msp^{\simp}_{r}(\BC[\pi_1(\Sigma_g)][x])$.  Since the fibres of $\varpi'\colon \Msp^{\simp}_{r}(\BC[\pi_1(\Sigma_g)][x])\rightarrow \Msp^{\simp}_{r}(\BC[\pi_1(\Sigma_g)])$ are isomorphic to $\BA^1$, it follows that $\underline{\DTS}^{\Betti}_{g,r}=\varpi'_*\underline{\DTS}_{B_g,\overline{r}}[-1]$ has a summand with full support on $\Msp^{\simp}_r(\BC[\pi_1(\Sigma_g)])$.  Since $\Msp^{\simp}_r(\BC[\pi_1(\Sigma_g)])$ is dense in $\Msp_r(\BC[\pi_1(\Sigma_g)])$ this completes the proof.
\end{proof}

The following is the analogue of Conjecture \ref{BorC1} for the category of $\BC[\pi_1(\Sigma_g)]$-modules:

\begin{conjecture}
\label{Freeg}
For $g\geq 2$ the Lie algebra $\mathfrak{g}^{\Betti}_{g,\bullet}=\bigoplus_{r\geq 1} \mathfrak{g}^{\Betti}_{g,r}$ is the free Lie algebra generated by $\bigoplus_{r\geq 1} \ICA(\Msp_{g,r,0}^{\Betti})$.  Furthermore, the Lie algebra object $\DTS^{\Betti}_{g,\bullet}$ is the free Lie algebra object in $\Perv(\Msp^{\Betti}_{g,\bullet})$ generated by $\bigoplus_{r\geq 1}\ICS_{\Msp_{g,r,0}^{\Betti}}$.  The same statements are true at the level of mixed Hodge structures and mixed Hodge modules, i.e. if we replace $\Fg^{\Betti}_{g,r}$ and $\ICA(\Msp_{g,r,0}^{\Betti})$ by their lifts $\underline{\Fg}^{\Betti}_{g,r}$ and $\underline{\ICA}(\Msp_{g,r,0}^{\Betti})$ (respectively) to mixed Hodge structures, and replace the perverse sheaves $\DTS^{\Betti}_{g,\bullet}$ and $\ICS_{\Msp_{g,r,0}^{\Betti}}$ by the mixed Hodge modules $\underline{\DTS}^{\Betti}_{g,\bullet}$ and $\underline{\ICS}_{\Msp_{g,r,0}^{\Betti}}$ (respectively).
\end{conjecture}

In view of Proposition \eqref{UEAB}, Conjecture \ref{Freeg} is equivalent to the conjecture 
\begin{conjecture}
\label{freeA}
For $g\geq 2$ the algebra $\FL^0\!\Coha^{\Betti}_{g,\bullet}$ is the free associative algebra generated by $\bigoplus_{r\geq 1}\ICA(\Msp_{g,r,0}^{\Betti})$.  Furthermore, the algebra object ${}^{\Fp}\tau^{\leq 0}\!\RCoha^{\Betti}_{g,\bullet}$ is the free algebra object in $\Perv(\Msp^{\Betti}_{g,\bullet,0})$ generated by $\bigoplus_{r\geq 1}\ICS_{\Msp_{g,r,0}^{\Betti}}$.  The same statements are true at the level of mixed Hodge structures and mixed Hodge modules, where the replacements are as in Conjecture \ref{Freeg}.
\end{conjecture}
While this paper was being revised, both conjectures were proved in joint work with Hennecart and Schlegel Mejia \cite{DHSM12}.
\begin{remark}
\label{trich_rem}
We explain, in outline, how in fact Conjecture \ref{freeA} follows from Conjecture \ref{BorC1} and the \'etale neighbourhood theorem (Theorem \ref{nbhd2CY}).  Assuming (as an extension of Theorem \ref{nbhd2CY}) that the relative CoHA $\RCoha^{\Betti}_{g,\bullet}$ is \'etale locally modelled by the preprojective CoHA for the quiver appearing in Theorem \ref{nbhd2CY}, the claim that $\Up\tau^{\leq 0}\RCoha^{\Betti}_{g,\bullet}$ is a free algebra generated by the subobjects $\ICS_{\Msp_{g,r,0}^{\Betti}}$ follows from the claims that for the quiver $Q$ appearing in Theorem \ref{nbhd2CY}, every vertex supports at least two loops, and any two vertices are joined by an edge.  
\smallbreak
For the first claim, we need to show that for each $\mathcal{F}$ a simple $\BC[\pi_1(\Sigma_g)]$-module, we have $\dim(\Ext^1(\CF,\CF))/2>1$.  By simplicity of $\CF$, and the 2-Calabi--Yau property, we have
\[
\dim(\Ext^2(\CF,\CF))=\dim(\Ext^0(\CF,\CF))=1
\]
and we have
\[
\chi(\CF,\CF)\coloneqq \sum_{i\in\BZ}(-1)^i\dim(\Ext^i(\CF,\CF))=2r^2(1-g)
\]
where $r=\dim(\CF)$.  So the required inequality follows from $g\geq 2$.  For the second claim, we note that if $\CF'$ and $\CF''$ are distinct simple objects of ranks $r'$ and $r''$, then $\Ext^0(\CF',\CF'')=\Ext^2(\CF',\CF'')=0$, so that 
\[
\chi(\CF',\CF'')=-\dim(\Ext^1(\CF',\CF''))=2r'r''(1-g)
\]
as required.
\end{remark}
\begin{remark}
\label{lowRank}
Conjecture \ref{Freeg} and Theorem \ref{CIT} give a precise way to calculate the mixed Hodge structure on $\underline{\ICA}(\Msp^{\Betti}_{g,r,0})$ from the mixed Hodge structures on $\HO^{\BM}(\Mst^{\Betti}_{g,r',0},\underline{\BQ}_{\vir})$ for $r'\leq r$, and similarly, a way to calculate the mixed Hodge structure on $\HO^{\BM}(\Mst^{\Betti}_{g,r,0},\underline{\BQ}_{\vir})$ via the mixed Hodge structures on $\underline{\ICA}(\Msp^{\Betti}_{g,r,0})$ for all $r'\leq r$.  Conversely, by independently calculating these mixed Hodge structures we can check Conjecture \ref{Freeg} in low rank.
\end{remark}
\subsection{$\chi$-independence conjecture}
We recall from \cite{Da16} a separate conjecture regarding the BPS Lie algebra $\underline{\Fg}^{\Betti}_{g,\bullet}$, along with the evidence for it from \cite{HLRV13}.
\smallbreak
Noting that $\Msp^{\Betti}_{g,r,1}$ is smooth of dimension $2r^2(g-1)+2$ if it is not empty, in Conjecture \ref{BCI} we set
\[
\BQ_{\vir}=\BQ_{\Msp^{\Betti}_{g,r,1}}[2r^2(g-1)+2],\quad\quad \ul{\BQ}_{\vir}=\ul{\BQ}_{\Msp^{\Betti}_{g,r,1}}\otimes\LLL^{-r^2(g-1)-1},
\]
the constant perverse sheaf and mixed Hodge module respectively.
\begin{conjecture}[Betti $\chi$-independence]
\label{BCI}
For all $g\in\BZ_{\geq 0}$ and $r\in\BZ_{> 0}$ there is a canonical isomorphism of cohomologically graded vector spaces
\begin{align}
\label{BCICI}
\Fg^{\Betti}_{g,r}\cong\HO(\Msp^{\Betti}_{g,r,1},\BQ_{\vir})
\end{align}
lifting to a canonical isomorphism of mixed Hodge structures
\[
\underline{\Fg}^{\Betti}_{g,r}\cong\HO(\Msp^{\Betti}_{g,r,1},\underline{\BQ}_{\vir}).
\]
\end{conjecture}
While this paper was being revised, the first part of this conjecture was proved in \cite{DHSM12}.  

We call this conjecture the $\chi$-independence conjecture by analogy with the Dolbeault side, where it corresponds to invariance of BPS cohomology under change of the degree of the underlying Higgs bundles under consideration, with rank fixed (equivalently, change of the Euler character, with rank fixed).  

We recall the evidence for the second part of Conjecture \ref{BCI} from \cite{HLRV13}.  Firstly note that Conjecture \ref{BCI} predicts, via \eqref{JAA}, that there is an isomorphism of $\BN$-graded, cohomologically graded mixed Hodge structures
\begin{align}\label{Echeck}
\underline{\Coha}^{\Betti}_{g,\bullet}\coloneqq &\bigoplus_{r\geq 0}\HO^{\BM}(\Mst^{\Betti}_{g,r,0},\ul{\BQ}_{\vir})\\ 
\label{Echeck}
\cong&\Sym\left(\bigoplus_{r\geq 1}\HO^{\BM}(\Msp^{\Betti}_{g,r,1},\ul{\BQ}_{\vir})\otimes\HO(\pt/\BC^*,\underline{\BQ})\right).
\end{align}
\smallbreak
Let $\mathcal{H}$ be a cohomologically graded mixed Hodge structure satisfying the condition
\begin{itemize}
\item[(*)] \label{lfit}$\Ho^i=0$ for $i\ll 0$, and $\Gr^W_j\!\Ho^i=0$ for $j<i$.  
\end{itemize}
The mixed Hodge structure on the Borel--Moore homology of a finite type global quotient stack satisfies condition (*).  We define the \textit{E-series}
\[
\EP(\CH)\coloneqq \sum_{i\in\BZ}(-1)^i\sum_{p,q\in\BZ}\dim(\Gr_p^F\Gr^W_{p+q}\CH^i)x^py^q.
\]
Our conditions on $\mathcal{H}$ ensure that this expression makes sense, i.e. for each $p,q\in\BZ$ the sum $\sum_{i\in\BZ} \dim(\Gr^F_p\Gr^W_{p+q}\mathcal{H}^i)$ is finite.  If $\Ho_{\bullet}$ is a $\BN$-graded, cohomologically graded mixed Hodge structure, where each $\BN$-graded piece satisfies (*), we define the E-series by
\[
\EP_t(\CH_{\bullet})\coloneqq\sum_{n\in\BZ}\EP(\CH_n)t^n.
\]
If $\CH_{n}$ vanishes for $n\leq 0$, and each $\CH_{n}$ satisfies condition (*), then each $\mathbb{N}$-graded piece of $\Sym(\CH_{\bullet})$ again satisfies (*), and
\[
\EP_t(\Sym(\CH_{\bullet}))=\EXP(\EP_t(\CH_{\bullet}))
\]
where on the RHS we take the plethystic exponential with respect to the variables $x,y,t$.  Explicitly, this operation is defined by
\[
\EXP\left(\sum_{i,j,k\in \BZ}a_{i,j,k}x^iy^jt^k\right)=\prod_{i,j,k\in\BZ}(1-x^iy^jt^k)^{-a_{i,j,k}}.
\]
In particular, the $E$-series of the RHS of \eqref{Echeck} is determined by the $E$-polynomials of $\HO(\Msp^{\Betti}_{g,r,1},\ul{\BQ})$.  In \cite{HRV08} the $E$-series of the mixed Hodge structures $\HO(\Msp^{\Betti}_{g,r,1},\ul{\BQ})$ and $\HO^{\BM}(\Mst^{\Betti}_{g,r,0},\ul{\BQ})$ are calculated, thereby determining the $E$-series of the LHS and RHS of \eqref{Echeck}, which agree.  This provides strong evidence in support of Conjecture \ref{BCI} (see \cite{Da16} for further discussion).
\subsection{Genus $\leq 1$}
\label{BG1}
We finish this section by proving the analogues of Conjectures \ref{Freeg} and \ref{BCI} in genus zero and one.  In genus zero, $\pi_1(\Sigma_g)$ is trivial.  There is a unique simple $\BC[\pi_1(\Sigma_g)]$-module, up to isomorphism: a one-dimensional vector space.  So $\Mst^{\Betti}_{0,r,0}\cong \pt/\Gl_r$, and $\Msp^{\Betti}_{0,r,0}\cong \pt$.  As per Remarks \ref{trich_rema} and \ref{trich_rem}, since now for the (unique) simple object $\CF$ there is an equality $\dim(\Ext^1(\CF,\CF))=0$, we expect the Lie algebra $\ul{\Fg}^{\Betti}_{0,\bullet}$ to be a Kac--Moody Lie algebra, generated by 
\[
\ul{\ICA}_{\Msp_{1}(\BC[\pi_1(\Sigma_0)])}\cong \ul{\BQ}.
\]
On the right hand side, $\ul{\BQ}$ is the weight zero one-dimensional mixed Hodge structure, placed in $\BN$-degree $1$, and on the left hand side we only take contributions corresponding to simple objects.  We have isomorphisms
\begin{align*}
\BD\ul{\BQ}_{\Mst^{\Betti}_{0,r,0},\vir}\cong &\ul{\BQ}_{\Mst^{\Betti}_{0,r,0}}\\
\JH_*\BD\ul{\BQ}_{\Mst^{\Betti}_{0,r,0},\vir}\cong& \HO(\pt/\Gl_r,\ul{\BQ})\otimes\ul{\BQ}_{\Msp_{0,r,0}^{\Betti}}.
\end{align*}
In particular, we find an isomorphism of pure Hodge structures
\begin{align*}
\FL^0 \!\Coha^{\Betti}_g\cong&\bigoplus_{r\geq 0}\HO^0(\pt/\Gl_r,\ul{\BQ})\\
\cong \ul{\BQ}[x]
\end{align*}
with $x$ in $\BN$-degree one.  By Proposition \ref{UEAB} there is an isomorphism of algebras 
\[
\FL^0 \!\Coha^{\Betti}_{0,\bullet}\cong \UEA\!\left(\Fg^{\Betti}_{0,\bullet}\right)
\]
and so by dimension counting we find that $\ul{\Fg}^{\Betti}_{0,\bullet}$ is a one-dimensional Lie algebra concentrated in $\BN$-degree one, and moreover $\ul{\Fg}^{\Betti}_{0,1}\cong\ul{\BQ}$.  This confirms the genus zero version of Conjecture \ref{Freeg}.  For Conjecture \ref{BCI}, it is enough to observe that
\[
\Msp^{\Betti}_{0,r,1}\cong \begin{cases} \pt& \textrm{if }r=1\\
\emptyset &\textrm{for }r\geq 2.\end{cases}
\]
\smallbreak
Next we turn to genus one.  Now we calculate that $\chi(\CF',\CF'')=0$, for all $\CF,\CF''$, so that (examining the Serre relations \eqref{SerreRel}) we expect the BPS Lie algebra to be Abelian, and furthermore we expect from Conjecture \ref{BorC2} that there is an isomorphism of Lie algebras
\begin{equation}
\label{gen1g}
\ul{\Fg}^{\Betti}_{1,\bullet}\cong \bigoplus_{r\geq 1}\ul{\ICA}(\Msp_{1}(\BC[\pi_1(\Sigma_1)]))
\end{equation}
where the right hand side is given the zero bracket.  By \cite[Sec.5.2]{preproj} there are isomorphisms
\[
\ul{\DTS}^{\Betti}_{1,r}\cong\Delta_{r,*}\ul{\BQ}_{(\BC^*)^2,\vir}
\]
where $\Delta_r\colon (\BC^*)^2\rightarrow \Sym^r((\BC^*)^2)\cong \Msp_r(\BC[\pi_1(\Sigma_1)])$ is the inclusion of the small diagonal, and 
\[
\ul{\BQ}_{(\BC^*)^2,\vir}\coloneqq \ul{\BQ}_{(\BC^*)^2}\otimes\LLL^{-1}.
\]
So by \eqref{JAR2} there is a PBW isomorphism
\begin{equation}
\label{relBPBW}
\JH_*\BD\ul{\BQ}_{\Mst_{\bullet}(\BC[\pi_1(\Sigma_1)]),\vir}\cong \Sym_{\oplus}\left(\bigoplus_{r\geq 1}\Delta_{r,*}\ul{\BQ}_{(\BC^*)^2,\vir}\otimes \HO(\pt/\BC^*,\ul{\BQ})\right)
\end{equation}
from which it follows that
\[
\tau^{\leq 0}\!\JH_*\BD\ul{\BQ}_{\Mst_{\bullet}(\BC[\pi_1(\Sigma_1)]),\vir}\cong \Sym_{\oplus}\left(\bigoplus_{r\geq 1}\Delta_{r,*}\ul{\BQ}_{(\BC^*)^2,\vir}\right)
\]
and so\footnote{In particular, the BPS invariants are independent of $r$, hence the non-appearance of $r$ inside the direct sum on the right hand side of the following isomorphism.} since there is an isomorphism $\Msp^{\Betti}_{1,1,0}\cong (\mathbb{C}^*)^2$ we obtain the isomorphism
\[
\FL^0\!\ul{\Coha}^{\Betti}_{0,\bullet}\cong \Sym\left(\bigoplus_{r\geq 1}\HO(\Msp^{\Betti}_{1,1,0},\ul{\BQ}_{\vir})\right).
\]
So the isomorphism \eqref{gen1g} exists, at least at the level of underlying mixed Hodge structures.  What remains is to show that the Lie bracket on $\Fg_{1,\bullet}^{\Betti}$ vanishes.  This follows by the arguments of \cite{preproj3}, we provide the details for completeness.  The Lie bracket 
\[
\Fg^{\Betti}_{1,n'}\otimes \Fg^{\Betti}_{1,n''}\rightarrow \Fg^{\Betti}_{1,n'+n''}
\]
is obtained by taking hypercohomology of a morphism
\begin{equation}
\label{LAL}
\oplus_{*}\!\left(\DTS^{\Betti}_{1,n'}\boxtimes \DTS^{\Betti}_{1,n''}\right)\rightarrow \DTS^{\Betti}_{1,n'+n''}.
\end{equation}
Each $\DTS^{\Betti}_{1,n'}$ is a simple perverse sheaf, while the morphism $\oplus$ is finite, so that the only support of the LHS of \eqref{LAL} is given by 
\[
\oplus(\supp(\DTS^{\Betti}_{1,n'})\times \supp(\DTS^{\Betti}_{1,n''})).
\]
In particular, the supports of the LHS and RHS of \eqref{LAL} are different, and so the morphism \eqref{LAL} is zero.  This completes the proof of the $g=1$ version of Conjecture \ref{Freeg}.  
\smallbreak
For Conjecture \ref{BCI}, we recall that for all $r\geq 1$ there is an isomorphism
\begin{align*}
(\BC^*)^{\times 2}&\rightarrow \Msp^{\Betti}_{1,r,1}\\
(z_1,z_2)&\mapsto \left(z_1\cdot \left(\begin{matrix} \zeta_r&0&\ldots &\\ 0&\zeta_r^2&0&\ldots\\ \ldots &&&\\ 0&\ldots&0&1\end{matrix}\right),z_2\cdot C\right)
\end{align*}
where $\zeta_r$ is a primitive $r$th root of unity, and $C$ is the cyclic permutation matrix.
So there are canonical isomorphisms
\[
\Msp^{\Betti}_{1,r,1}\cong (\BC^*)^{\times 2},
\]
and
\[
\ul{\Fg}^{\Betti}_{1,r}\cong \HO(\Msp^{\Betti}_{1,r,1},\ul{\BQ}_{\vir})
\]
as required (comparing with \eqref{gen1g}).
\section{The Dolbeault side}
\subsection{BPS algebra}
In analogy with the Betti side, we denote by $\JH\colon \Mst_{r,d}^{\Dol}(C)\rightarrow \Msp_{r,d}^{\Dol}(C)$ the morphism from the stack of degree $d$ rank $r$ semistable Higgs bundles to its coarse moduli space.  Then by Theorem \ref{p2CY} there is an isomorphism in $\Dub^+(\Perv(\Msp^{\Dol}_{r,d}(C)))$
\[
\JH_*\BD \mathbb{Q}_{\Mst_{r,d}^{\Dol}(C),\vir}\cong \bigoplus_{i\in 2\cdot \BZ_{\geq 0}}\!\!\Up\Ho^i(\JH_*\BD \mathbb{Q}_{\Mst_{r,d}^{\Dol}(C),\vir})[-i]
\]
and each
\begin{equation}
\label{HPU}
\Up\Ho^i\!\left(\JH_*\BD\mathbb{Q}_{\Mst_{r,d}^{\Dol}(C),\vir}\right)
\end{equation}
is a semisimple perverse sheaf.  Moreover, the lift of \eqref{HPU} is a pure weight $i$ mixed Hodge module.  The following is a corollary of these results.
\begin{theorem}[\cite{Da21a}; Sec.7.4]\label{absPt}
The complex of mixed Hodge modules $\Hit^{\Sta}_*\BD\ul{\mathbb{Q}}_{\Mst_{\bullet,0}^{\Dol}(C),\vir}$ is pure, as is the cohomologically graded mixed Hodge structure $\HO^{\BM}(\Mst_{\bullet,0}^{\Dol}(C),\ul{\mathbb{Q}}_{\vir})$.
\end{theorem}

The relative integrality conjecture for degree zero Higgs bundles, recently proved by Tasuki Kinjo and Naoki Koseki \cite{KiKo21} states that there is an isomorphism 
\begin{equation}
\label{HIC1}
\JH_*\BD \ul{\mathbb{Q}}_{\Mst_{\bullet,0}^{\Dol}(C),\vir}\cong \Sym_{\oplus}\left( \bigoplus_{r\geq 1}\ul{\DTS}_r^{\Dol}(C)\otimes \HO(\pt/\mathbb{C}^*,\BQ) \right)
\end{equation}
for certain mixed Hodge modules $\ul{\DTS}^{\Dol}_r(C)\in\MHM(\Msp_{r,0}^{\Dol}(C))$.  
\begin{remark}
\label{DolUn}
By Theorem \ref{p2CY}, we may consider the left hand side of \eqref{HIC1} as an object in a semisimple Abelian category, namely the category of cohomologically graded pure mixed Hodge modules.  It follows, formally, that if there are two isomorphisms
\begin{align*}
\JH_*\BD \ul{\mathbb{Q}}_{\Mst_{\bullet,0}^{\Dol}(C),\vir}\cong\Sym_{\oplus}(\mathcal{V}\otimes \HO(\pt/\mathbb{C}^*,\BQ)_{\vir})\\
\JH_*\BD \ul{\mathbb{Q}}_{\Mst_{\bullet,0}^{\Dol}(C),\vir}\cong\Sym_{\oplus}(\mathcal{V}'\otimes \HO(\pt/\mathbb{C}^*,\BQ)_{\vir})
\end{align*}
then $\mathcal{V}$ and $\mathcal{V}'$ are isomorphic, so that, up to isomorphism, $\ul{\DTS}_r^{\Dol}(C)$ is uniquely defined as the object fitting into an isomorphism \eqref{HIC1}.
\end{remark}

The (absolute) integrality theorem, obtained by applying the derived global sections functor to \eqref{HIC1}, states that there is a $\mathbb{N}$-graded isomorphism of cohomologically graded mixed Hodge structures.
\begin{equation}
\label{HIC2}
\HO^{\BM}(\Mst_{\bullet,0}^{\Dol}(C),\ul{\mathbb{Q}}_{\vir})\cong \Sym\left( \bigoplus_{r\geq 1}\ul{\DT}^{\Dol}_r(C)\otimes \HO(\pt/\mathbb{C}^*,\underline{\BQ})\right)
\end{equation}
for a collection of cohomologically graded mixed Hodge structures $\ul{\DT}^{\Dol}_r(C)$, with finite-dimensional total cohomology.  
\smallbreak
We denote by $\Hit_{r,0}\colon\Msp^{\Dol}_{r,0}(C)\rightarrow \Lambda_{r}$ the morphism to the Hitchin base, where
\[
\Lambda_r=\prod_{i=1}^r\HO^0(C,\omega_C^{\otimes i})
\]
and $\Hit_{r,0}((\CF,\eta\colon \CF\rightarrow \omega_C\otimes\CF))=(\Tr(\eta),\ldots,\Tr(\wedge^r\eta))$.  We define
\[
\Hit^{\Sta}_{r,0}\colon\Mst^{\Dol}_{r,0}(C)\rightarrow \Lambda_{r}
\]
by setting $\Hit^{\Sta}_{r,0}=\Hit_{r,0}\circ \JH$.  
\smallbreak
We give $\Lambda_{\bullet}$ the structure of a commutative monoid as follows.  By the theory of symmetric functions, for each $l\in\mathbb{Z}_{\geq 1}$ there are universal polynomials $P_l(X_1,\ldots,X_l)$ and $E_l(X_1,\ldots,X_l)$ such that if $G\in\Mat_{r\times r}(\BC)$ then $\Tr(G^l)=P_l(\Tr(G),\ldots,\Tr(\wedge^l G))$ and $\Tr(\wedge^l G)=E(\Tr(G),\ldots,\Tr(G^l))$, and we define the isomorphism
\begin{align*}
\mathbf{P}_r\colon &\Lambda_r\rightarrow \Lambda_r\\
&(s_1,\ldots,s_r)\mapsto (P_1(s_1),P_2(s_1,s_2),\ldots,P_r(s_1,\ldots,s_r)).
\end{align*}
In words, $\mathbf{P}_r$ replaces elementary symmetric functions in the eigenvalues of the Higgs field with power sums.  Fix $m,n$.  For each $j\geq 1$ there is a universal polynomial $A_j(X_1,\ldots,X_m)\in\BQ[X_1,\ldots,X_m]$ such that if $G$ is a $m\times m$ matrix, there is an equality $\Tr(G^j)=A_j(\Tr(G),\ldots,\Tr(G^m))$.  We define the embedding $\iota_m$
\begin{align*}
\iota_m\colon &\Lambda_m\hookrightarrow \Lambda_{m+n}\\
&(s_1,\ldots, s_m)\mapsto (A_1(s_1,\ldots,s_m),\ldots,A_{m+n}(s_1,\ldots,s_m)).
\end{align*}
Note that $A_i(s_1,\ldots,s_n)=s_i$ for $i\leq m$, so $\iota_m$ is indeed an embedding.  We define an embedding $\iota_n\colon \Lambda_n\hookrightarrow \Lambda_{m+n}$ similarly.  Then we define the morphism
\begin{align}
\label{sc_def}
\Lambda_m\times\Lambda_n&\rightarrow \Lambda_{m+n}\\ \nonumber
(s,s')&\mapsto \mathbf{P}^{-1}_{m+n}(\iota_{m}\mathbf{P}_m(s)+\iota_n\mathbf{P}_n(s')).
\end{align}
Clearly this morphism is symmetric, in the sense that the morphism
\[
\cup\colon \Lambda_{\bullet}\times\Lambda_{\bullet}\rightarrow\Lambda_{\bullet}
\]
defined by taking the union over all morphisms \eqref{sc_def} is symmetric in the two factors of $\Lambda_{\bullet}$.
\begin{proposition}\label{CupFin}
The morphism $\cup$ is finite.
\end{proposition}
\begin{proof}
Since $\cup$ is a morphism of affine varieties, finiteness follows from properness, which we now prove.
\smallbreak
For fixed $r\in\BZ_{\geq 1}$, pick $N$ distinct points $x_1,\ldots,x_N\in C$.  We denote by $i_t\colon x_t\hookrightarrow C$ the inclusion, and set 
\[
\Lambda^{\circ}_{r,N}=\prod_{t=1}^N \prod_{j=1}^r i_{t}^*\omega_{C}^{\otimes j}.
\]
We define the morphism of vector spaces
\begin{align*}
&\Lambda_r\rightarrow \Lambda^{\circ}_{r,N}\\
&(s_1,\ldots,s_r)\mapsto (s_1\lvert_{x_1},s_2\lvert_{x_1},\ldots,s_r\lvert_{x_1}, s_1\lvert_{x_2},\ldots, s_r\lvert_{x_N}).
\end{align*}
For $N\gg 0$ this is an embedding of vector spaces, and so it is closed.  Picking trivializations of $\omega_C$ at each point $x_t$ defines an isomorphism
\[
\Lambda^{\circ}_{r,N}\cong \Sym^r(\mathbb{A}^1)^{\times N},
\]
and a commutative diagram
\[
\xymatrix{
\Lambda_m\times\Lambda_n\ar@{^{(}->}[d]\ar[r]^-\cup&\Lambda_{m+n}\ar@{^{(}->}[d]\\
\Sym^{m}(\BA^1)^{\times N}\times \Sym^n(\BA^1)^{\times N}\ar[r]^-{\cup^{\times N}}&\Sym^{m+n}(\BA^1)^{\times N}
}
\]
in which the vertical maps are closed embeddings, and we wish to show that the top horizontal morphism is proper.  This follows from properness of
\[
\Sym^{m}(\BA^1)\times \Sym^n(\BA^1)\xrightarrow{\cup}\Sym^{m+n}(\BA^1),
\]
which is an easy exercise (as well as being the special case of \cite[Lem.2.1]{Meinhardt14} for the one loop quiver).
\end{proof}
Applying $\Hit^{\Sta}_*$, instead of $\HO(-)$, to \eqref{HIC1}, yields an isomorphism of cohomologically graded semisimple mixed Hodge modules on the Hitchin base $\Lambda_{\bullet}$
\begin{equation}
\label{HIC3}
\Hit^{\Sta}_{\bullet,0,*}\BD\ul{\mathbb{Q}}_{\Mst_{\bullet,0}^{\Dol}(C),\vir}\cong \Sym_{\cup}\!\left( \bigoplus_{r\geq 1}\Hit_{r,0,*}\ul{\DTS}_r^{\Dol}(C)\otimes \HO(\pt/\mathbb{C}^*,\BQ)_{\vir} \right).
\end{equation}
The above forms of the integrality theorem are linked in a very natural way, which we recall:
\begin{lemma}\cite[Cor.5.20]{KiKo21}
The existence of the isomorphism \eqref{HIC1} implies the existence of the isomorphisms \eqref{HIC2} and \eqref{HIC3}.
\end{lemma}
\begin{proof}
Since we obtain \eqref{HIC2} from \eqref{HIC1} by taking derived global sections, the claim regarding \eqref{HIC2} is obvious.  Next, \eqref{HIC3} is obtained by applying $\Hit_{\bullet,0,*}$ to \eqref{HIC1}, and the claim is then a formal consequence of the claim that $\Hit_{\bullet,0,*}$ respects the symmetric monoidal structures on $\Dub^+(\MHM(\Msp_{\bullet,0}^{\Dol}(C)))$ and $\Dub^+(\MHM(\Lambda_{\bullet}))$.  This in turn follows from commutativity of the diagram
\[
\xymatrix{
\Msp^{\Dol}_{m,0}(C)\times\Msp^{\Dol}_{n,0}(C)\ar[r]^-\oplus\ar[d]^{\Hit_{m,0}\times\Hit_{n,0}}&\Msp^{\Dol}_{m+n,0}(C)\ar[d]^{\Hit_{m+n,0}}\\
\Lambda_m\times\Lambda_n\ar[r]^-\cup&\Lambda_{m+n}.
}
\]
\end{proof}
Set
\begin{align*}
\Coha^{\Dol}_{\bullet}(C)\coloneqq &\bigoplus_{r\in\BZ_{\geq 0}}\HO^{\BM}(\Mst^{\Dol}_{r,0}(C),\BQ_{\vir})\\
\RCoha^{\Dol}_{\bullet}(C)\coloneqq &\JH_*\BQ_{\Mst^{\Dol}_{\bullet,0}(C),\vir}.
\end{align*}
In \cite{SS20} (following \cite{Mi20} for the case of torsion Higgs sheaves) a Hall algebra structure is defined on $\Coha^{\Dol}_{\bullet}(C)$.  As in the definitions of $\Coha_{B_g,W,\bullet}$ and $\Coha^{\Betti}_{g,\bullet}$, the multiplication is defined in terms of morphisms of complexes of perverse sheaves that descend to the respective coarse moduli spaces, so that we may lift the CoHA to a Hall algebra structure on $\RCoha_{\bullet}^{\Dol}(C)$, and moreover lift both algebra objects to the categories of mixed Hodge structures, or mixed Hodge modules, respectively.
\smallbreak
We do not yet have a definition of the BPS Lie algebra on the Dolbeault side.  On the other hand, we have an analogue of the universal enveloping algebra of the BPS Lie algebra from the Betti side; by temporary abuse of notation we write
\begin{align*}
\UEA(\Fg^{\Dol}_{\bullet}(C))\coloneqq &\mathfrak{L}^0\!\Coha^{\Dol}_{\bullet}(C)\\
\coloneqq&\HO\!\left(\Msp^{\Dol}_{\bullet,0}(C),\Up\tau^{\leq 0}\JH_*\BD\BQ_{\Mst^{\Dol}_{\bullet,0}(C),\vir}\right).
\end{align*}

For now we set $g(C)\geq 2$.  From the proof of \cite[Thm.7.35]{Da21a} we deduce
\begin{theorem}
There is an inclusion of perverse sheaves 
\[
\ICS_{\Msp^{\Dol}_{r,0}(C)}\hookrightarrow \Up\tau^{\leq 0}\!\RCoha^{\Dol}_r(C)
\]
which is moreover primitive, in the sense that this is the inclusion associated to a (canonical) decomposition
\[
\Up\tau^{\leq 0}\!\RCoha^{\Dol}_r(C)\cong \ICS_{\Msp^{\Dol}_{r,0}(C)}\oplus\mathcal{C}_r
\]
and for $r'+r''=r$ with $r',r''\in\BZ_{>0}$, the multiplication morphism
\[
\Up\tau^{\leq 0}\!\RCoha^{\Dol}_{r'}(C)\boxtimes_{\oplus}\Up\tau^{\leq 0}\!\RCoha^{\Dol}_{r''}(C)\rightarrow \Up\tau^{\leq 0}\!\RCoha^{\Dol}_{r}(C)
\]
factors through the inclusion of $\mathcal{C}_r$.
\end{theorem}
Taking derived global sections, the theorem implies that there are canonical inclusions 
\[
\ICA(\Msp^{\Dol}_{r,0}(C))\hookrightarrow\UEA(\Fg^{\Dol}_{\bullet}(C))
\]
which can be extended to form the inclusion of a minimal generating space of the algebra $\UEA(\Fg^{\Dol}_{\bullet}(C))$.  The following conjecture states that the images of these inclusions already generate.  This is the Dolbeault version of Conjecture \ref{freeA}, and is motivated the same way, via Conjecture \ref{BorC2}.  As in the case of Conjecture \ref{freeA}, it was proved during revision of this paper in \cite{DHSM12}.
\begin{conjecture}
\label{HGC}
If $g(C)\geq 2$ then the algebra $\UEA(\Fg^{\Dol}_{\bullet}(C))$ is the free associative algebra object generated by the subspaces $\bigoplus_{r\geq 1}\ICA(\Msp^{\Dol}_{r,0}(C))$.    Similarly, the algebra object $\Up\:\tau^{\leq 0}\RCoha^{\Dol}_{\bullet}(C)$ is the free associative algebra object generated by $\bigoplus_{r\geq 1}\ICS_{\Msp^{\Dol}_{r,0}(C)}$.  Moreover these statements lift to (pure) mixed Hodge structures/modules, respectively.
\end{conjecture}

\subsection{$\chi$-independence}
The following result was very recently proved by Naoki Koseki and Tasuki Kinjo.  The proof builds, in particular, on the $\chi$-independence result of \cite{MaSh20}, the dimensional reduction isomorphism of \cite{Kin20}, and the interpretation of deformed Calabi--Yau completions of total spaces of line bundles of \cite{KiMa21}.  It is a cohomological version of the $\chi$-independence conjecture for Gopakumar--Vafa invariants of Calabi--Yau threefolds (see \cite{MaTo18,To17} for details), and is indeed obtained as a special case of the $\chi$-independence result for general local curves in Calabi--Yau threefolds, also proved in \cite{KiKo21}.
\begin{theorem}\cite{KiKo21}
\label{HiggsChi}
Let $d\in\BZ$ satisfy $(r,d)=1$.  Then there is a canonical isomorphism of cohomologically graded mixed Hodge modules
\begin{equation}
\label{chi_ind}
\Hit_{r,d,*}\ul{\BQ}_{\Msp^{\Dol}_{r,d}(C),\vir}\cong \Hit_{r,0,*}\ul{\DTS}_{r}^{\Dol}(C)
\end{equation}
\end{theorem}
Setting $d=1$ in Theorem \ref{HiggsChi} and taking derived global sections produces
\begin{theorem}\cite{KiKo21}
\label{ACCCon}
There is an isomorphism of cohomologically graded mixed Hodge structures
\begin{equation}
\label{ACCC}
\Coha^{\Dol}_{\bullet}(C)\cong \Sym\left(\HO(\Msp^{\Dol}_{\bullet,1}(C),\ul{\BQ}_{\vir})\otimes\HO(\pt/\BC^*)\right)
\end{equation}
which moreover respects the perverse filtrations on both sides.
\end{theorem}

The purity statements of Theorem \ref{p2CY} and \ref{absPt} imply the following
\begin{corollary}
The complexes of mixed Hodge modules $\ul{\DTS}_r^{\Dol}(C)$ and $\Hit_{r,0,*}\ul{\DTS}_r^{\Dol}(C)$ are pure, as is the cohomologically graded mixed Hodge structure $\ul{\DT}_r^{\Dol}(C)$.
\end{corollary}


\begin{remark}
Similarly to Remark \ref{lowRank}, Conjecture \ref{HGC}, along with Theorem \ref{ACCCon} and the integrality isomorphism \eqref{HIC2} give a precise way of calculating one out of $\ICA(\Msp^{\Dol}_{r,0}(C))$, $\HO(\Msp^{\Dol}_{r,1}(C),\BQ_{\vir})$ or $\HO^{\BM}(\Mst^{\Dol}_{r',0}(C),\BQ_{\vir})$, for fixed $r$, by calculating one of the others for all $r'\leq r$.  Again, by independently calculating these three types of cohomology, we can begin to confirm Conjecture \ref{HGC} in sufficiently low rank; see \cite{SM21}, building on \cite{Fe21,Ma21} for results in this direction.  In particular, in \cite{SM21} Schlegel-Mejia proves that Conjecture \ref{HGC}, and the integrality isomorphism \eqref{HIC2} together give the right prediction for the graded dimensions of $\HO^{\BM}(\Mst^{\Dol}_{r,0}(C),\BQ_{\vir})$ up to $r=2$, and for all $g\geq 2$, and thus confirmation of Conjecture \ref{HGC} for $r\leq 2$.
\end{remark}
In the remainder of this section we check appropriate modifications of Conjecture \ref{HGC} for the cases $g(C)\leq 1$.  


\subsection{Genus zero}
First we deal with the case $g(C)=0$, i.e. we set $C=\mathbb{P}^1$.  By Grothendieck's theorem, a rank $r$ vector bundle $V$ on $C$ splits as a direct sum
\[
V=\bigoplus_{1\leq p\leq r}\mathscr{O}_{\BP^1}(a_p)
\]
for some integers $a_p$.  If $V$ has degree zero then $\sum_{p=1}^r a_p=0$.  Furthermore 
\[
\Hom(\mathscr{O}_{\BP^1}(a_{p'}),\mathscr{O}_{\BP^1}(a_{p''})\otimes\omega_{\BP^1})=0
\]
unless $a_{p''}\geq a_{p'}+2$.  Putting these facts together, the only degree zero rank $r$ semistable Higgs bundle is $\mathscr{O}_{\mathbb{P}^1}^{\oplus r}$, with zero Higgs field, and so 
\begin{align*}
\Mst^{\Dol}_{r,0}(C)\cong &\pt/\Gl_r\\
\Msp^{\Dol}_{r,0}(C)\cong &\pt.
\end{align*}
Similarly, for $(r,d)=1$, we deduce that there is a semistable Higgs bundle of rank $r$ and degree $d$ if and only if $r=1$, and for $r=1$, $(\mathscr{O}_C(d),0)$ is the unique such Higgs bundle.  Thus, for $(r,d)=1$
\begin{equation}
\label{rdcde}
\Msp^{\Dol}_{r,d}(C)=\begin{cases} \pt& \textrm{if }r=1\\
\emptyset & \textrm{if }r\geq 2.
\end{cases}
\end{equation}
Arguing as in the Betti case, we find 
\begin{align*}
\JH_*\BD\ul{\BQ}_{\Mst^{\Dol}_{r,0}(C),\vir}\cong &\HO(\pt/\Gl_r,\ul{\BQ})\otimes \ul{\BQ}_{\Msp^{\Dol}_{r,0}(C)}\\
\tau^{\leq 0}\JH_*\BD\ul{\BQ}_{\Mst^{\Dol}_{r,0}(C),\vir}\cong &\ul{\BQ}_{\Msp^{\Dol}_{r,0}(C)}
\end{align*}
and so there is an isomorphism of $\BZ_{\geq 0}$-graded vector spaces
\begin{equation}
\label{teffi}
\FL^0\!\Coha^{\Dol}_{0,\bullet}(C)\cong\BC[x].
\end{equation}
Since we have not actually shown that the LHS of \eqref{teffi} is a $\mathbb{N}^{Q_0}$-graded Borcherds algebra, this is not enough to conclude that \eqref{teffi} is an isomorphism of algebras.  On the other hand, it follows by construction \cite{SS20} that the algebra $\Coha^{\Dol}_{0,\bullet}(C)$ is isomorphic to the CoHA $\Coha_{\BC Q,0,\bullet}$, for $Q$ the quiver with one vertex and one loop, which is commutative.  In fact this algebra admits an explicit description, given in \cite{KS2}.  We have $\Coha_{\BC Q,0,d}\cong \mathbb{Q}[x_1,\ldots,x_d]^{\mathfrak{S}_d}$, the ring of symmetric polynomials in $d$ variables, with product given by the (commutative) shuffle product; see \cite[Sec.2.5]{KS2} for details.  This implies that \eqref{teffi} is an isomorphism of algebras, proving the genus zero version of Conjecture \ref{HGC} \footnote{Moreover, this yields the cohomological integrality conjecture for $\Mst_{\bullet,0}^{\Dol}(C)$.}.  The genus zero case of Theorems \ref{HiggsChi} and \ref{ACCCon}) are confirmed explicitly by \eqref{rdcde}.
\subsection{Genus one}
We first recall some results on Fourier-Mukai functors for elliptic curves (see \cite{Po03} for full details).  Let $C$ be a genus one smooth projective curve.  We denote by $\hat{C}$ the dual Abelian variety, i.e. $\hat{C}=\Pic_0(C)$.  The variety $\hat{C}$ is another elliptic curve, isomorphic to $C$.  Since $\hat{C}$ represents the $\Pic_0$-functor, there is a universal bundle $\mathcal{P}\in\Coh(C\times \hat{C})$ such that $\mathcal{P}\lvert_{C\times \{x\}}\cong \mathcal{L}_x$, where $\mathcal{L}_x$ is the degree zero line bundle corresponding to $x\in\hat{C}$.  The Fourier-Mukai functor with kernel $\mathcal{P}$, which is defined to be
\[
\FoM_{\mathcal{P}}=\pi_{2,*}\circ(-\otimes \mathcal{P})\circ\pi_1^*\colon\Db(\Coh(C))\rightarrow \Db(\Coh(\hat{C}))
\]
is an equivalence of categories.  An inverse equivalence is provided by the Fourier-Mukai kernel $\mathcal{Q}=(\iota\times\id_{\hat{C}})^*\mathcal{P}[1]$, where $\iota$ is the group inverse operation.  
\smallbreak
Given a variety $X$ we denote by $\Mst_n(X)$ the stack of length $n$ coherent sheaves on $X$.  The Fourier--Mukai functor $\FoM_{\mathcal{Q}}$ induces a morphism of stacks
\[
\FFM_{\mathcal{Q}}\colon \Mst_{n}(\hat{C})\rightarrow \mathfrak{C}oh(C)
\]
in the standard way.  Explicitly, to an object $\CF\in\Coh(T\times \hat{C})$, flat over $T$, we assign the object 
\[
\pi_{1,3,*}((\mathscr{O}_T\boxtimes\mathcal{Q})\otimes \pi_{1,2}^*\CF)
\]
where $\pi_{i,j}\colon T\times \hat{C}\times C$ is the projection onto the $i$th and $j$th factors.  The claim that the resulting complex of coherent sheaves is in fact a coherent sheaf follows as in \cite{Bri99}.  The functor $\FoM_{\mathcal{Q}}$ maps structure sheaves of points to degree zero line bundles on $C$, and so since all coherent sheaves on $\hat{C}$ with zero-dimensional can be written as extensions of structure sheaves on points, $\FoM_{\mathcal{Q}}$ maps degree zero coherent sheaves to semistable degree zero coherent sheaves on $C$.  By \cite{Po03}, this morphism is in fact an isomorphism at the level of points, so that (again by arguments of \cite{Bri99}) $\FFM_{\mathcal{Q}}$ induces the isomorphism of stacks
\[
\FFM_{\mathcal{Q}}\colon \FM_n(\hat{C})\rightarrow \mathfrak{C}oh^{\sstab}_{n,0}(C),
\]
with target the stack of degree zero, rank $n$ semistable coherent sheaves on $C$.  Considering instead the Fourier--Mukai kernel $\mathcal{Q}^+=\mathcal{Q}\boxtimes\mathscr{O}_{\Delta_{\BA^1}}\in\Coh(\hat{C}\times C\times \BA^1\times \BA^1)$, the same arguments yield that the induced morphism of stacks
\[
\FFM_{\mathcal{Q}^+}\colon \FM_r(\hat{C}\times\BA^1)\rightarrow \Mst^{\Dol}_{r,0}(C)
\]
is an equivalence of stacks.  Matching these moduli stacks with their respective coarse moduli spaces yields the commutative diagram
\begin{equation}
\label{FM_diag}
\xymatrix{
\FM_r(\hat{C}\times\BA^1)\ar[d]^{p_r}\ar[r]^-{\FFM_{\mathcal{Q}^+}}&\Mst^{\Dol}_{r,0}(C)\ar[d]^{\JH}\\
\Sym^r(\hat{C}\times\BA^1)\ar[r]^-{\beta}&\Msp^{\Dol}_{r,0}(C)
}
\end{equation}
where $\beta$ exists by universality of $p_r$, and is an isomorphism since $\FFM_{\mathcal{Q}^+}$ is.  
\smallbreak
In conclusion, in order to calculate $\JH_*\underline{\BQ}_{\Mst^{\Dol}_{r,0}(C)}$ it is enough to calculate $p_{r,*}\underline{\BQ}_{\FM_r(\hat{C}\times\BA^1)}$.  We calculate this derived direct image for all surfaces (not just $\hat{C}\times\mathbb{A}^1$) in the appendix (\S \ref{DZ_sec}).  From Proposition \ref{KVredux} and the diagram \eqref{FM_diag} we deduce
\begin{proposition}
\label{gen1to}
Let
\begin{align*}
\Delta_{r}\colon &\Msp_{1,0}^{\Dol}(C)\rightarrow \Msp_{r,0}^{\Dol}(C)
\\ &\tilde{\CF}\mapsto\tilde{\CF}\oplus\ldots\oplus\tilde{\CF}
\end{align*}
be the inclusion of the small diagonal.  There is an isomorphism of cohomologically graded mixed Hodge modules
\[
\JH_*\BD\ul{\BQ}_{\Mst^{\Dol}_{r,0}(C),\vir}\cong \Sym_{\oplus}\left(\bigoplus_{r\geq 1}\Delta_{r,*}\ul{\BQ}_{\Msp_{1,0}^{\Dol}(C),\vir}\otimes\HO(\pt/\BC^*,\ul{\BQ})\right).
\]
There is moreover an isomorphism of algebra objects
\[
\tau^{\leq 0}\JH_*\BD\ul{\BQ}_{\Mst^{\Dol}_{r,0}(C),\vir}\cong \Sym_{\oplus}\left(\bigoplus_{r\geq 1}\Delta_{r,*}\ul{\BQ}_{\Msp_{1,0}^{\Dol}(C),\vir}\right).
\]
\end{proposition}
Applying the derived global sections functor to the second statement of Proposition \ref{gen1to}, we deduce that there is an isomorphism of algebras
\begin{align*}
\UEA(\Fg^{\Dol}_{\bullet}(C))\coloneqq &\HO\!\left(\Msp^{\Dol}_{\bullet,0}(C),\Up\tau^{\leq 0}\JH_*\BD\BQ_{\Mst^{\Dol}_{\bullet,0}(C),\vir}\right)\\
\cong &\Sym\left(\bigoplus_{r\geq 1} \HO(\Msp^{\Dol}_{1,0}(C),\BQ_{\vir})\right),
\end{align*}
proving the genus one analogue of Conjecture \ref{HGC}.  As on the Betti side, the fact that $r$ does not appear in the sum on the right hand side of this isomorphism tells us that the cohomological BPS invariants for $g=1$ are independent of $r$.

\section{Nonabelian Hodge theory for stacks and the P=W conjectures}
\label{final_sec}
\subsection{Nonabelian Hodge isomorphism}
\label{naht_sec}
We now have all of the ingredients together to explain the extension of nonabelian Hodge theory to stacks, via Hall algebras.
\smallbreak
To start with, we concentrate on the case $g\geq 2$, where things are partly conjectural, before returning to the cases $g\leq 1$, where we prove the analogues of these conjectures.  Recall that for each $r\in\BZ_{\geq 1}$ we have a homeomorphism
\[
\Gamma_{r,0}\colon \Msp^{\Betti}_{g,r,0}\rightarrow \Msp^{\Dol}_{r,0}(C).
\]

By construction, the morphism $\Gamma_{\bullet,0}$ respects the monoidal structures on the domain and target.  Since the intersection complex of a space is invariant under homeomorphism, there is an induced isomorphism
\[
\Gamma_{r,0,*}\ICS_{\Msp^{\Betti}_{g,r,0}}\xrightarrow{\cong}\ICS_{\Msp^{\Dol}_{r,0}(C)}.
\]
Assuming Conjecture \ref{Freeg}, via the (Betti side) integrality isomorphism \eqref{JAR2} there is an isomorphism
\[
p_*\BD\BQ_{\Mst^{\Betti}_{g,\bullet,0},\vir}\cong \Sym_{\oplus}\!\left(\Free_{\Lie}\left(\bigoplus_{r\geq 1} \ICS_{\Msp^{\Betti}_{g,r,0}}\right)\otimes\HO(\pt/\BC^*,\BQ)\right).
\]
Assuming Conjecture \ref{HGC} the (Dolbeault side) integrality isomorphism \eqref{HIC1} yields
\[
p_*\BD\BQ_{\Mst^{\Dol}_{\bullet,0}(C),\vir}\cong \Sym_{\oplus}\!\left(\Free_{\Lie}\left(\bigoplus_{r\geq 1} \ICS_{\Msp^{\Dol}_{r,0}(C)}\right)\otimes\HO(\pt/\BC^*,\BQ)\right).
\]
Putting all of this together, we construct an isomorphism $\Phi_{\bullet}\colon\Gamma_{\bullet,0,*}p_*\BQ_{\Mst^{\Betti}_{g,\bullet,0},\vir}\rightarrow p_*\BQ_{\Mst^{\Dol}_{\bullet,0}(C),\vir}$ by composing the isomorphisms
\begin{align*}
\Gamma_{\bullet,0,*}p_*\BD\BQ_{\Mst^{\Betti}_{g,\bullet,0},\vir}\cong &\Gamma_{\bullet,0,*}\left(\Sym_{\oplus}\!\left(\Free_{\Lie}\left(\bigoplus_{r\geq 1} \ICS_{\Msp^{\Betti}_{g,r,0}}\right)\otimes\HO(\pt/\BC^*,\BQ)\right)\right)
\\\cong & \Sym_{\oplus}\!\left(\Free_{\Lie}\left(\bigoplus_{r\geq 1} \Gamma_{\bullet,0,*}\ICS_{\Msp^{\Betti}_{g,r,0}}\right)\otimes\HO(\pt/\BC^*,\BQ)\right)
\\\cong &\Sym_{\oplus}\!\left(\Free_{\Lie}\left(\bigoplus_{r\geq 1} \ICS_{\Msp^{\Dol}_{r,0}(C)}\right)\otimes\HO(\pt/\BC^*,\BQ)\right)
\\\cong&p_*\BD\BQ_{\Mst^{\Dol}_{\bullet,0}(C),\vir}.
\end{align*}
Taking global sections, we finally arrive at a (canonical) isomorphism:
\begin{equation}
\label{CPD}
\HO(\Phi_{\bullet})\colon\HO^{\BM}(\Mst^{\Betti}_{g,\bullet,0},\BQ_{\vir})\cong\HO^{\BM}(\Mst^{\Dol}_{\bullet,0}(C),\BQ_{\vir}).
\end{equation}
\begin{proof}[Finishing the proof of Theorem \ref{ThmA}]
Using the analogues of the above results for $g(C)\leq 1$, we deduce Theorem \ref{ThmA}.  We split this up into two cases
\subsubsection{$g=0$}
As we recalled in \S \ref{cnaht}, there are isomorphisms of stacks 
\[
G\colon \Mst^{\Dol}_{r,0}(C)\cong \Mst^{\Betti}_{0,r,0}.
\]
There are isomorphisms \eqref{BMforH} allowing us to substitute Borel--Moore homology for cohomology.  Applying $\HO(-,\BQ)$ to $G$, we find the nonabelian Hodge isomorphism $\Upsilon$ required in Conjecture \ref{conj_c}.  Furthermore there are isomorphisms of varieties
\[
\Msp^{\Dol}_{r,0}(C)\cong \Lambda_r\cong \pt\cong \Msp^{\Betti}_{0,r,0}
\]
for all $r$.  It follows that the perverse filtrations $\FH^{\bullet}$ and $\FL^{\bullet}$ on $\HO(\Mst^{\Dol}_{r,0}(C),\BQ)$ both agree with the filtration by cohomological degree.  Likewise, the ``less'' perverse filtration $\FL^{\bullet}$ on $\HO(\Mst^{\Betti}_{0,r,0},\BQ)$ agrees with the filtration by cohomological degree, as does the weight filtration.  In particular the isomorphism $\Upsilon$ respects the two filtrations denoted $\FL^{\bullet}$.  The isomorphism in Conjecture \ref{conj_b} is obtained by restricting to the zeroth pieces of these filtered objects, and Conjecture \ref{conj_c} follows from the fact that all of the filtrations considered on $\HO(\Mst^{\Betti}_{0,r,0},\BQ)$ and $\HO(\Mst^{\Dol}_{r,0}(C),\BQ)$ are the same.
\subsubsection{$g=1$}
Applying the global sections functor to \eqref{relBPBW} we deduce that there is a PBW isomorphism
\begin{equation}
\label{bPBW}
\HO^{\BM}(\Mst_{1,\bullet,0},\BQ_{\vir})\cong \Sym\left( \bigoplus_{r\geq 1}\HO((\mathbb{C}^*)^2,\mathbb{Q}_{\vir})\otimes \HO(\pt/\mathbb{C}^*,\mathbb{Q})\right).
\end{equation}
Identifying $\mathbb{Q}[u]=\HO(\pt/\mathbb{C}^*,\mathbb{Q})$, we write elements of the right hand side as symmetric linear combinations of elements $\alpha_1u^{i_i}\otimes\ldots \otimes \alpha_lu^{i_l}$ with $\alpha_m\in \bigoplus_{r\geq 1}\HO((\mathbb{C}^*)^2,\mathbb{Q}_{\vir})$.  Since $\HO^i(\pt/\mathbb{C}^*,\underline{\mathbb{Q}})$ is pure of weight $i$, and $\HO^i((\mathbb{C}^*)^2,\underline{\mathbb{Q}}_{\vir})$ is pure of weight $2(i+1)$, the weight filtration admits a natural splitting: we set the weight of a tensor $\alpha_1u^{i_i}\otimes\ldots \alpha_lu^{i_l}$ to be $2\sum_{m\leq l}i_m+2\sum_{m\leq l}(\lvert \alpha_i\lvert+1)$.  From \eqref{relBPBW} again we deduce that $\FL^0\HO^{\BM}(\Mst_{1,\bullet,0},\BQ_{\vir})$ is spanned by symmetric tensors with $i_1=\ldots =i_l=0$.  

Turning to the Dolbeault moduli stack, we identify $\Msp^{\Dol}_{1,0}=\mathrm{T}^*C\cong C\times \mathbb{A}^1$, and we identify the Hitchin map $\Hit$ with the projection to $\mathbb{A}^1$.  It follows that the perverse filtration induced by the Hitchin map is split: we have $\Fh^i\HO(\Msp^{\Dol}_{1,0},\mathbb{Q}_{\vir})=\bigoplus_{j\leq i-1}\HO^j(\Msp^{\Dol}_{1,0},\mathbb{Q}_{\vir})$.  By the classical nonabelian Hodge isomorphism (or direct inspection) there is an isomorphism 
\begin{equation}
\label{soph}
\HO(\Msp^{\Dol}_{1,0},\mathbb{Q}_{\vir})\cong \HO((\mathbb{C}^*)^2,\mathbb{Q}_{\vir})
\end{equation}
sending the perverse filtration to the weight filtration, after doubling indices.  Applying the global sections functor to Proposition \ref{gen1to} there is an isomorphism
\begin{equation}
\label{dPBW}
\HO^{\BM}(\Mst^{\Dol}_{\bullet,0}(C),\mathbb{Q}_{\vir})\cong \Sym\left( \bigoplus_{r\geq 1}\HO(\Msp_{1,0}^{\Dol}(C),\mathbb{Q}_{\vir})\otimes\HO(\pt/\mathbb{C}^*,\mathbb{Q})\right).
\end{equation}
Comparing \eqref{bPBW} and \eqref{dPBW} we obtain the isomorphism of Conjecture \ref{conj_c} via the isomorphism \eqref{soph}.  Writing elements of the right hand side of \eqref{dPBW} as symmetric linear combinations of tensors $\beta_1u^{i_1}\otimes\ldots \beta_lu^{i_l}$ where $\beta_m\in\bigoplus_{r\geq 1}\HO(\Msp_{1,0}^{\Dol}(C),\mathbb{Q}_{\vir})$, we deduce from Proposition \ref{gen1to} that the two perverse filtrations on $\HO^{\BM}(\Mst^{\Dol}_{\bullet,0}(C),\mathbb{Q}_{\vir})$ admit natural splittings: we give $\beta_1u^{i_1}\otimes\ldots \beta_lu^{i_l}$ bidegree $(\sum_{m\leq l}(\lvert \beta_m\lvert+1+2i_m),2\sum_{m\leq l}i_m)$.  Then $\FL^j\HO^{\BM}(\Mst^{\Dol}_{\bullet,0}(C),\mathbb{Q}_{\vir})$ is spanned by symmetric tensors of bidegree $(a,b)$ with $b\leq j$, while $\Fh^j\HO^{\BM}(\Mst^{\Dol}_{\bullet,0}(C),\mathbb{Q}_{\vir})$ is spanned by symmetric tensors of bidegree $(a,b)$ with $a\leq j$.  In particular, $\FL^0\HO^{\BM}(\Mst^{\Dol}_{\bullet,0}(C),\mathbb{Q}_{\vir})$ is spanned by symmetric tensors with $i_1=\ldots=i_l=0$ and Conjecture \ref{conj_b} follows.  The combined filtration $\FF^{\bullet}$ is thus also split: $\FF^j\HO^{\BM}(\Mst^{\Dol}_{\bullet,0}(C),\mathbb{Q}_{\vir})$ is spanned by elements of bidegree $(a,b)$ with $a-b/2\leq j$, i.e. symmetrizations of tensors satisfying
\[
\sum_{m\leq l}(\lvert\beta_m\lvert +1+i_m)\leq j.
\]
Comparing with the description of the weight filtration on the Betti side, we deduce the PS=WS conjecture in genus 1.
\end{proof}
\subsection{PS=WS $\Leftrightarrow$ PI=WI}
\label{penusec}
\begin{theorem}
\label{fsta}
Assuming that the morphism $\HO(\Phi_{\bullet})$ from \eqref{CPD} is well-defined and an isomorphism, the PI=WI conjecture is equivalent to both Conjecture \ref{conj_b} and Conjecture \ref{conj_c}.
\end{theorem}
\begin{proof}
First we show that Conjecture \ref{conj_c} implies Conjecture \ref{conj_b}.  The proposed isomorphism $\Upsilon=\HO(\Phi_{\bullet})$ satisfies 
\[
\Upsilon\left(\mathfrak{L}^0\HO^{\BM}(\Mst^{\Betti}_{g,r,0},\BQ_{\vir})\right)=\mathfrak{L}^0\HO^{\BM}(\Mst^{\Dol}_{r,0}(C),\BQ_{\vir})
\]
by construction, so that Conjecture \ref{conj_b} follows from Conjecture \ref{conj_c} by restricting to the zeroth graded piece with respect to the less perverse filtration.  There are inclusions of mixed Hodge structures
\[
\underline{\ICA}(\Msp^{\Betti}_{g,r,0})\hookrightarrow \mathfrak{L}^0\HO^{\BM}(\Mst^{\Betti}_{g,r,0},\underline{\BQ}_{\vir})
\]
and of vector spaces respecting perverse filtrations
\[
\ICA(\Msp^{\Dol}_{r,0}(C))\hookrightarrow \mathfrak{L}^0\HO^{\BM}(\Mst^{\Dol}_{r,0}(C),\BQ_{\vir})
\]
by \cite{Da21a}.  Again by construction, $\Upsilon(\ICA(\Msp^{\Betti}_{g,r,0}))=\ICA(\Msp^{\Dol}_{r,0}(C))$, and so either of Conjectures \ref{conj_b} or \ref{conj_c} implies the PI=WI conjecture.  

Next, there is a natural isomorphism of filtered vector spaces
\begin{align}
\label{Firra}
&W_{\bullet}\left(\Sym_{\oplus}\!\left(\Free_{\Lie}\left(\bigoplus_{r\geq 1} \ICA(\Msp^{\Betti}_{g,r,0})\right)\otimes \HO(\pt/\BC^*,\BQ)\right)\right)\cong\\\nonumber
&\Sym_{\oplus}\!\left(\Free_{\Lie}\left(\bigoplus_{r\geq 1} W_{\bullet}\ICA(\Msp^{\Betti}_{g,r,0})\right)\otimes W_{\bullet}\HO(\pt/\BC^*,\BQ)\right)
\end{align}
i.e. the weight filtration on $\HO^{\BM}(\Mst^{\Betti}_{g,r,0},\underline{\BQ}_{\vir})$ is induced by the the one on $\ICA(\Msp^{\Betti}_{g,r,0})$.
\smallbreak
Let $\CF$ be a pure complex of mixed Hodge modules on $\Msp^{\Dol}_{r,0}(C)$.  Then the cohomology $\HO(\Msp^{\Dol}_{r,0}(C),\CF)$ carries two filtrations
\begin{align*}
\FH^i\HO(\Msp^{\Dol}_{r,0}(C),\CF)\coloneqq &\HO(\Lambda_r,\tau^{\leq i}\Hit_*\CF)\\
\mathfrak{L}^i\HO(\Msp^{\Dol}_{r,0}(C),\CF)\coloneqq &\HO(\Msp^{\Dol}_{r,0}(C),\tau^{\leq i}\CF)
\end{align*}
and
\begin{align*}
\FH^i\HO(\Msp^{\Dol}_{r,0}(C),\CF\otimes\LLL^k)\cong&\FH^{i-2k}\HO(\Msp^{\Dol}_{r,0}(C),\CF)\otimes\LLL^k\\
\mathfrak{L}^i\HO(\Msp^{\Dol}_{r,0}(C),\CF\otimes \LLL^k)\cong &\mathfrak{L}^{i-2k}\HO(\Msp^{\Dol}_{r,0}(C),\CF)\otimes \LLL^k.
\end{align*}
If $\CG$ and $\CF$ are pure complexes of mixed Hodge modules on $\Msp^{\Dol}_{r,0}(C)$, it follows from finiteness of the morphisms (see Proposition \ref{CupFin} for $\cup$)
\[
\oplus\colon \Msp^{\Dol}_{\bullet,0}(C)\times\Msp^{\Dol}_{\bullet,0}(C)\rightarrow \Msp^{\Dol}_{\bullet,0}(C);\quad\quad
\cup\colon\Lambda_{\bullet}\times\Lambda_{\bullet}\rightarrow \Lambda_{\bullet}
\]
that there are natural isomorphisms
\begin{align*}
\FH^{\bullet}\HO(\Msp_{\bullet,0}^{\Dol}(C),\CF\boxtimes_{\oplus} \CG)\cong &\FH^{\bullet}\HO(\Msp_{\bullet,0}^{\Dol}(C),\CF)\otimes \FH^{\bullet}\HO(\Msp_{\bullet,0}^{\Dol}(C),\CG)\\
\FL^{\bullet}\HO(\Msp_{\bullet,0}^{\Dol}(C),\CF\boxtimes_{\oplus} \CG)\cong &\FL^{\bullet}\HO(\Msp_{\bullet,0}^{\Dol}(C),\CF)\otimes \FL^{\bullet}\HO(\Msp_{\bullet,0}^{\Dol}(C),\CG).
\end{align*}
Since $\Free_{\Lie}(\CF)$ is a direct sum of summands of tensor powers of $\CF$, as is $\Sym(\CF)$, we finally deduce that there are natural isomorphisms
\begin{align}
\label{Secra}
&\FF^{\bullet}\left(\HO\left(\Msp^{\Dol}_{\bullet,0}(C),\Sym_{\oplus}\!\left(\Free_{\Lie}\left(\bigoplus_{r\geq 1} \ICS_{\Msp^{\Dol}_{r,0}(C)}\right)\otimes \HO(\pt/\BC^*,\BQ)\right)\right)\right)\cong\\\nonumber
&\Sym\!\left(\Free_{\Lie}\left(\bigoplus_{r\geq 1} \FF^{\bullet}\ICA(\Msp^{\Dol}_{r,0}(C))\right)\otimes \FF^{\bullet}\HO(\pt/\BC^*,\BQ)\right)
\end{align}
for $\FF^{\bullet}=\FH^{\bullet}$ or $\FF^{\bullet}=\FL^{\bullet}$, where for both filtrations we take the cohomological filtration on $\HO(\pt/\BC^*,\BQ)$, and the filtration $\FL^{\bullet}$ on $\ICA(\Msp^{\Dol}_{r,0}(C))$ is the trivial filtration for which the associated graded object is concentrated entirely in degree zero.  By the same reasoning, we deduce that
\begin{align}
\label{Secrt}
&\FL^{\bullet}\left(\HO\left(\Msp^{\Dol}_{\bullet,0}(C),\Sym_{\oplus}\!\left(\Free_{\Lie}\left(\bigoplus_{r\geq 1} \ICS_{\Msp^{\Betti}_{r,0}(C)}\right)\otimes \HO(\pt/\BC^*,\BQ)\right)\right)\right)\cong\\\nonumber
&\Sym\!\left(\Free_{\Lie}\left(\bigoplus_{r\geq 1} \FL^{\bullet}\ICA(\Msp^{\Betti}_{r,0}(C))\right)\otimes \FL^{\bullet}\HO(\pt/\BC^*,\BQ)\right)
\end{align}
With $\FL^{\bullet}$ again denoting the trivial filtration on $\ICA(\Msp^{\Betti}_{r,0}(C))$.  It follows from \eqref{Firra}, \eqref{Secra} and \eqref{Secrt} that Conjecture \ref{IPIW} implies Conjecture \ref{conj_c}.

\end{proof}

\subsection{$\chi$-independence $\Rightarrow$ (PS=WS $\Leftrightarrow$ P=W)}
We finish by using the same constructions to give a conjectural equivalence between Conjectures \ref{conj_c} and the original P=W conjecture: i.e. Theorem \ref{PW}.  We show that Conjectures \ref{PW} and \ref{conj_b} are equivalent, and then use the equivalence of Conjectures \ref{conj_c} and \ref{conj_b} from the previous subsection.  As in \S \ref{penusec} we will assume that the proposed morphism
\[
\Phi_{\bullet}\colon\Gamma_{\bullet,0,*}p_*\BQ_{\Mst^{\Betti}_{g,\bullet,0},\vir}\rightarrow p_*\BQ_{\Mst^{\Dol}_{\bullet,0}(C),\vir}.
\]
is an isomorphism.  Then we can write
\[
p_*\BQ_{\Mst^{\Dol}_{\bullet,0}(C),\vir}\cong \Sym_{\oplus}\left(\bigoplus_{r\geq 1}\Gamma_{r,0,*}\DTS^{\Betti}_{g,r}\otimes \HO(\pt/\BC^*,\BQ)\right)
\]
and by Remark \ref{DolUn} we may as well \textit{define}
\[
\DTS^{\Dol}_r(C)\coloneqq \Gamma_{r,0,*}\DTS^{\Betti}_{g,r}.
\]
Now the (Dolbeault) $\chi$-independence theorem (Theorem \ref{HiggsChi}) states that there is an isomorphism
\[
\Hit_{*}\DTS^{\Dol}_r(C)\cong \Hit_*\BQ_{\Msp^{\Dol}_{r,1}(C),\vir}
\]
so that there is an isomorphism of filtered vector spaces
\[
\FH^{\bullet}\HO(\Msp_{r,0}^{\Dol}(C),\DTS^{\Dol}_r(C))\cong \FH^{\bullet}\HO(\Msp_{r,1}^{\Dol}(C),\BQ_{\vir}).
\]
On the other hand, the Betti $\chi$-independence conjecture (Conjecture \ref{BCI}) yields the isomorphism of filtered vector spaces
\[
W_{\bullet}\HO(\Msp^{\Betti}_{g,r,0},\DTS^{\Betti}_{g,r})\cong W_{\bullet}\HO(\Msp^{\Betti}_{g,r,1},\ul{\BQ}_{\vir}).
\]
Arguing as in \S \ref{penusec}, via finiteness of $\cup\colon \Lambda_{\bullet}\times \Lambda_{\bullet}\rightarrow \Lambda_{\bullet}$ and $\oplus \colon\Msp^{\Dol}_{\bullet,0}(C)\times \Msp^{\Dol}_{\bullet,0}(C)\rightarrow \Msp^{\Dol}_{\bullet,0}(C)$, there is an equality of filtrations
\begin{align*}
&\FH^{\bullet}\left(\HO(\Msp^{\Dol}_{\bullet,0}(C), \Sym_{\oplus}\left(\DTS^{\Dol}_{\bullet}(C)\right)\right)\cong \Sym\left(\FH^{\bullet}\HO\left(\Msp^{\Dol}_{\bullet,0}(C),\DTS^{\Dol}_{\bullet}(C)\right)\right),
\end{align*}
as well as the usual equality of filtrations
\begin{align*}
&W_{\bullet}\left(\HO(\Msp^{\Dol}_{\bullet,0}(C), \Sym_{\oplus}\left(\DTS^{\Betti}_{g,\bullet}\otimes \HO(\pt/\BC^*,\BQ)\right)\right)\cong \\&\Sym\left(W_{\bullet}\HO\left(\Msp^{\Betti}_{g,\bullet},\DTS^{\Dol}_{\bullet}(C)\right)\otimes W_{\bullet}\HO(\pt/\BC^*,\BQ)\right).
\end{align*}
Moreover, by construction, the proposed isomorphism 
\[
\FL^0\Upsilon=\FL^0\HO(\Gamma_{\bullet,0})\colon \FL^0\HO^{\BM}(\Mst^{\Betti}_{g,\bullet,0},\BQ_{\vir})\xrightarrow{\cong}\FL^0\HO^{\BM}(\Mst^{\Dol}_{\bullet,0}(C),\BQ_{\vir})
\]
satisfies
\[
\Upsilon(\HO\left(\Msp^{\Betti}_{g,\bullet},\DTS^{\Dol}_{\bullet}(C)\right))=\HO\left(\Msp^{\Dol}_{\bullet,0}(C),\DTS^{\Betti}_{g,\bullet}\right).
\]
Still arguing as in \S \ref{penusec}, the isomorphism $\Upsilon$ satisfies
\[
\Upsilon(W_{2i}\FL^0\HO^{\BM}(\Mst^{\Betti}_{g,r,0},\BQ_{\vir}))=\FH^i\FL^0\HO^{\BM}(\Mst^{\Dol}_{r,0}(C),\BQ_{\vir})
\]
if and only if
\[
\Upsilon(W_{2i}\HO(\Msp^{\Betti}_{g,r,1},\BQ_{\vir}))=\FH^i\HO(\Msp^{\Dol}_{r,1}(C),\BQ_{\vir}).
\]
This finishes the proof of Theorem \ref{thm_b}.

\appendix
\section{Zero-dimensional sheaves on surfaces}
\label{DZ_sec}
In this appendix we give a relatively quick proof of Proposition \ref{KVredux}, which is the modification of one of the main results of \cite{KV19} that we need in the paper.
\smallbreak
For $X$ an algebraic variety we denote by
\[
\Delta_n\colon X\rightarrow \Sym^n(X)
\]
the inclusion of the small diagonal.  
\smallbreak
We begin with a couple of preliminaries.
\begin{proposition}\cite{preproj,QEAs}
\label{PBWA2}
There is a PBW isomorphism
\begin{equation}
\label{PBee}
p_{\bullet,*}\BD\ul{\BQ}_{\Mst_{\bullet}(\BA^2)}\cong\Sym_{\cup}\left(\bigoplus_{m\geq 1}\Delta_{m,*}\ul{\BQ}_{\BA^2,\vir}\otimes\HO(\pt/\BC^*,\BQ)\right).
\end{equation}
\end{proposition}
\begin{proposition}\cite{Bar01}\label{BaraProp}
Let $\CN_n\subset \FM_n(\BA^2)$ denote the stack of pairs of commuting nilpotent matrices.  Then $\CN_n$ is irreducible and $-1$-dimensional.
\end{proposition}
We give a quick alternative proof of Proposition \ref{BaraProp}, using Proposition \ref{PBWA2}.
\begin{proof}
Let $\iota_n\colon0_n\hookrightarrow \Sym^n(\BA^2)$ denote the inclusion of the origin.  Then by base change, we have an isomorphism
\begin{align*}
\HO^{\BM}(\CN_n,\BQ)\cong &\iota_n^! p_{n,*}\BD\BQ_{\Mst_n(\BA^2)}\\
\cong &\iota_n^!\Sym_{\cup}\left(\bigoplus_{m\geq 1}\Delta_{m,*}\BQ_{\BA^2,\vir}\otimes\HO(\pt/\BC^*,\BQ)\right)\\
\cong &\Sym_{\cup}\left(\bigoplus_{m\geq 1}\iota_m^!\Delta_{m,*}\BQ_{\BA^2,\vir}\otimes\HO(\pt/\BC^*,\BQ)\right)\lvert_{0_n}\\
\cong &\Sym_{\cup}\left(\bigoplus_{m\geq 1}\BQ_{0_m}[-2]\otimes\HO(\pt/\BC^*,\BQ)\right)\lvert_{0_n}\\
\cong &\Sym_{\cup}\left(\bigoplus_{m\geq 1}\bigoplus_{i\geq 1}\BQ_{0_m}[-2i]\right)\lvert_{0_n}
\end{align*}
where the second isomorphism follows from Proposition \ref{PBWA2}.  Since the compactly supported cohomology of $\CN_n$ is the graded vector dual of the Borel--Moore homology, we deduce that the largest $i$ such that $\HO_c^i(\CN_n,\BQ)\neq 0$ is $-2$, and $\HO_c^{-2}(\CN_n,\BQ)\cong\BQ$.
\end{proof}
For $X$ a surface, let $\mathcal{N}_n(X)\subset \Mst_{n}(X)$ denote the substack of the stack of coherent sheaves with punctual support, which we \textit{define}\footnote{I.e. we are ignoring the fact that the object representing the moduli functor of coherent sheaves with set-theoretically punctual support is actually a formal completion of the stack we define here.} via the Cartesian diagram
\[
\xymatrix{
\mathcal{N}_n(X)\ar[r]^-{\tilde{\Delta}_n}\ar[d]^{q_n}&\Mst_{n}(X)\ar[d]^{p_n}
\\
X\ar[r]^-{\Delta_n}&\Sym^n(X).
}
\]
The morphism $\tilde{\Delta}_n$ is a closed embedding, since $\Delta_n$ is.  
\begin{proposition}
\label{findBPS}
The natural morphism $\Up\tau^{\leq 0}\left(\Delta_{n,*}q_{n,*}\BD\BQ_{\CN_n(\BA^2)}\xrightarrow{\alpha} p_{n,*}\BD\BQ_{\Mst_m(\BA^2)}\right)$ is a split inclusion of perverse sheaves.
\end{proposition}
\begin{proof}
First we construct the morphism $\alpha$.  Applying the functor $p_{n,!}$ to the natural morphism $\BQ_{\Mst_n(\BA^2)}\rightarrow \tilde{\Delta}_{n,*}\BQ_{\CN_n(\BA^2)}$, and applying base change, we obtain the morphism
\[
p_{n,!}\BQ_{\Mst_n(\BA^2)}\rightarrow \Delta_{n,*}q_{n,!}\BQ_{\CN_n(\BA^2)}
\]
from which we obtain $\alpha$ by applying $\BD$.  I.e. we obtain $\alpha$ by applying the natural transformation $\id\rightarrow \Delta_{n,*}\Delta_n^*$ to $p_{n,!}\BQ_{\Mst_n(\BA^2)}$ and then taking the Verdier dual.  So to prove the proposition it is enough to show that 
\begin{equation}
\label{APSh}
\Up\tau^{\geq 0}\left(\id\rightarrow \Delta_{n,*}\Delta_n^*\right)p_{n,!}\BQ_{\Mst_n(\BA^2)}
\end{equation}
is a split surjection of perverse sheaves.  
\smallbreak
Firstly we show that the domain and target of \eqref{APSh} are actually perverse sheaves.  The domain is a perverse sheaf, since, taking Verdier duals of \eqref{PBee}, $\Up\Ho^i(p_{n,!}\BQ_{\Mst_{n}(\BA^2)})=0$ for $i>0$.  We write 
\[
q_{n,!}\BQ_{\CN_n(\BA^2)}\cong \BQ_{\BA^2}\otimes\HO_c(\CN_n,\BQ)
\]
and so $\Up\tau^{\geq 0}q_{n,!}\BQ_{\CN_n(\BA^2)}\cong \BQ_{\BA^2,\vir}$ is a perverse sheaf by Proposition \ref{BaraProp}.  Since $\Delta_{n,*}$ is exact with respect to the perverse t-structure, it follows that 
\[
\Up\tau^{\geq 0}\Delta_{n,*}\Delta_n^*p_{n,!}\BQ_{\Mst_n(\BA^2)}\cong \Up\tau^{\geq 0}\Delta_{n,*}q_{n,!}\BQ_{\CN_n(\BA^2)}\cong\Delta_{n,*}\mathbb{Q}_{\mathbb{A}^2,\vir}
\]
is indeed perverse.  
\smallbreak
We write
\begin{align*}
p_{n,!}\BQ_{\Mst_{\bullet}(\BA^2)}\cong &\Sym_{\cup}\left(\bigoplus_{m\geq 1}\Delta_{m,*}\BQ_{\BA^2,\vir}\otimes\HO(\pt/\BC^*,\BQ)^{\vee} \right)\Lvert_{\Sym^n(\BA^2)}\\
\cong & \left(\Delta_{n,*}\BQ_{\BA^2,\vir}\otimes\HO(\pt/\BC^*,\BQ)^{\vee}\right)\oplus \\ &\oplus\Sym^{\geq 2}_{\cup}\left(\bigoplus_{m\geq 1}\Delta_{m,*}\BQ_{\BA^2,\vir}\otimes\HO(\pt/\BC^*,\BQ)^{\vee} \right)\Lvert_{\Sym^n(\BA^2)}
\end{align*}
We have shown that the functor $\Up\tau^{\geq 0}\Delta_{n,*}\Delta_n^*p_{n,!}$ applied to the first summand, is isomorphic to $\Up\tau^{\geq 0}\Delta_{n,*}\Delta_n^*p_{n,!}$ applied to the LHS, so that $\Up\tau^{\geq 0}\Delta_{n,*}\Delta_n^*p_{n,!}$ applied to the second summand is zero.  Applying the functor $\id\rightarrow \Delta_{n,*}\Delta_n^*$ to the first summand yields an isomorphism of perverse sheaves, and so the morphism $\Up\tau^{\geq 0}(\id\rightarrow \Delta_{n,*}\Delta_n^*)p_{n,!}\BQ_{\Mst_{\bullet}(\BA^2)}$ is a split surjection of perverse sheaves, as required.
\end{proof}

\begin{proposition}
\label{KVredux}
Let $X$ be a smooth surface.  Let $\Mst_{n}(X)$ denote the stack of coherent sheaves on $X$ with zero-dimensional support and length $n$, and let 
\[
p_n\colon \Mst_{n}(X)\rightarrow \Sym^n(X)
\]
denote the morphism taking a coherent sheaf to its support.  Then there is a PBW isomorphism of cohomologically graded mixed Hodge modules
\begin{equation}
\label{VT1}
p_*\BD\underline{\BQ}_{\Mst_{\bullet}(X)}\cong \Sym_{\cup}\left(\bigoplus_{n\geq 1}\Delta_{n,*}\underline{\BQ}_{X,\vir}\otimes \HO(\pt/\BC^*,\BQ)\right).
\end{equation}
Moreover, there is an isomorphism of algebra objects in the tensor category $\MHM(\Sym^{\bullet}(X))$ with the convolution product provided by $\cup\colon \Sym^{\bullet}(X)\times\Sym^{\bullet}(X)\rightarrow \Sym^{\bullet}(X)$
\begin{equation}
\label{VT2}
\tau^{\leq 0}p_*\BD\underline{\BQ}_{\Mst_{\bullet}(X)}\cong \Sym_{\cup}\left(\bigoplus_{n\geq 1}\Delta_{n,*}\underline{\BQ}_{X,\vir}\right).
\end{equation}
\end{proposition}
For the proof we give a simplified version of the proof of a related result (regarding hypercohomology) in \cite{KV19}.  We note that the isomorphism \eqref{VT1} is typically never realised via an isomorphism of algebra objects.  Indeed for the case $X=\mathbb{A}^2$ it is shown in \cite{Da22a} that there is an isomorphism of Hall algebras
\[
\bigoplus_{d\geq 0}\HO^{\BM}(\Mst_d(\mathbb{A}^2),\mathbb{Q})\cong \UEA(L)
\]
where $L$ is a $\mathbb{Z}_{\geq 1}$-graded, cohomologically graded Lie algebra with
\[
L_n^i=\begin{cases} \mathbb{Q}\cdot \alpha_{n}^{(i/2+1)}&\textrm{if }n\geq 1\textrm{ and } i\in 2\cdot \mathbb{Z}_{\geq -1}\\
0& \textrm{otherwise.}
\end{cases}
\]
and $[\alpha^{(i)}_n,\alpha_m^{(j)}]=(im-jn)\alpha_{m+n}^{(i+j-1)}$.  We set $\alpha_n^{(-1)}=0$ for all $n$.  In particular, while the subalgebra generated by $\alpha_n^{(0)}$ for $n\geq 1$ is commutative, and corresponds to taking the derived global sections of either side of \eqref{VT2}, the entire algebra is not.
\begin{proof}[Proof of Proposition \ref{KVredux}]
We start by constructing \eqref{VT1}.  By Proposition \ref{BaraProp} the stack $\mathcal{N}_n(X)$ is an irreducible $1$-dimensional stack.  Writing $X$ analytically locally as an open ball $U$ in $\BA^2$, we find that (locally) $\mathcal{N}_n(X)$ can be written as $U\times \mathcal{N}_n$.  By Proposition \ref{BaraProp}, we have $\HO_c^{-2}(\CN_n,\BQ)\cong \BQ$ and $\HO_c^{i}(\CN_n,\BQ)=0$ for $i>-2$.
\smallbreak
It follows that 
\begin{align*}
\tau^{\leq 0}q_{n,*}\BD\ul{\BQ}_{\mathcal{N}_n(X)}\cong &\BD\tau^{\geq 0}q_{n,!}\ul{\BQ}_{\mathcal{N}_n(X)}\\
\cong &\BD(\tau^{\geq 0}\left(\ul{\BQ}_X\otimes \HO_c(\CN_n,\ul{\BQ})\right))\\
\cong &\BD(\ul{\BQ}_{X}\otimes\LLL^{-1})\\
\cong &\ul{\BQ}_{X}\otimes\LLL^{-1}.
\end{align*}
Taking the Verdier dual of the morphism
\[
\ul{\BQ}_{\Mst_{n}(X)}\rightarrow \tilde{\Delta}_{n,*}\ul{\BQ}_{\mathcal{N}_n(X)}
\]
we obtain the morphism 
\[
\beta\colon \tilde{\Delta}_{n,*}\BD\ul{\BQ}_{\mathcal{N}_n(X)}\rightarrow \BD\ul{\BQ}_{\Mst_{n}(X)}.
\]
Applying ${}^{\Fp}\!\tau^{\leq 0}p_{n,*}$ to $\beta$ yields the morphism
\[
\Delta_{n,*}\ul{\BQ}_{X}\otimes\LLL^{-1}\rightarrow \tau^{\leq 0}p_{n,*}\BD\underline{\BQ}_{\Mst_{n}(X)}
\]
which we compose with the natural transformation
\[
\tau^{\leq 0}p_{n,*}\BD\BQ_{\Mst_{n}(X)}\rightarrow p_{n,*}\BD\BQ_{\Mst_{n}(X)}
\]
to obtain the morphism
\[
\Phi_{X,n}\colon\Delta_{n,*}\ul{\BQ}_{X}\otimes\LLL^{-1}\rightarrow \ul{\RCoha}_{X,n}
\]
where the target is the $n$th piece of $\ul{\RCoha}_{X,\bullet}\coloneqq p_{\bullet,*}\BD\underline{\BQ}_{\Mst_{\bullet}(X)}$, the relative Hall algebra of zero-dimensional coherent sheaves on $X$.  Then via the multiplication and $\HO(\pt/\BC^*,\BQ)$-action on $\ul{\RCoha}_{X,\bullet}$, we obtain a morphism
\begin{equation}
\label{PBC2}
\Psi_X\colon\Sym_{\cup}\left(\bigoplus_{n\geq 1}\Delta_{n,*}\ul{\BQ}_X\otimes\LLL^{-1}\otimes \HO(\pt/\BC^*,\ul{\BQ})\right)\rightarrow \ul{\RCoha}_{X,\bullet}
\end{equation}
which we claim is an isomorphism.  To see this, it is enough to work analytically locally on $\Sym^n(X)$, and check the statement at the level of the underlying cohomologically graded perverse sheaves.  The variety $\Sym^n(X)$ admits an analytic cover by subspaces $\Sym^n(\coprod_{j\leq r} U_j)$, where $U_1,\ldots,U_r$ are disjoint open balls in $X$, and we assume we are given embeddings $U_i\subset \BA^2$.  Via the natural isomorphism
\[
\Sym^{\bullet}\left(\coprod_{j\leq r}U_j\right)\cong \prod_{j\leq r}\Sym^{\bullet}(U_j)
\]
we can reduce the claim that \eqref{PBC2} is an isomorphism to the special case in which $X=U_i$.  This case in turn follows from the special case in which $X=\BA^2$, by restriction.  By Corollary \ref{findBPS} each of the morphisms $\Phi_{\BA^2,n}$ is a split inclusion, and the inclusion of the BPS sheaf, and the claim regarding $\Psi_{\BA^2}$ follows from the PBW theorem (Proposition \ref{PBWA2}).  We have thus produced the isomorphism \ref{VT1}.
\smallbreak
Applying $\tau^{\leq 0}$ to \eqref{PBC2} produces a PBW isomorphism
\[
\Sym_{\cup}\left(\bigoplus_{n\geq 1}\Delta_{n,*}\underline{\BQ}_{X,\vir}\right)\cong \tau^{\leq 0}p_*\BD\underline{\BQ}_{\Mst_{\bullet}(X)}
\]
and the claim that this is in fact an isomorphism of algebras reduces to the claim that $\tau^{\leq 0}p_*\BD\ul{\BQ}_{\Mst_{\bullet}(X)}$ is a commutative algebra object.  Again, we may work analytically locally and reduce to the claim that $\Up\tau^{\leq 0}p_*\BD\BQ_{\Mst_{\bullet}(\BA^2)}$ is a commutative algebra object.  With notation as in \S \ref{propsec}, this is equivalent to the claim that $\DTS_{\Pi_Q,\bullet}$ is an Abelian Lie algebra object, for $Q$ the Jordan quiver with one vertex and one loop.  For this the argument is given in \S \ref{BG1}: the Lie algebra morphism is induced by a morphism
\begin{equation}
\label{finlo}
\cup_*(\Delta_{m,*}\underline{\BQ}_{\mathbb{A}^2,\vir}\boxtimes\Delta_{n,*}\underline{\BQ}_{\mathbb{A}^2,\vir})\rightarrow \Delta_{m+n,*}\underline{\BQ}_{\mathbb{A}^2,\vir}
\end{equation}
where the source and the target are semisimple mixed Hodge modules with differing supports, so that \eqref{finlo} is indeed zero.
\end{proof}

\bibliographystyle{alpha}
\bibliography{Literatur}

\end{document}